%% file: main.tex
\documentclass{article}

\usepackage{amsmath,amsfonts,amsthm,epsfig,latexsym,graphicx,amssymb}

\usepackage[english]{babel}

\usepackage{xfrac}

\usepackage[dvipsnames]{xcolor}
\usepackage{mathtools}
\usepackage{graphpap} 

\usepackage{version}
\usepackage[toc,page]{appendix} 

\usepackage{subfiles} 

\usepackage[
  colorlinks = true,  
  linkcolor = PineGreen, 
  citecolor = PineGreen 
  ]{hyperref} 

\usepackage{sectsty} 
\sectionfont{\color{Salmon!70!red}}
\subsectionfont{\color{Salmon!70!red}}

\DeclareTextFontCommand{\emph}{\color{ForestGreen}\em} 

\usepackage{caption}
\DeclareCaptionFont{myblue}{\color{RoyalBlue}}
\definecolor{theocolor}{rgb}{0, 0.5, 0}

\newcommand{\NN}{\mathbb{N}}
\newcommand{\RR}{\mathbb{R}}
\newcommand{\ZZ}{\mathbb{Z}}
\newcommand{\Class}{\mathcal{C}}
\newcommand{\TT}{\mathbb{T}}
\newcommand{\Intt}[1]{\mathring{#1}}
\newcommand{\algcap}{\mathrlap{\hspace{2.3pt}\cdot}{\cap}}
\newcommand{\wb}{\overline}

\newtheoremstyle{colorplain}%
{\topsep}   
{\topsep}   
{\itshape}  
{0pt}       
{} 
{.}         
{5pt plus 1pt minus 1pt} 
{\textbf{\textcolor{RoyalBlue}{\textbf{\thmname{#1} \thmnumber{#2}}}}\thmnote{ (#3)}}
{}

\newtheoremstyle{colorremark}%
{\topsep}   
{\topsep}   
{\it}  
{0pt}       
{\itshape} 
{.}         
{5pt plus 1pt minus 1pt} 
{\textcolor{Purple}{\thmname{#1} \thmnumber{#2}}\thmnote{ (#3)}}
{}

\newtheoremstyle{colordefinition}%
{\topsep}   
{\topsep}   
{}  
{0pt}       
{} 
{.}         
{5pt plus 1pt minus 1pt} 
{\textcolor{Purple}{\textbf{\thmname{#1} \thmnumber{#2}}}\thmnote{ (#3)}}
{}

\theoremstyle{colorplain}

\newtheorem{maintheorem}{Theorem}
\newtheorem{maincorollary}[maintheorem]{Corollary}

\newtheorem{theorem}{Theorem}
\newtheorem{lemma}[theorem]{Lemma}
\newtheorem{proposition}[theorem]{Proposition}

\newtheorem{corollary}[theorem]{Corollary}

\newtheorem*{claim*}{Claim}

\theoremstyle{colorremark}

\newtheorem{remark}[theorem]{Remark}

\theoremstyle{colordefinition}
\newtheorem{definition}[theorem]{Definition}

\newcommand{\Sec}{\Section}

\def\ra{{\rightarrow}}

\def\del{\partial}
\def\Sec{\Section}

\def\hsp{{\hspace{3mm}}}

\DeclareMathOperator{\Int}{Int}

\DeclareMathOperator{\Fix}{Fix}

\newcommand\mcal[1]{{\mathcal{#1}}}
\newcommand\wt[1]{{\widetilde{#1}}}
\newcommand\wtc[1]{{\wt{\mcal{#1}}}}
\newcommand\wh[1]{{\widehat{#1}}}

  \def\cH{{\mathcal H}}   
    
   \def\cP{{\mathcal P}} \def\cV{{\mathcal V}}

\def\Si{\Sigma}
\def\Ga{\Gamma}
\def\De{\Delta}

\title{Oriented Birkhoff sections of Anosov flows}
\author{Masayuki Asaoka, Christian Bonatti, Th\'eo Marty }
\date{}

\newcommand{\flow}{{\phi}}
\newcommand{\returnmap}{{\Phi}}
\newcommand{\orbitspace}{{\mathcal{P}_\phi}}
\newcommand{\orbitspaceG}{\mathcal{P}_{\phi,\Gamma}}
\newcommand{\orbitspaceS}{\mathcal{P}_{\phi,\del \Section}}

\newcommand{\Section}{S}
\newcommand{\surectangle}{\mathcal{R}}
\newcommand{\link}{link}
\newcommand{\mult}{mult}
\newcommand{\ofoliation}{{\mathcal{L}}}
\newcommand{\foliation}{{\mathcal{F}}}
\newcommand{\lfoliation}{{\wt{\mathcal{F}}}}

\begin{document}

\maketitle

\begin{abstract}
This paper gives $3$ different proofs (independently obtained by the $3$ authors) of the following fact: 
given an Anosov flow on an oriented $3$ manifold, the existence of a positive Birkhoff section is equivalent to the fact that the flow is $\RR$-covered positively twisted. 
\end{abstract}

{
    \hypersetup{linkcolor=black} 
    \tableofcontents \label{ToC}
}

\section*{Introduction}
\addcontentsline{toc}{section}{Introduction}

\subfile{Section/Introduction}

%
%

\paragraph{Acknowledgments.}  
MA was supported by the JSPS Kakenki Grants 18K03276.
TM is grateful to P.Dehornoy and to the Max Plank Institute in Bonn.
\bigskip

\subfile{Section/Preliminar}  

\section{Oriented partial sections}
\label{Section:Marty}

\subfile{SectionMarty/SectionMarty}

\section{Drift along Birkhoff sections}
\label{Section:Asaoka}

\subfile{SectionAsaoka/SectionAsaoka}

\section{Holonomy in the bi-foliation of pseudo-Anosov map.}
\label{Section:Bonatti}

\subfile{SectionBonatti/SectionBonatti}




\input{biblio-unified}
\addcontentsline{toc}{section}{References}

\end{document}

%% file: Section/Introduction.tex


The geodesic flow of hyperbolic surfaces has been the object of extensive studies since the beginning of the 20th century, in particular with  works of Hadamard, Birkhoff and Hopf. In the 1950s and 1960s, in order to extend the results to the geodesic flow of surfaces of negative but variable curvature, Anosov extracted an essential structure, which he called the U-systems. In the language of Smale’s theory, the notion of U-systems has become what is now called \emph{hyperbolicity}, and vector fields with hyperbolic dynamics on all the ambient manifold 
are now called \emph{Anosov flows}. One of the main features of Anosov flows is that they are structurally stable: any vector field close, in the $C^1$-topology, to an Anosov flow $X$ is Anosov itself and is orbitally (topologically) equivalent to $X$.  In this paper we consider Anosov flows up to orbital equivalence.

In addition to the geodesic flows in negative curvature, the classical examples, known as \emph{algebraic}, include suspension of hyperbolic automorphisms of the torus $\mathbb{T}^n$, or more generally (in higher dimensions)  quotients of a Lie algebra by a lattice.

There are many  examples of non-algebraic Anosov flows on $3$-manifolds starting in 1980 with \cite{FrWi} and \cite{HT} who presented $2$ different processes for building Anosov flows on $3$-manifolds. These two processes have been made more powerful and flexible, making it possible to build new flows either by surgery from already built flows (see \cite{Go,Fri,Sh}), or as  construction game by gluing hyperbolic plugs along their boundary (see \cite{BeBoYu}). 

One of the main tools for understanding the dynamics of an Anosov flow on a $3$-manifold $M$  are the \emph{Birkhoff sections} $\Sigma$: this means that $\Sigma$ is a compact surface whose interior is embedded in $M$ transversely to $X$ and whose boundary $\partial \Sigma$ consists of finitely many periodic orbits of $X$; furthermore, $\Sigma$ meets every orbit in a bounded time. In 1983, Fried \cite{Fri} proves that any transitive Anosov flow on an oriented $3$-manifold admits a Birkhoff section denoted by $\Sigma$, and the first return map of the flow on $\Sigma$ induces a pseudo-Anosov homeomorphism $h$ on the closed surface $\tilde \Sigma$ obtained from $\Sigma$ by collapsing every component of $\partial \Sigma$. Furthermore, the flow $X$ can be recovered from the suspension of $h$ by simple surgeries which depends on two numbers, the \emph{multiplicity} and the \emph{linking number},  associated to every component of $\partial \Sigma$ and which describe the topological position of $\Sigma$ in a small neighborhood of that component :
\begin{itemize}
\item the connected component $c$ of $\Sigma$ (oriented by the boundary orientation) covers the periodic orbit $\gamma$ and the multiplicity $m(\Sigma,\gamma)$ is the topological degree of this covering.
 \item Now the stable (or equivalently the unstable) manifold of $\gamma$ defines a notion of parallel along $\gamma$ and the linking number (see Section~\ref{subsection:LocalInvariant}) described how $\Sigma$ turns around $\gamma$, with respect to this notion of parallel.
\end{itemize}

Another main tool is \emph{the bi-foliated plane} associated to the Anosov flow $X$:
\cite{Ba1,Fe1} show that, given an Anosov flow $X$ on a $3$-manifold $M$, the lift $\tilde X$ of $X$ on the universal cover $\tilde M$ is conjugated to the translation flow $\frac{\partial}{\partial x}$ on~$\RR^3$. Thus the orbit space
$\cP_X$ of $\tilde X$ is a plane (diffeomorphic to $\RR^2$).  The lifts of the weak stable and unstable ($2$-dimensional) foliations of $X$ 
project on $\cP_X$ in two transverse, $1$-dimensional foliations, $\ofoliation^s$ and $\ofoliation^u$, respectively. The triple $(\cP_X,\ofoliation^s,\ofoliation^u)$ is the bi-foliated plane associate to $X$.

A foliation on $\RR^2$ is trivial, that is conjugated to the horizontal foliation $\cH$ directed by $\frac{\partial}{\partial x}$, if and only if its leaf space  is Hausdorff. Therefore it is quite natural to ask if the leaves spaces of $\ofoliation^s$ or $\ofoliation^u$ are Hausdorff. The answer to these questions leads to the main splitting of the set of (equivalence classes of) Anosov flows in $3$ subsets with distinct behaviors.
Indeed, \cite{Ba1,Fe1} show that the leaves space of $\ofoliation^s$ is Hausdorff if and only if the leaves space of $\ofoliation^u$ is Hausdorff. Such Anosov flows are called \emph{$\RR$-covered} (or else \emph{alignable} in the french terminology of Barbot).
Furthermore, in this $\RR$-covered case, there are only two possibilities for the bi-foliated plane $(\cP_X,\ofoliation^s,\ofoliation^u)$:
\begin{itemize}
 \item either $(\cP_X,\ofoliation^s,\ofoliation^u)$ is conjugated to $(\RR^2, \cH,\cV)$ where $\cH$ and $\cV$  are  directed by $\frac{\partial}{\partial x}$ and $\frac{\partial}{\partial y}$ respectively; in that case Solodov proved that $X$ is topologically equivalent to the suspension of a hyperbolic linear automorphism of $\TT^2$,

\item or  $(\cP_X,\ofoliation^s,\ofoliation^u)$ is conjugated to the restriction of  $(\cH,\cV)$ to the strip $\{(x,y)\in\RR^2, |x-y|<1\}$.  In that case one says that $X$ is \emph{twisted} (or equivalently \emph{skewed}) \emph{$\RR$-covered}.
\end{itemize}

By fixing an orientation on the ambient manifold, and thus on the orbit space, one can distinguish positively twisted and negatively twisted $\RR$-covered Anosov flows.
Our main result exhibit a relation between  Birkhoff sections and bi-foliated plane:

\begin{maintheorem}\label{theorem:rcoveredcondition}
    Let~$M$ be an oriented closed~$3$-manifold and~$\phi$ be an Anosov flow on~$M$. The flow~$\phi$ admits a positive Birkhoff section if and only if it is positively twisted~$\RR$-covered.
\end{maintheorem}

Theorem~\ref{theorem:rcoveredcondition} is an important step in \cite{Ma} towards proving that every twisted $\RR$-covered Anosov flow is orbitally equivalent to a Reeb-Anosov flow.

\begin{maincorollary}
\label{Acor:cor B}
Let $\flow$ be a topologically transitive Anosov flow on an orientable 3-manifold. It can be obtained from an $\RR$-covered Anosov flow by a Goodman-Fried surgery along finitely many periodic orbits. If $S$ is a Birkhoff section of $\flow$, the Goodman-surgeries in the previous statement can be chosen on the boundary orbits of $S$.
\end{maincorollary}

Corollary \ref{Acor:cor B} is reminiscent of one from Bonatti and Iakovoglou~\cite{BoIa}.
As a special case of Theorem \ref{theorem:rcoveredcondition}, we also obtain
\begin{maincorollary}
\label{Acor:cor D}
Let $\flow$ be a topologically transitive Anosov flow on an orientable 3-manifold. If $\flow$ admits a Birkhoff section with just one boundary component, then $\flow$ is a twisted $\RR$-covered flow.
\end{maincorollary}


One approach to prove the main theorem is to generalize the classification of $\RR$-covered flows by the existence of some Birkhoff sections. Non-necessarily Anosov flows on an oriented closed manifold of dimension~3 are classified in four families, which correspond to flows admitting a positive Birkhoff section, a global section, a negative Birkhoff section, and the other flows. We call the flows satisfying the first three properties respectively positively twisted flows, flat flows and negatively twisted flows. For Anosov flows, they correspond respectively to the three cases: positively twisted $\RR$-covered, trivially $\RR$-covered (or suspensions), and negatively twisted $\RR$-covered.

\begin{maintheorem}\label{Ttheorem:BirkhoffSectionRestriction}
    Let~$M$ be an oriented closed 3-dimensional manifold and $\flow$ a smooth flow on $M$. If $\flow$ admits a positive Birkhoff section, then it does not admit any negative partial section.

    If $\flow$ admits a Birkhoff section without boundary (a global section), then~$\flow$ does not admit any  partial section with only positive or only negative boundary components.
\end{maintheorem}

\begin{maincorollary}\label{Tcorollary:NatureOfAFlow}
    For a flow on an oriented closed 3-dimensional manifold, being flat, positively twisted or negatively twisted are mutually exclusive cases.
\end{maincorollary}

The mutual exclusivity of oriented partial sections, combined with a previous work of Barbot who defined a family of oriented partial sections, gives one proof of the main theorem.

\cite{BoIa} shows that surgeries along the finite set of pivot orbits of a non-$\RR$-covered Anosov flow does not change the fact that the flow is not $\RR$-covered. In the same spirit, \cite{BoIa} shows that surgeries on the geodesic flow of a hyperbolic surface along the orbits corresponding to simple disjoint geodesics does not change the property of being $\RR$-covered negatively twisted. Thus one deduces that some finite set of orbits cannot bound a Birkhoff section:

\begin{maincorollary}
    \label{Acor:cor E}
    Let $\flow$ be a topologically transitive Anosov flow, and $\Section$ a Birkhoff section of $\flow$.  Then the boundary $\partial\Section$ is not contained
    \begin{itemize}
     \item either in the finite set of pivots orbits (assuming that $\flow$ is not $\RR$-covered)
     \item or in a set corresponding to disjoint simple geodesics (if $\flow$ is the geodesic flow of a hyperbolic surface).
    \end{itemize}
\end{maincorollary}

\paragraph*{Plan of the article.}

In Section \ref{Asec:foliation sec} and  \ref{Section:PartialAndBirkhoffSection}, we describe the notions of orbit space and Birkhoff section. We give several properties common to the three last sections. We also describe Fried-Goodman surgeries and their relations with Theorem \ref{theorem:rcoveredcondition}. The last three sections consist of three independent points of view giving three proofs of the implication $\Rightarrow$ of the Theorem \ref{theorem:rcoveredcondition}. The implication $\Leftarrow$ can be found in Section \ref{Section:Marty}. The three following parts follow respectively T.Marty, M.Asaoka and C.Bonatti's proofs.

Section \ref{Section:Marty} uses a weaker version of Birkhoff section:  positive and negative partial sections. On the one hand we prove that positive partial sections and negative partial sections are mutually excluding. On the other hand we relate the existence of some of these partial sections to the nature of the orbit space of the flow, using previous works from Barbot.

In Section \ref{Section:Asaoka}, we associate to a curve in a Birkhoff section a curve in the orbit space, and vice versa. This procedure, called drifting, allows one to determine the sign of some boundary component of a Birkhoff section when drifting along a subset of the orbit space, called lozenge. 

The last Section \ref{Section:Bonatti} uses the notion of complete quadrants in the orbit space: roughly speaking a region of the orbit space given by two stable/unstable half leaves $l^s,l^u$ (starting at the same point) so that every stable leaf that intersects $l^u$ crosses every unstable leaf that instersects $l^s$. A flow which is not a suspension is positively twisted if and only of all its quadrant $+,+$ are complete. To determine if a quadrant is complete, we look at the holonomy of half a stable leaf in the orbit space, along the unstable foliation. We compare this holonomy to a generalized holonomy we define on a given Birkhoff section. It allows one to prove the completeness of a quadrant when the Birkhoff section is of the positive type.

%% file: Section/Preliminar.tex
\section{The orbit space and lozenges}
\label{Asec:foliation sec}


\subsection{Anosov flows}

Let $M$ be a closed three-dimensional manifold and $\flow$ a smooth flow on $M$. The flow $\flow$ is called \emph{Anosov} if there is a $\phi$-invariant splitting of $TM$ into three line bundles $TM=E^s\oplus X\oplus E^u$ and two real numbers $A,B>0$ such that for one/any Riemannian norm $\|.\|$ on $TM$, we have:
\begin{itemize}
    \item $X$ is tangent to the flow,
    \item for all $t\geq 0$, $\|d_{E^s}\phi_t\|\leq A\exp^{Bt}$,
    \item for all $t\leq 0$, $\|d_{E^u}\phi_t\|\leq A\exp^{B|t|}$. 
\end{itemize}
We fix a smooth Anosov flow $\phi$ on a closed three-dimensional manifold $M$. The flow is called \emph{topological transitive} if there exists an orbit of $\phi$ dense inside $M$. Anosov flows which admit a Birkhoff section or which are twisted $\RR$-covered are topologically transitive.

The bundles $E^s$, $E^u$ and $X$ are integrable to three 1-foliations, which we denote respectively by $\foliation^{ss}$, $\foliation^{uu}$ and $\mcal{O}$. The first two ones are called the strong stable and unstable foliations. The plane bundles $E^s\oplus X$ and $E^u\oplus X$ are also integrable to two 2-foliations, invariant by $\phi$. They are denoted by $\foliation^s$ and $\foliation^u$ and called weak stable and unstable foliations.

Every leaf of the foliations $\foliation^{ss}$ and $\foliation^{uu}$
 is homeomorphic to $\RR$
 and each leaf of the foliations $\foliation^s$ and $\foliation^u$
 is homeomorphic to either
 $\RR/\ZZ \times \RR$ or an open Möbius strip if it contains a periodic orbit,
 or $\RR^2$ otherwise.

 \begin{remark}[on the orientation]
     In this article, we only work with orientable 3-manifolds. Given two local co-orientations of the stable and unstable foliations, and the orientation of the orbit foliation by the flow, their sum is a local orientation on the 3-manifold. Hence, if one of the stable/unstable foliation is co-orientable, the other is also co-orientable. Likewise, given a closed orbit, its stable and unstable are simultaneously orientable or non-orientable. We usually focus on the orientable case, and deal with the non-orientable case using the so called \emph{bundles-orientations} cover. It is a degree two cover constructed as the set of triple $(x,o^s,o^u)$ where $o^s$ and $o^u$ are local co-orientations of the stable and unstable foliations in a neighborhood of $x\in M$, modulo $(x,o^s,o^u)\simeq(x,-o^s,-o^u)$.
 \end{remark}

Given two Anosov flow $\phi$ and $\phi'$ on two 3-manifolds $M$ and $M'$, an \emph{orbit equivalence} (also called topological equivalence) between the flows is a homeomorphism $h\colon M\to M'$ such that for each orbit $\gamma$ of $\phi$, $h(\gamma)$ is an orbit of $\phi'$ and $\gamma\xrightarrow[]{h} h(\gamma)$ preserves the orientations of the orbits by the flows. For orientation purpose, when orientation are fixed on $M$ and $M'$, we suppose that $h$ is additionally orientation preserving. The main notions of this article, that is the topology of the orbit space and existence of Birkhoff section with positive boundary, are invariant by orbit equivalence.

\subsection{The orbit space of an Anosov flow}
\label{Asec:foliation}

 
 Let $P\colon \wt{M} \ra M$ be the universal covering of $M$
 and $\wt{\flow}$ the lift of the Anosov flow $\flow$ to $\wt{M}$.
Let $\wtc{O}$, $\lfoliation^s$, $\lfoliation^u$, $\lfoliation^{ss}$, and $\lfoliation^{uu}$
 be the lifts of the foliations
 $\mcal{O}$, $\foliation^s$, $\foliation^u$, $\foliation^{ss}$, and $\foliation^{uu}$
 respectively.
The foliation $\wtc{O}$ is the orbit foliation of the lifted flow $\wt{\flow}$
 and the flow $\wt{\flow}$ preserves
 $\lfoliation^s$, $\lfoliation^u$, $\lfoliation^{ss}$, and $\lfoliation^{uu}$.
Leaves of the foliations $\wtc{O}$, $\lfoliation^{ss}$, and $\lfoliation^{uu}$
 are homeomorphic to $\RR$
 and leaves of $\lfoliation^s$ and $\lfoliation^u$ are homeomorphic to $\RR^2$.

We call \emph{orbit space} of the flow, or its \emph{bi-foliated plane}, the quotient $\wt{M}/\wtc{O}$ of the universal covering by the orbit foliation. It is denoted by $\orbitspace$. Let $Q\colon \wt{M} \ra \orbitspace$ the quotient map.
The orbit space $\orbitspace$ is known to be a good topological space.

\begin{theorem}[Barbot \cite{Ba1}, Fenley \cite{Fe1}]
The orbit space $\orbitspace$ is homeomorphic to $\RR^2$
 and $Q\colon \wt{M} \ra \orbitspace$ is a locally trivial fibration.
\end{theorem}

Since any leaf of $\lfoliation^s$ is saturated by the orbits of $\wt{\flow}$,
 the two-dimensional foliation $\lfoliation^s$ projects to
 a one-dimensional foliation $\ofoliation^s$ on $\orbitspace$
 such that $Q^{-1}(L)$ is a leaf of $\lfoliation^s$
 for any leaf $L$ of $\ofoliation^s$.
Similarly, the foliation $\lfoliation^u$ projects to a foliation
 $\ofoliation^u$ on $\orbitspace$.

 We remark that an orbit equivalence $h\colon M\to N$ between two Anosov flows $\phi$ and $\psi$ induces an homeomorphism $\wt h\colon P_\phi\to P_\psi$ which sends the weak stable and unstable foliations of $\phi$ to those of $\psi$ \cite{BG10}. It then goes down to an homeomorphism between the orbit spaces preserving the stable and unstable foliations. Hence the homeomorphism class of the bi-foliated plane depends only on the orbit equivalence class of flow.

\begin{figure}[h]
    \begin{center}
        \begin{picture}(120,38)(0,0)
            \put(0,0){\includegraphics[width=12cm]{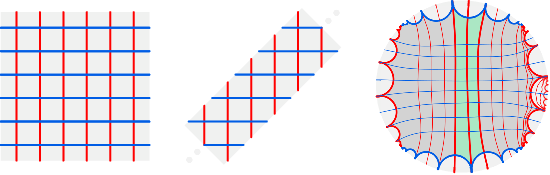}}
            \put(12,-0.5){$\ofoliation^u$}
            \put(33.5,25){$\ofoliation^s$}
        \end{picture}
    \end{center}
    \caption{The bi-foliated plane of three Anosov flow. From left to right: an Anosov suspension, a (positively) twisted Anosov flow and a non-$\RR$-covered Anosov flow. In the third case, the gray region is bounded by four $\ZZ$-families of pairwise non-separated leaves (for either the stable or unstable foliation). This region contains adjacent positive and negative lozenges (see Section \ref{section:lozenges} for the definition of lozenges), two of which are represented in green.}
    \label{Tfigure:orbitSpace}
\end{figure}

It is known that the weak stable and unstable foliations $\foliation^s$
 and $\foliation^u$ are $\Class^1$ foliations (see {\it e.g.} Corollary 9.4.11
 of \cite{FiHa}).
Their $\Class^1$-structures induce transverse $\Class^1$-structures of
 $\ofoliation^s$ and $\ofoliation^u$.
Since $\ofoliation^s$ and $\ofoliation^u$ are mutually transverse,
 we can endow $\orbitspace$ with a $\Class^1$-structure
 such that $\ofoliation^s$ and $\ofoliation^u$ are $\Class^1$ foliations.
For $\xi \in \orbitspace$,
 we denote the leaves of $\ofoliation^s$ and $\ofoliation^u$
 containing $\xi$ by $\ofoliation^s(\xi)$ and $\ofoliation^u(\xi)$.
Since $\orbitspace$ is homeomorphic to $\RR^2$,
Poincaré-Bendixson's theorem implies that
 leaves of $\ofoliation^s$ and $\ofoliation^u$ are homeomorphic to $\RR$
 and the intersection of $\ofoliation^s(\xi)$ and $\ofoliation^u(\eta)$
 contains at most one point for any $\xi,\eta \in \orbitspace$. If one of the foliation $\ofoliation^s$ or $\ofoliation^u$ is trivial, then the other foliation is also trivial, and the flow is called \emph{$\RR$-covered} (see Figure \ref{Tfigure:orbitSpace} for an illustration). 

\begin{theorem}[Barbot \cite{Ba1}, Fenley \cite{Fe1}]
    \label{Theorem:RcoveredOrbitSpaceNonOriented}
    If an Anosov flow $\flow$ is $\RR$-covered, then it is transitive and up to conjugacy by a homeomorphism, the bi-foliated  plane $(\orbitspace,\ofoliation^s,\ofoliation^u)$ is one of the $2$ following possibilities:
    \begin{itemize}
     \item the trivial bi-foliated plane $(\RR^2,H,V)$ where $H$ is directed by $\frac{\partial}{\partial x}$ and $V$ is directed by $\frac{\partial}{\partial y}$,
     \item the restriction of the trivially bi-foliated plane to $\{(x,y)\in\RR^2, |x-y|<1\}$. In this case, the ambient manifold is orientable and one says that $\flow$ is \emph{twisted} (or skewed).
    \end{itemize}
\end{theorem}

Solodov proved that the first case is equivalent to $\flow$ being orbitally equivalent to the suspension flow of a linear Anosov automorphism of $\TT^2$. The geodesic flows on hyperbolic surfaces are examples of twisted $\RR$-covered flows. In Section~\ref{Asec:orientation}, we additionally use the orientation of the ambient manifold to distinguish between two types of twisted flows.

\subsection{Action on the orbit space}

We review three elementary properties on the action of the fundamental group of the ambient manifold on the orbit space.
We denote the fundamental group of $M$ by $\pi_1(M)$
 and its identity element by $e$.
The group $\pi_1(M)$ acts on $\wt{M}$ as deck transformations.
We write $g \cdot \wt{p}$ for the action of $g \in \pi_1(M)$ to
 $\wt{p} \in \wt{M}$.
Since the action preserves the orbit foliation $\wt{O}$,
 this action projects to the quotient space $\orbitspace$ to an action preserving $\ofoliation^s$ and $\ofoliation^u$.
We write $g * \xi$ for the action of~$g \in \pi_1(M)$ to
 $\xi \in \orbitspace$.
\begin{lemma}
Let $g$ be a non-identity element of $\pi_1(M)$
 and $\xi \in \orbitspace$ a fixed point of the action of $g$.
Then, 
 $P(Q^{-1}(\xi))$ is a periodic orbit of $\flow$
 and there exists unique $T \neq 0$ such that
 any $\wt{p} \in Q^{-1}(\xi)$ satisfies that
 $g \cdot \wt{p}=\wt{\flow}_T(\wt{p})$
\end{lemma}
\begin{proof}
Take $\wt{p} \in Q^{-1}(\xi)$.
Then, $Q(g \cdot \wt{p})=g * \xi =\xi =Q(\wt{p})$.
This implies that $g \cdot \wt{p}=\wt{\flow}_T(\wt{p})$ for some $T \in \RR$.
Recall that
 the leaves of the orbit foliation $\wtc{O}$ is homeomorphic to $\RR$.
This means that $\wt{\flow}$ has no periodic orbit,
 and hence, the above $T$ is uniquely determined.
Since $g$ acts on $\wt{M}$ as a non-trivial deck transformation,
 we have $g \cdot \wt{p} \neq \wt{p}$.
This implies that $T \neq 0$.
By projecting to $M$, we have $P(\wt{p})=\flow_T(P(\wt{p}))$.
Therefore, $P(\wt{p}) \in P(Q^{-1}(\xi))$ is a periodic point of $\flow$.
We check that the above $T$ does not depend on the choice of
 $\wt{p} \in Q^{-1}(\xi)$.
For $\wt{q} \in Q^{-1}(\xi)$, there exists $t \in \RR$ such that
 $\wt{q}=\wt{\flow}_t(\wt{p})$.
Then, 
\begin{equation*}
 \wt{\flow}_T(\wt{q})=\wt{\flow}_T(\wt{\flow}_t(\wt{p}))
 =\wt{\flow}_t(\wt{\flow}_T(\wt{p}))=\wt{\flow}_t(g.\wt{p})=g.\wt{\flow}_t(\wt p)=g.\wt{q}.
\qedhere
\end{equation*}
\end{proof}
For $g \in \pi_1(M)$, we define a subset $\Fix(g)_+$ of $\orbitspace$ by
\begin{equation*}
 \Fix(g)_+=\{\xi \in \orbitspace \mid
 g \cdot \wt{p}=\wt{\flow}_T(\wt{p})
 \text{ for some }  \wt{p} \in Q^{-1}(\xi) \text{ and } T>0 \}.
\end{equation*}

\begin{proposition}
\label{Aprop:g action}
For any $g \in \pi_1(M)$ and $\xi \in \Fix(g)_+$,
 there exists a homeomorphism $\alpha\colon  \RR \ra \ofoliation^s(\xi)$
 such that $\alpha(0)=\xi$
 and $g * \alpha(u)= \alpha(2u)$.
\end{proposition}
\begin{proof}
Fix $\wt{p} \in Q^{-1}(\xi)$
 and take  $T>0$ such that $g \cdot \wt{q}=\wt{\flow}_T(\wt{q})$
 for any $\wt{q} \in Q^{-1}(\xi)$.
Since $Q(\wt{p})=Q(\wt{\flow}_{-T}(\wt{p}))=\xi$, we have
\begin{equation*}
 g \cdot \wt{\flow}_{-T}(\lfoliation^{ss}(\wt{p}))
 =\lfoliation^{ss}(g\cdot \wt{\flow}_{-T}(\wt{p}))
 =\lfoliation^{ss}(\wt{\flow}_T \circ \wt{\flow}_{-T}(\wt{p}))
 =\lfoliation^{ss}(\wt{p}).
\end{equation*}
Hence, the map $g \cdot \wt{\flow}_{-T}$ is a homeomorphism
 of $\lfoliation^{ss}(\wt{p})$.
The restriction of $P$ to $\lfoliation^{ss}(\wt{p})$
 is a homeomorphism onto $\foliation^{ss}(P(\wt{p}))$
 and $P \circ (g \cdot \wt{\flow}_{-T})=\flow_{-T} \circ P$.
The restriction of $Q$ to $\lfoliation^{ss}(\wt{p})$
 is a homeomorphism onto $\ofoliation^s(\xi)$
 and $Q \circ (g \cdot \wt{\flow}_{-T})=g * Q$.
These restrictions conjugate
 the action of $\flow_T$ on $\foliation^{ss}(P(\wt{p}))$
 and the action of $g$ on $\ofoliation^s(\xi)$.
More precisely, there exists a homeomorphism
 $H\colon \foliation^{ss}(P(\wt{p})) \ra \ofoliation^s(\xi)$
 such that $H(P(\wt{p}))=\xi$ and $g * H(u) =H \circ \flow_{-T}(u)$
 for any $u \in \foliation^{ss}(P(\wt{p}))$.
Since $\flow$ is an Anosov flow,
 the homeomorphism $\flow_{-T}$ on $\foliation^{ss}(P(\wt{p}))$
 is a uniform expansion.
Therefore, there exists a homeomorphism
 $\beta\colon \RR \ra \foliation^{ss}(P(\wt{p}))$ such that
 $\beta(0)=P(\wt{p})$ and $\flow_{-T} \circ \beta(u)=\beta(2u)$
 for any $u \in \RR$.
Define a homeomorphism  $\alpha\colon \RR \ra \ofoliation^s(\xi)$
 by $\alpha=H \circ \beta$.
Then,  $\alpha(0)=H(\wt{p})=\xi$ and
\begin{equation*}
g * \alpha(u)= g * H(\beta(u))= H \circ \flow_{-T} \circ \beta(u)
 =H \circ \beta(2u)=\alpha(2u)
\end{equation*}
 for any $u \in \RR$.
\end{proof}
\begin{corollary}
Let $g$ be a non-trivial element of $\pi_1(M)$
 and $\xi \in \orbitspace$ a fixed point of the action of $g$.
Then, $P(Q^{-1}(\ofoliation^s(\xi)\setminus\{\xi\}))$ contains
 no periodic orbit of $\flow$.
\end{corollary}

%
%

\subsection{Orientations}
\label{Asec:orientation}

Let $M$ be an oriented closed three-dimensional manifold 
 and $\flow$ a smooth flow on $M$.
An embedded surface $\Sec$ transverse to $\flow$
 is endowed with an orientation induced by
 the orientation of $M$
 and a co-orientation of $\Sec$ given by the flow $\flow$.
An embedded surface $\wt{\Sec}$ in the universal cover $\wt{M}$
 transverse to $\wt{\flow}$ is also endowed with
 an orientation in the same way.

Suppose that $\flow$ is an Anosov flow.
Then, the orbit space $\orbitspace$ of $\wt{\flow}$ admits
 an orientation such that the projection, from any embedded surface
 in $\wt{M}$ transverse to $\wt{\flow}$, to $\orbitspace$
 is orientation preserving.

\begin{remark}
Given a twisted $\RR$-covered Anosov flow on an oriented
 three dimensional manifold $M$,
 we can distinguish whether the twist is positive or negative.
Let $\phi$ be a twisted $\RR$-covered Anosov flow on $M$.
Up to conjugacy by an orientation-preserving homeomorphism,
the bi-foliated plane $(\orbitspace,\ofoliation^s,\ofoliation^u)$
 is one of the $2$ following possibilities:
\begin{itemize}
\item The restriction of the trivially bi-foliated plane to
 $\{(x,y)\in\RR^2, |x-y|<1\}$. 
In this case the $\RR$-covered Anosov flow $\flow$
 is called \emph{positively twisted}.
\item The restriction of the trivially bi-foliated plane to
 $\{(x,y)\in\RR^2, |x+y|<1\}$.
In this case $\flow$ is called \emph{negatively twisted}.
\end{itemize}
Inverting the orientation on $M$ exchanges
 the positively and negatively twisted flows.
The geodesic flow on a hyperbolic flow is known to be negatively twisted
 for the orientation given in local coordinates
 by $(\partial_x,\partial_y,\partial_\theta)$\footnote{The geodesic flow is positively twisted for the orientation provided by the natural Liouville symplectic form of the unitary cotangent bundle.}.
\end{remark}

Since the orbit space $\orbitspace$ is homeomorphic to $\RR^2$,
 the one-dimensional foliations $\ofoliation^s$
 and $\ofoliation^u$ are orientable.
We say that orientations of $\ofoliation^s$ and $\ofoliation^u$
 are \emph{coherent} if  the positive vectors $v^s$ and $v^s$
 along leaves $\ofoliation^s(\xi)$ and $\ofoliation^u(\xi)$
 at a point $\xi$ in $\orbitspace$
 gives a basis $(v^s,v^u)$ which is compatible with
 the orientation of~$\orbitspace$ induced from that of $M$.
We remark that there are two choices of
 a pair of orientations of $\ofoliation^s$ and $\ofoliation^u$.
Given an orientation of $\ofoliation^s$,
 we define the connected components
 $\ofoliation^s_+(\xi)$ and $\ofoliation^s_-(\xi)$
 of $\ofoliation^s(\xi) \setminus \{\xi\}$ for $\xi \in \orbitspace$ by
\begin{align*}
 \ofoliation^s_+(\xi) & = \{\alpha(t) \mid t>0\}, &
 \ofoliation^s_-(\xi) & = \{\alpha(t) \mid t<0\},
\end{align*}
 where $\alpha\colon \RR \ra \ofoliation^s(\xi)$ is an orientation-preserving
 homeomorphism such that $\alpha(0)=\xi$.
Given an orientation of $\ofoliation^s$,
 we also define the connected components
 $\ofoliation^u_+(\xi)$ and $\ofoliation^u_-(\xi)$)
 of $\ofoliation^u(\xi) \setminus \{\xi\}$ in the same way.

\subsection{Lozenges}\label{section:lozenges}
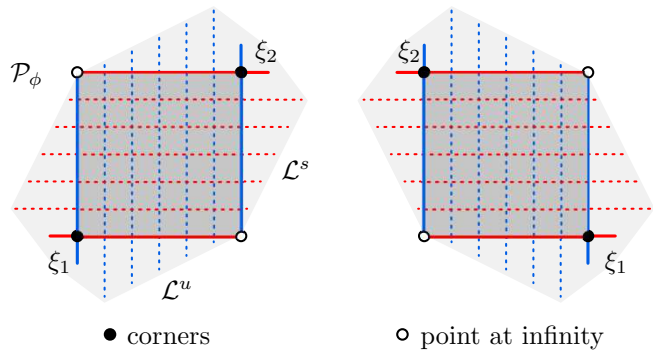
\begin{figure}
    \begin{center}
        \begin{picture}(86,45)(0,0)
            \put(0,0){\includegraphics[height=4.5cm]{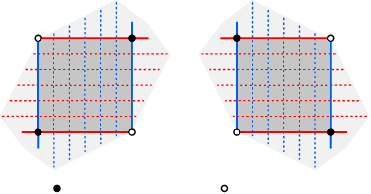}}
            \put(0,35){$\orbitspace$}
            \put(5,10){$\xi_1$}
            \put(32.5,38){$\xi_2$}
            \put(79,10){$\xi_1$}
            \put(51.5,38){$\xi_2$}
            \put(36,22){$\ofoliation^s$}
            \put(20,6){$\ofoliation^u$}
            \put(15.5,0){corners}
            \put(54.5,0){point at infinity}
        \end{picture}
    \end{center}
    \caption{A lozenge of type $(++)$ on the left: $S_{++}(\xi_1)=S_{--}(\xi_2)$ and of type $(--)$ on the right: $S_{-+}(\xi_1)=S_{+-}(\xi_2)$.}
    \label{Afig:lozenge def}
\end{figure}

In \cite{Fe95a}, Fenley
 introduced a lozenge
 as a key tool to investigate topology of the foliations
 $\ofoliation^s$ and $\ofoliation^u$.
For $\xi \in \orbitspace$ and $\sigma,\tau \in \{+,-\}$,
 we define a subset $S_{\sigma\tau}(\xi)$ of $\orbitspace$ by
\begin{alignat*}{8}
S_{\sigma\tau}(\xi)
 =\left\{\left.\eta \in \orbitspace \; \right|\;
 \ofoliation^u(\eta) \cap \ofoliation^s_\sigma(\xi) \neq \emptyset,\,
 \ofoliation^s(\eta) \cap \ofoliation^u_\tau(\xi) \neq \emptyset
 \right\}.
\end{alignat*}
By the transversality of $\ofoliation^s$ and $\ofoliation^u$,
 $S_{\sigma\tau}(\xi)$ is a non-empty open subsets of $\orbitspace$.
We call a subset $L$ of $\orbitspace$ a \emph{lozenge}\footnote{
Lozenges were used in \cite{Fe1} implicitly
 and the explicit definition was given in \cite{Fe95a}.
We adopt the definition in \cite{Fe2},
 which is equivalent to the definition in \cite{Fe95a}.}
 if there exists a pair $(\xi_1,\xi_2)$ of points of $\orbitspace$
 such that $L=S_{++}(\xi_1)=S_{--}(\xi_2)$ or $L=S_{+-}(\xi_1)=S_{-+}(\xi_2)$.
We say that lozenge $L$ is \emph{of type $(++)$}
 in the former case and \emph{of type $(+-)$} in the latter case.
The points $\xi_1$ and $\xi_2$ are called \emph{the corners} of the lozenge $L$.
See Figure~\ref{Afig:lozenge def}.

If a corner $c$ of a lozenge $L$ corresponds to a closed orbit of $\flow$,
 then the orbit is not homotopically trivial \cite{Verjovsky}. So there exists a non-trivial element $g$ of $\pi_1(M)$
 whose action fixes the corner. Up to taking $g^2$ instead of $g$, we can suppose that $g$ preserves the orientation of the stable and unstable leaves of the corner $c$, Then $g$ preserves the lozenges $L$.
The $g$-invariance of the lozenge $L$ implies that
 the action also fix the other corner, and hence,
 the other corner of $L$ corresponds to a closed orbit.
We say that such a lozenge $L$ is \emph{bounded by closed orbits}.
In~\cite{Fe2},
 he gave a characterization of $\RR$-coveredness of an Anosov flow
 in terms of lozenges bounded by closed orbits.

\begin{theorem}
[{\cite[Corollary 4.4]{Fe2}}]
\label{Athm:adjacent lozenge}
If the Anosov flow $\flow$ is not $\RR$-covered,
 then there exist lozenges $L$ and $L'$ bounded by one common closed orbit,
 whose types are $(++)$ and $(+-)$ respectively.\footnote{
 In \cite{Fe2}, Fenley proved a more stronger result. 
 In fact, there exists two lozenges which are `adjacent'.
 See Proof of Claim 3 in Corollary 4.4 of \cite{Fe2}.}
\end{theorem}
\begin{corollary}
\label{Acor:lozenge R-covered}
A topologically transitive Anosov flow is
 positively twisted $\RR$-covered Anosov flow
 if and only if it admits a lozenge of type $(++)$ bounded by closed orbits
 but no lozenge of type $(+-)$ bounded by closed orbits.
\end{corollary}
\begin{proof}
If the flow is positively twisted $\RR$-covered Anosov flow,
 any point of the orbit space is a corner of a lozenge of type $(++)$.
This implies the ``only if'' part of the corollary.

To see the ``if'' part,
 suppose that the flow admits a lozenge of type $(++)$ bounded
 by closed orbits but no lozenge of type $(+-)$ bounded by closed orbits.
By the above theorem, the flow must be $\RR$-covered.
If the flow is non-twisted, then the pair $(\ofoliation^s,\ofoliation^u)$
 has a global product structure, and hence, there exists no lozenges.
If the flow is negatively twisted, then any lozenges are of type $(+-)$.
Therefore, the flow must be positively twisted.
\end{proof}

\section{Partial Sections and Birkhoff sections}
\label{Section:PartialAndBirkhoffSection}

Partial sections is the second main tool which is useful to express topological properties of the flow.
We use the blow up along an orbit to define the invariants of a boundary of a partial section, which are the multiplicity and the linking number. These invariants are also used in the Section \ref{subsection:surgery} to define the Fried-Goodman surgeries.

\subsection{Blow up manifold and local invariant}\label{subsection:Blowing-up}

Let $M$ be a oriented closed three-dimensional manifolds
 and $\flow$ a smooth flow on $M$. For a union $\Gamma$ of finitely many periodic orbits of $\flow$. Let $\pi_\Gamma\colon  M_\Gamma \ra M$ be the blow up along $\Gamma$, obtained by replacing each curve in $\Gamma$ by its unit normal bundle with the suitable topology (also called oriented blow up or dynamical blow up).

We denote by $\TT_\gamma=\pi_\Gamma^{-1}(\gamma)\subset M_\Gamma$ the boundary torus for each curve $\gamma$ of $\Gamma$.
More precisely, $\TT_\gamma$ is the unit normal bundle of $\gamma$,
 that is, $\TT_\gamma =\bigcup_{p \in \gamma}(T_p M \setminus \{0\})/\sim$
 where $v \sim w$ for
 $v,w \in \bigcup_{p \in \gamma}(T_p M \setminus \{0\})$
 if and only if there exist $p \in \gamma$ and $\lambda>0$
 such that $v,w \in T_p M$ and $v-\lambda w \in T_p \gamma$.
By construction, the restriction of $\pi_\Gamma$ to $\TT_\gamma\to\gamma$ is a circle bundle.
For an oriented immersed curve $\gamma^*$ in $\TT_\gamma$, we define the \emph{multiplicity} of $\gamma^*$ by the mapping degree of the restriction of $\pi_\Gamma$ to $\gamma^* \ra \gamma$, where the orientation of $\gamma$ is given by the flow.
We remark that $|\mult(\gamma^*)|$ is the covering degree of $\pi_\Gamma\colon \gamma^* \ra \gamma$ when the blow-down is a covering map from $\gamma^*$ to $\gamma$.


We denote by $\phi^\Gamma$ the lifted flow on $M^\Gamma$. Suppose that the flow $\flow$ is Anosov and
 the stable and unstable leaves of $\gamma$ are orientable.
Then, the restriction of $\flow^\Gamma$ to $\TT_\gamma$ contains
 four periodic orbits that correspond to
 the stable and unstable directions, that is the lifts of $T_\gamma\foliation^{s}(\gamma)$ and $T_\gamma\foliation^{u}(\gamma)$ to the unit normal bundle of $\gamma$.
Let $\lambda_\gamma$ be one of these periodic orbit, 
 and $\mu_\gamma$ be a fiber of $\pi_\Gamma\colon T_\gamma\to\gamma$.
The closed curves $\lambda_\gamma$ and $\mu_\gamma$ are called
 the \emph{parallel} and the \emph{meridian} of $\TT_\gamma$.
We orient $\lambda_\gamma$ by the direction of the flow $\phi$
 and $\mu_\gamma$ by the orientation of the normal bundle of $\gamma$
 induced from the co-orientation of the flow $\phi$.
Then, the homology classes $[\lambda_\gamma]$ and $[\mu_\gamma]$
 do not depend on the choice of the curves
 $\lambda_\gamma$ and $\mu_\gamma$,
 and they form a basis of $H_1(\TT_\gamma,\ZZ)$.
For an {\it oriented} closed curve $\gamma^*$ in $\TT_\gamma$,
 we can write the homology class $[\gamma^*]$ of $\gamma^*$
 in $H_1(\TT_\gamma,\ZZ)$ as $\mult(\gamma^*)[\lambda_\gamma]+k[\mu_\gamma]$
 for a unique integer $k$. We call \emph{linking number} of $\gamma^*$ the integer $k$,
 and denote it by $\link(\gamma^*)$.

 When the stable and unstable leaves of $\gamma$ are not-orientable, one could also define the linking number similarly. It turns out to be less efficient to work with, so we usually resolve any question on these orbits by lifting to the bundles-orientations covering.

\subsection{Partial sections and Birkhoff sections}\label{subsection:Partial sections}



We introduce partial sections, Birkhoff sections and their first return map. Partial section are only relevant for Section \ref{Section:Marty} and the first return map is used in Sections \ref{Section:Asaoka},\ref{Section:Bonatti}. 

\paragraph*{Partial sections.}

Let $M$ be an oriented three-dimensional manifold
 and $\flow$ a smooth flow on $M$.
We call a subset $\Sec$ of $M$ a \emph{partial section}
 if it is the image of a smooth immersion 
 $\iota$ from a compact surface $\wh{\Sec}$ to $M$
 such that the restriction of $\iota$ to $\Int \wh{\Sec}$ 
 is an embedding transverse to the flow $\flow$ 
 and $\iota(\del\wh{\Sec})$ is a finite union of closed orbits of $\flow$,
 where $\Int \wh{\Sec}$ and $\del\wh{\Sec}$ are the interior
 and the boundary of the surface $\wh{\Sec}$. 
We call the immersion $\iota$ a \emph{parametrization} of $\Sec$
 and the sets $\iota(\Int \wh{\Sec})$ and $\iota(\del \wh{\Sec})$
 are called the \emph{interior} and the \emph{boundary} of~$\Sec$
 and denoted by $\Int\Sec$ and $\del\Sec$ respectively.
The parametrization $\iota$ lifts 
 to an immersion $\iota_{\del\Sec}$ into the blow up $M_{\del\Sec}$
 along $\del \Sec$, such that $\iota=\pi_{\partial \Sec}\circ\iota_{\partial\Sec}$.
The image $\iota_{\del\Sec}(\wh{\Sec})$ is called
 the \emph{resolution} of $\Sec$ and denoted by $\Sec^*$.

 Each closed orbit $\gamma$ in $\del \Sec$ is the image $\iota(\wh{\gamma})$
  for some connected component $\wh{\gamma}$ of $\del \wh{\Sec}$.
 The lift $\iota_{\del \Sec}(\wh{\gamma})$
  is called a \emph{boundary component} immersed in $\gamma$.
 We remark that each boundary component $\gamma^*$ immersed in an orbit $\gamma$ is transverse to the fibers of the restriction of $\pi_{\del\Sec}$ to $\TT_\gamma\to\gamma$.
 Since $\Int \Sec$ is transverse to the flow $\flow$,
  it admits the orientation given in Subsection \ref{Asec:orientation}.
 It induces an orientation of $\Sec^*$ as an immersed surface,
  we endow each boundary component 
  with an orientation as a boundary curve of $\Sec^*$.

\paragraph*{Birkhoff sections.}

Let $\flow$ be a topologically transitive Anosov
 on a oriented closed three-dimensional manifold $M$.
A \emph{Birkhoff section} is a partial section of $\flow$
 such that the length of orbit segments disjoint from $\Sec$ is bounded,
 {\it i.e.}, $\bigcup_{t \in [0,T]}\flow_t(\Sec)=M$ for some $T>0$. Birkhoff section are the partial section for which a first return map, as define further, is well-defined. Any topologically transitive Anosov flows admits a Birkhoff section \cite{Fri}. 

 We say a Birkhoff section $\Sec$ is \emph{tame} if each boundary component in $M_{\partial \Sec}$ is transverse to the closed orbits in the boundary torus
  which corresponds to the stable and the unstable directions.
 In this paper, we assume that the Birkhoff sections are smooth and tame.
 Given a topological surface with Birkhoff like properties,
  under a tame hypothesis on the surface
  one can find an isotopy along the flow from the interior
  of the surface to the interior of a smooth and tame Birkhoff section.
 See \cite{BG10}.

\begin{figure}[h]
    \begin{center}
        \begin{picture}(120,35)(0,0)
            \put(0,0){\includegraphics[width=12cm]{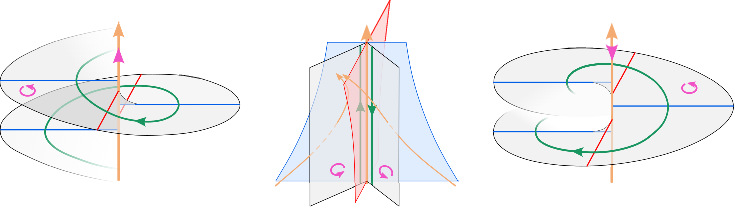}}
            \put(5.5,0){$\link<0$, $\mult>0$}

            \put(54,-3){$\link=0$}
            \put(38.5,24){$\mult>0$}
            \put(66,24){$\mult<0$}
            
            \put(86,0){$\link<0$, $\mult<0$}
        \end{picture}
    \end{center}
    \caption{Linking numbers and multiplicity of four boundaries of partial sections. If $\gamma^*$ denote the boundary component of the surface, then is represented in green a small perturbation of $\gamma^*$ by pushing it inside the interior of the surface.}
    \label{Tfigure:LK}
\end{figure}
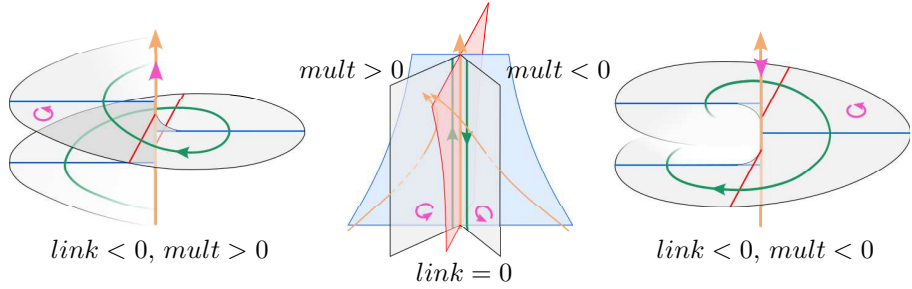

 Around a boundary component, one may think of a Birkhoff section as finitely many parallel helicoids supported by a periodic orbit. Local models for partial section along a boundary components are illustrated in Figure~\ref{Tfigure:LK}, together with their local invariants that are defined below.
 \paragraph*{First return map.}
 

 Let $\flow$ be a transitive Anosov flow with oriented bundles, endowed with a Birkhoff section $\Section$. We denote by $\wh \Section_\partial$ the closed surface obtained from $\Section$ by collapsing to a point each boundary component of $\Section$, and $\Sigma\subset \wh \Section_\partial$ is the finite set corresponding to $\partial \Section$. \cite{Fri}  noticed that the return map $\phi\colon \Section\to \Section$  is well-defined and induces on $\wh \Section_\partial$ a pseudo-Anosov homeomorphism  $\returnmap$ on $(\wh \Section_\partial,\Si)$ so that $\returnmap(\Si)=\Si$, and  whose singular points belong to $\Si$.

 \begin{remark}
 \begin{itemize}
    \item The pseudo-Anosov map obtained
    above has one restriction: there are exactly two $\returnmap$-orbits of \emph{prongs} (also called \emph{separatrixes}) are each singularity.  
    \item a point $\sigma\in \Si$ is periodic of period $k$ if and only if $k$ boundary components of $\Section$ are mapped on the periodic orbit $\gamma$ of $\flow$ corresponding to $\sigma$. 
 \end{itemize}
 \end{remark}

 \subsection{Local invariants}\label{subsection:LocalInvariant}

 A boundary component $\gamma^*$ of $\Sec$ 
  has an invariant: the \emph{multiplicity} $\mult(\gamma^*)$.
 Since $\pi_{\del\Sec}$ is a covering map from $\gamma^*$
  to a closed orbit in $\del \Sec$, the multiplicity $\mult(\gamma^*)$
  is non-zero.
 We say that $\gamma^*$ is \emph{positive} ({\it resp.} \emph{negative})
  if we have $\mult(\gamma^*)>0$ ({\it resp.} $\mult(\gamma^*)<0$).
 By extension a partial section is called \emph{positive} ({\it resp.}
  \emph{negative}) when all its boundary components are positive
  ({\it resp. negative}).
Notice that $\gamma^*$ is a positive if and only the restriction of $\pi_{\partial S}$ to $\gamma^*\to\gamma$ is orientation preserving.

 When the stable and unstable foliations of $\flow$ is orientable,
  a boundary component $\gamma^*$ of $\Sec$
  has another invariant; the \emph{linking number} $\link(\gamma^*)$.
  
  We remark that given $h\colon M\to N$ an orbit equivalence between two Anosov flow $\phi$ and $\psi$, and $\Section$ is a tame Birkhoff section of $M$, one can smooth $h(\Sigma)$ into a tame Birkhoff section of $\phi$ bounded by $h(\partial \Sigma)$ and with the same local invariant at each boundary component.

 \begin{lemma}
 \label{Lemma:NonPositiveLink}
 Let $\Sec$ be a partial section and $\gamma^*$ be a boundary component
  such that $\link(\gamma^*)$ is well-defined (the stable and unstable leaves of $\pi_{\partial S}(\gamma^*)$ are orientable).
 Then, $\link(\gamma^*) \leq 0$.
 If $\Sec$ is a Birkhoff section then $\link(\gamma^*)<0$.
 \end{lemma}

 \begin{proof}
Denote by $S^*$ the resolution of $S$ into the blow up manifold $\subset M_{\partial S}$, and $\lambda_\gamma\subset\TT_\gamma$ one closed curve corresponding to the stable direction in the torus boundary. If $\gamma^*$ and $\lambda_\gamma$ are disjoint, then the linking number is zero. So we only need to consider the case when they do intersect.

 Because we assume $\Sec$ to be tame, $\gamma^*$ intersects transversally $\lambda_\gamma$. Let $x$ be an intersection point between these curves. We fix three vectors based on $x$, $u^*$ tangent to $\gamma^*$ and positively oriented, $v^s$ tangent to $\lambda_\gamma$ and positively oriented, and $v_{out}$ normal to the torus $\TT_\gamma$ and going out of $M_{\partial S}$. 
 
 We choose the orientations such that $(v_{out},u^*)$ is a positive basis of $T_xS^*$, and such that $(u^*,v^s,-v_{out})$ is a positive  basis if and only if the algebraic intersection of $\gamma^*$ and $\lambda_\gamma$ at $x$ is positive. Since $S^*$ is transverse to the flow in a neighborhood of $x$, the basis $(v_{out},u^*,v^s)$ is positive. So the basis $(u^*,v^s,-v_{out})$ is negative. Hence, the algebraic intersection at $x$ between $\gamma^*$ and $\lambda_\gamma$ is negative. Hence, $\link(\gamma^*)=\gamma^*\algcap \lambda_\gamma$ is non-positive.
 
 When $S$ is a Birkhoff section, take $y\in M$ in the half stable leaf of $\lambda_\gamma$. Then there is a sequence of time $t_n$ going to $+\infty$ such that $\phi_{t_n}(y)$ is inside $S$ for all $n$. By compactness, $\pi_{\partial S}^{-1}(\phi_{t_n}(y))$ converges toward a point inside $\partial S^*\cap \lambda_\gamma$. Hence, this intersection is non-empty and the linking number is non-zero. 
 
     
 \end{proof}
 

 \begin{lemma}\label{Lemma:MultLinkWellDefined}
     For a Birkhoff section $\Sec$ and a closed orbit $\gamma$ in $\del\Sec$,
      the multiplicity $\mult(\gamma^*)$ does not depend on the choice of a boundary component $\gamma^*$ immersed in $\gamma$. 
      
      When the flow is Anosov and the linking number $\link(\gamma^*)$ is well-defined, then it does not depend on the choice of a boundary component $\gamma^*$ immersed in~$\gamma$.
 \end{lemma}
 
 \begin{proof}[Proof in a special case]
     We suppose that the flow is Anosov and that the stable and unstable leaves of $\gamma$ are orientable. The general case is proved in Section~\ref{Tsection:GeneralFlow}. Suppose that two boundary components $\gamma_1^*$ and $\gamma_2^*$
     immersed in $\gamma$ have different multiplicities or linking numbers.    
     We denote by $[\gamma_i^*]$ the homology of the curve $\gamma_i^*$
     in the torus $T$ obtained by blowing up $\gamma$. Since a partial section is immersed in its interior, the curve $\gamma_i^*$ are homotopic to simple curves.
     The linking number and the multiplicity form a coordinates system
     in $H_1(T,\ZZ)$, so $\mult(\gamma_i^*)$ and $\link(\gamma_i^*)$ are co-prime. 
     Hence $[\gamma_1^*]$ and $[\gamma_2^*]$ are either opposite
     or non-proportional. 
 
     Since $\Sec$ is a Birkhoff section, according to Lemma \ref{Lemma:NonPositiveLink}
     the linking numbers of $\gamma_1^*$ and $\gamma_2^*$ are negative,
     so $[\gamma_1^*]$ and $[\gamma_2^*]$ are not opposite one to the other. 
     Hence they are not proportional, and their algebraic intersection inside $T$
     is non-zero.
     This means that annular neighborhoods of $\gamma_1*$ and $\gamma_2^*$
     have non-trivial intersection in the interior of $\Sec$.
     It contradicts that $\pi_\Gamma(\Int \Sec^*)=\Int \Sec$ is an embedded surface in $M$.
 \end{proof}
 
 A complete proof is provided in Lemma \ref{Lemma:MultLinkWellDefined}, in Section \ref{Tsection:GeneralFlow}, for Anosov flows with non-orientable foliations and for non-Anosov flows.
 
     So, when $\Sec$ is a Birkhoff section, we define the \emph{multiplicity} $\mult(\gamma)$
      and the \emph{linking number} $\link(\gamma)$
      of a closed orbit $\gamma$ in $\del \Sec$ by
      $\mult(\gamma)=\mult(\gamma^*)$
      and $\link(\gamma)=\link(\gamma^*)$
      with a boundary component $\gamma^*$ of $\Sec$ immersed in $\gamma$.
     When~$\Sec$ is only a partial section, the lemma is false. Indeed two boundary components with linking number zero and multiplicity $+1$ and $-1$ can be immersed in the same closed orbit. Hence for partial section, the multiplicity and linking number depends on the immersed boundary component.

\subsection{Lift of a Birkhoff section to the universal cover}

Sections \ref{Section:Asaoka} and \ref{Section:Bonatti} uses crucially the relation between the lift of a Birkhoff sections in the universal cover and the orbit space. We give here the key property.

Let $S$ be a Birkhoff section. By definition there exists $T>0$ such that
 $\{\flow_t(p) \mid 0 \leq t \leq T\}$ intersects with $\Section$
 for any $p \in M$. For $p \in \Int \Section$,
 we define the \emph{return time sequence}
 $(\tau_n(p))_{n \in \ZZ}$ by 
\begin{alignat*}{5}
 \tau_0(p) & = 0, & \hsp
 \tau_{n+1}(p)=\min\{t>\tau_n(p) \mid \flow_t(p) \in \Int \Section \}
\end{alignat*}

The return time sequence exists for any $p \in \Int\Section$ and it satisfies that
\begin{alignat*}{5}
 \tau_n(p) & <\tau_{n+1}(p)\leq \tau_n(p)+T, & \hsp
 \lim_{n \ra \pm \infty} \tau_n(p) & = \pm \infty.
\end{alignat*}

Since the surface $\Int \Section$ is transverse to the orbit foliation $\mcal{O}$,
 the function $\tau_n$ on $\Int \Section$ is continuous for each $n$. We define {\it the return map} $\returnmap_\Section\colon \Int \Section \ra \Int \Section$ by $\returnmap_\Section(p)=\flow_{\tau_1(p)}(p)$. This is a homeomorphism of $\Int \Section$ and it satisfies that $\returnmap_\Section^n(p)=\flow_{\tau_n(p)}(p)$ for any $n \in \ZZ$ and $p \in \Int \Section$.

 \vspace{1\baselineskip}

 Let $\Gamma\subset M$ be a finite union of closed orbits. We denote by $\wt{\pi}_\Gamma\colon  \wt{M}_\Gamma \ra \wt{M}$
  the blow up of the universal covering $\wt{M}$ along $P^{-1}(\Gamma)$ 
  and by $\Pi_\Gamma\colon \orbitspaceG \ra \orbitspace$ 
  the blow up of the orbit space $\orbitspace$ at $Q(P^{-1}(\Gamma))$.
 Then, the projections $P\colon \wt{M} \ra M$ and $Q\colon \wt{M} \ra \orbitspace$
  lifts to $P_\Gamma\colon \wt{M}_\Gamma \ra M_\Gamma$
  and $Q_\Gamma\colon \wt{M}_\Gamma \ra \orbitspaceG$ between the blow ups.

Recall that $P\colon \wt{M} \ra M$ is the universal covering of $M$. Define the {\it lifts} $\Int \wt{\Section}$ and $\del \wt{\Section}$ of $\Int \Section$ and $\del \Section$ to $\wt{M}$ by
\begin{alignat*}{4}
 \Int \wt{\Section} & = P^{-1}(\Int \Section), & \hsp
 \del \wt{\Section} & = P^{-1}(\del \Section).
\end{alignat*}
We also define the subsets $\Int\wt{\Section}/\wtc{O}$ and $\del\wt{\Section}/\wt{O}$ of $\orbitspace$ by
\begin{alignat*}{4}
\Int \wt{\Section}/\wtc{O} & = Q(\Int \wt{\Section}) &, \hsp
\del \wt{\Section}/\wtc{O} & = Q(\del \wt{\Section}).
\end{alignat*}
Then, we have
\begin{alignat*}{4}
\Int \wt{\Section}
 & = Q^{-1}(\Int \wt{\Section}/\wt{O}) \cap P^{-1}(\Section), & \hsp
\del \wt{\Section} & = Q^{-1}(\del\wt{\Section}/\wtc{O})
\end{alignat*}
 and the orbit space $\orbitspace$ of $\wt{\flow}$ is the disjoint union of $\Int \wt{\Section}/\wtc{O}$ and $\del\wt{\Section}/\wtc{O}$ (recall that $\Sec$ is a Birkhoff section). Since $\del \Section$ consists of finite periodic orbits, the set $\del \wt{\Section}$ is the union of isolated orbit lines of $\wt{\flow}$. Hence, $\del \wt{\Section}/\wt{O}$ is a discrete subset of~$\orbitspace$.

Define a function $\wt{\tau}_n\colon \Int \wt{\Section} \ra \RR$ and a homeomorphism $\wt{\returnmap}_\Section$ of $\Int \wt{\Section}$ by
\begin{alignat*}{5}
 \wt{\tau}_n(\wt{p}) & = \tau_n(P(\wt{p})), & \hsp
 \wt{\returnmap}_\Section(\wt{p}) & = \wt{\flow}_{\wt{\tau}_1(\wt{p})}(\wt{p})
 =\wt{\flow}_{\tau_1(P(\wt{p}))}(\wt{p})
\end{alignat*}
 for $\wt{p} \in \Int \wt{\Section}$.
We denote the restriction of $Q$ to $\Int \wt{\Section}$ by $Q_\Section$.
The homeomorphism $\wt{\returnmap}_\Section$ satisfies that
\begin{alignat*}{5}
 P \circ \wt{\returnmap}_\Section & = \returnmap_\Section \circ P, \hsp
 Q_\Section \circ \wt{\returnmap}_\Section & = Q_\Section.
\end{alignat*}
In particular, $\wt{\returnmap}_\Section$ is a lift of the return map $\returnmap_\Section$.

\begin{proposition}
\label{Aprop:covering}
The map $Q_\Section$ is a covering map onto $\Int \wt{\Section}/\wtc{O}\subset\orbitspace$.
The group of deck transformations of the covering is
 generated by $\wt{\returnmap}_\Section$.
\end{proposition}

\begin{proof}
Take $\xi \in \Int \wt{\Section}/\wtc{O}$
 and $\wt{p} \in Q_\Section^{-1}(\xi)$.
Since $\Int \wt{\Section}$ is transverse to $\wt{\flow}$,
 there exists a neighborhood $U$ of $\wt{p}$ in $\Int \wt{\Section}$
 such that the restriction of $Q_\Section$ to $U$
 is a homeomorphism onto an open subset $V$ of $\orbitspace$. 
By replace $U$ with its connected component if necessary,
 we may assume that $U$ and $V$ are connected. 
For each $\eta \in V$ and $\wt{q} \in Q_\Section^{-1}(\eta) \cap U$, we have 
 \begin{equation*}
 Q_\Section^{-1}(\eta)
 =\{\wt{\returnmap}_\Section^n(\wt{q}) \mid n \in \ZZ\} 
\end{equation*}
 by the definition of the return time sequence $(\tau_n(\wt{q}))_{n \in \ZZ}$.
This implies that
 $Q_\Section^{-1}(V)=\bigcup_{n \in \ZZ}\wt{\returnmap}_\Section^n(U)$. 
We claim that
 $\wt{\returnmap}_\Section^m(U) \cap \wt{\returnmap}_\Section^n(U)=\emptyset$
 if $m \neq n$. 
Suppose that
 $\wt{\returnmap}_\Section^m(\wt{q}_1)=\wt{\returnmap}_\Section^n(\wt{q}_2)$
 for some $\wt{q}_1,\wt{q}_2 \in U$ and $m,n \in \ZZ$. 
This implies that $Q(\wt{q}_1)=Q(\wt{q}_2)$. 
Since the restriction of $Q$ to $U$ is homeomorphism, 
we have $\wt{q}_1=\wt{q}_2$. 
Since the flow $\wt{\flow}$ has no periodic orbit, we have $m=n$. 
Therefore, 
$\wt{\returnmap}_\Section^m(U) \cap \wt{\returnmap}_\Section^n(U)=\emptyset$
 if $m \neq n$. 
In particular, $\wt{\returnmap}_\Section^n(U)$ is
 a connected component of $Q_\Section^{-1}(V)$. 
Recall that the restriction of $Q_\Section$ to $U$ is
 a homeomorphism onto $V$ and $\returnmap_\Section$ is
 a homeomorphism of $\Int \wt{\Section}$. 
This implies that the restriction of
 $Q_\Section=Q_\Section \circ \wt{\returnmap}_\Section^{-n}$
 to $\wt{\returnmap}_\Section^n(U)$ is a homeomorphism onto $V$. 
Therefore, the map $Q_\Section$ is a covering map. 
Since $Q_\Section \circ \wt{\returnmap}_\Section=Q_\Section$, 
 the map $\wt{\returnmap}_\Section$ is a deck transformation
 of the covering. 
The equation
 $Q_\Section^{-1}(\eta)
 =\{\wt{\returnmap}_\Section^n(\wt{q}) \mid n \in \ZZ\}$
 implies that the set of deck transformations of the covering is
 $\{\wt{\returnmap}_\Section^n \mid n \in \ZZ\}$.
\end{proof}

\subsection{Fried-Goodman surgeries}\label{subsection:surgery}

Let $M$ be an oriented closed three-dimensional manifold
 and $\flow$ a smooth Anosov flow on $M$. We take $\gamma$ a closed orbit and $M_\gamma$ is the blow up of $M$ along $\gamma$. Recall that $\pi_\gamma\colon M_\gamma\to M$ is the blowing-down projection and $T_\gamma = \pi^{-1}(\gamma)$ is a torus boundary of $M_\gamma$. Also denote by $\flow^\gamma$ the lift of the flow $\flow$ to $M_\gamma$.
 
\paragraph*{Fried-Goodman surgery on orientable orbits.}
We suppose that $\gamma$ has orientable stable and unstable leaves.
We give the definition of Fried surgeries which have been introduced in \cite{Fri}.

A Fried surgery along $\gamma$ with coefficient $k \in \ZZ$
 is the blow-down of $M_\gamma$ along the circle foliation
 on $\TT_\gamma$ which is transverse to
 the flow $\flow^\gamma$ and each circular leaf represents
 the homology class $k[\lambda_\gamma]+[\mu_\gamma]$\footnote{
In the standard terminology of three-dimensional manifold theory,
 a Fried-Goodman surgery of coefficient $k$ along $\gamma$
 is the Dehn surgery of coefficient $1/k$ along $\gamma$
 with respect to the framing give by $\flow$
 (see \cite[Section 5.3]{GS} for instance).}.
Let $M_{\gamma,k}$ be the manifold obtained  by the surgery
 and  $\pi_{\gamma,k}\colon M_\gamma \ra M_{\gamma,k}$ the projection.
As described in \cite{Fri}, $\flow^\gamma$ blows down to a
 topological Anosov flow on $M_{\gamma,k}$.
Shannon \cite{Sh} proved that Fried surgery is equivalent to
 another surgery given by Goodman in \cite{Go} up to orbit equivalence,
 and hence, the topological flow $\flow'$ is orbit equivalent to
 a smooth Anosov flow $\psi$ when the original flow $\flow$ is topologically transitive. We say that $\psi$ is obtained by \emph{Fried-Goodman} surgery from $\phi$ along $\gamma$.

The following lemma describe
 how the multiplicity and longitude of a closed curve in $\TT_\gamma$
 change by a Fried-Goodman surgery.
\begin{lemma}
\label{Alemma:surgery slope} 
For a periodic orbit $\gamma$ of $\flow$ and an integer $k$,
 let $M_\gamma$ be the blow up along $\gamma$,
 $\pi_{\gamma,k}\colon M_\gamma \ra M_{\gamma,k}$ the projection,
 and $\flow^{\gamma,k}$ be the Anosov flow obtained by
 the Fried-Goodman surgery along $\gamma$ with coefficient $k$.
Let $\gamma_k$ be the periodic orbit $\pi_{\gamma,k}(\TT_\gamma)$
 of $\flow^{\gamma,k}$.
For any closed curve $\gamma^*$
 in $\TT_\gamma=\TT_{\gamma,k}=\pi_{\gamma_k}^{-1}(\gamma_k)$,
 we have
\begin{align*}
\mult(\gamma^*,\flow^{\gamma,k})
 & = \mult(\gamma^*,\flow)-k\cdot\link(\gamma^*,\flow), &
\link(\gamma^*,\flow^{\gamma,k}) & = \link(\gamma^*,\flow).
\end{align*}
\end{lemma}
\begin{proof}
The projection $\pi_{\gamma,k}\colon \TT_\gamma \ra \gamma_k$
 maps a curve which represents the homology class
 $k[\lambda_\gamma]+[\mu_\gamma]$ to a point.
This implies that
 $[\mu_{\gamma_k}]=[\mu_\gamma]+k[\lambda_\gamma]$.
Since the parallel $\lambda_\gamma$ corresponds the stable direction
 of $\flow^{\gamma,k}$, we also have
 $[\lambda_{\gamma_k}]=[\lambda_\gamma]$.
Hence,
\begin{align*}
 [\gamma^*] & =\mult(\gamma^*,\flow)[\lambda_\gamma]
 +\link(\gamma^*,\flow)[\mu_\gamma]\\
 & = (\mult(\gamma^*,\flow)-k\cdot \link(\gamma^*,\flow))[\lambda_{\gamma_k}]
 +\link(\gamma^*,\flow)[\mu_{\gamma_{k}}].
\end{align*}

This means that
\begin{align*}
\mult(\gamma^*,\flow^{\gamma,k})
 & = \mult(\gamma^*,\flow)-k\cdot\link(\gamma^*,\flow), &
\link(\gamma^*,\flow^{\gamma,k}) & = \link(\gamma^*,\flow).
\qedhere
\end{align*}
\end{proof}
 
\begin{proposition}
    Let $\flow$ be a topologically transitive Anosov flow
     on an oriented closed three-dimensional manifold, with orientable stable and unstable foliations,
     and $\Sec$ a Birkhoff section of $\flow$.
    Then, there is a Fried-Goodman surgery along $\del \Sec$ with
     positive coefficients which makes $\Sec$ a positive Birkhoff section.
\end{proposition}

\begin{proof}
    By Lemma \ref{Alemma:surgery slope},
     a $(1/k)$-Fried-Goodman surgery along a closed orbit in $\del \Sec$
     increase the multiplicity of a boundary component $\gamma^*$
     by $-k \cdot \link(\gamma^*)$.
    Since $\link(\gamma^*)$ is negative by the above lemma,
     the multiplicity becomes positive if $k$ is sufficiently large.
\end{proof}

\paragraph*{Fried-Goodman surgeries on periodic orbits with negative eigenvalues.}
Let $\flow$ be a transitive Anosov flow on an oriented manifold $M$, whose strong stable and unstable bundles are not orientable.  Thus the bundles-orientations cover is a two-sheets cover $\wb M$ and the lift of $\flow$ on $\wb M$ is a transitive Anosov flow~$\wb\flow$ with oriented stable and unstable bundles. 


Consider a closed orbit $\gamma$ of the flow. If the stable and unstable leaf of $\gamma$ are orientable, then Fried surgeries and Goodman surgeries can be defined on~$\gamma$ as in the previous section. Suppose that the stable and unstable leaves are not orientable. We need to adapt the previous definition. The homology of the torus~$\TT_\gamma$ does not have a coordinate system given by the multiplicity and the linking number anymore. Indeed the close 1-manifold in the stable direction of~$\gamma$ on $\TT_\gamma$ is one connected curve with multiplicity 2 in this case. 

To do a Fried surgery, we consider one simple closed curve $\delta\subset\TT_\gamma$ transverse to the induced flow, and which intersects the stable direction of $\TT_\gamma$ exactly twice. Then the blow down of $M_\gamma$ along one foliation of transverse curves parallel to $\delta$ gives a closed manifold with a flow which is only continuous on a neighborhood of the blow down orbit. Shannon \cite{Sh} proved that this flow is topologically equivalent to a smooth Anosov flow.

Equivalently one can consider the bundles-orientations cover $\wb M\to M$ which is a degree two cover and take $\wb\gamma$ the preimage of $\gamma$ ($\wb\gamma\to\gamma$ is a degree 2 map). We do a Fried surgery (as defined in the previous section) along a foliation on $\TT_\gamma$ which is invariant by the deck transformation map $f\colon \wb M\to \wb M$ of the covering map $\wb M\to M$. The manifold $\wb N$ obtained by surgery admits an involution $g\colon \wb N\to \wb N$, which is equal to the map $f$ when restricted to $\wb N$ minus the surgery orbit. The manifold $N=\wb N/g$ is equivalently obtained by a Fried surgery on the orbit $\gamma$ in the initial manifold $M$, as described in previous paragraph.

With the same idea, Goodman surgeries can be extended to orbits with non orientable stable/unstable foliations by taking the bundles-orientations cover, doing two Goodman surgeries which are invariant by the deck transformation, and taking the quotient by the new deck transformation.

The sign of a Fried-Goodman surgery in the non orientable case is well-defined. Since the initial manifold is oriented, the bundles-orientations cover is endowed with an orientation such that the cover map is orientation preserving. For any foliation used to do a Fried surgery invariant by the deck transformation, the sign on the foliation is invariant by the deck transformation. So for a Fried surgery on $M$, its sign is defined as the sign of the induced Fried surgery on the bundles-orientations covering manifold. 

\vspace{\baselineskip}

Assuming the main Theorem \ref{theorem:rcoveredcondition}, we can prove its corollary: any transitive Anosov flow on an oriented 3-manifold can be obtained by Fried-Goodman surgeries from a positively twisted $\RR$-covered Anosov flow.
\begin{proof}[Proof of Corollary \ref{Acor:cor B}]
Take $\phi$ a transitive Anosov flow on an oriented 3-manifold $M$. Fried proved that $\phi$ admits a Birkhoff section $\Section$ \cite{Fri}with non-empty boundary. Suppose first that each orbit in $\partial \Section$ has orientable stable and unstable leaves. 

For each orbit $\gamma$ in the boundary of $\Section$, fix an integer $k_\gamma\in\ZZ$ such that $\mult(\gamma,\flow)-k_\gamma\cdot\link(\gamma,\flow)>0$. Denote by $N$ and $\psi$ the 3-manifolds and the Anosov flow obtained by Fried-Goodman surgeries with coefficients $k(\gamma)$ along the orbit $\gamma$, for each orbit $\gamma\subset\partial \Section$. Denote by $\Section_N\subset N$ the image of the surface $\Section$ under the surgery, which we can suppose smooth after a small deformation. For an orbit $\gamma\in M$, we denote by $\gamma_N$ the corresponding orbit inside $N$ after surgeries.

According to Lemma \ref{Alemma:surgery slope}, for every orbit $\gamma\subset\partial\Sigma$, the multiplicity of $\Section_N$ along $\gamma_N$ is equal to $\mult(\gamma,\flow)-k_\gamma.\link(\gamma,\flow)$ which is positive. Hence $\Section_N$ has only positive boundary components. Which implies that $\psi$ is positively twisted $\RR$-covered thanks to Theorem \ref{theorem:rcoveredcondition}.

Consider now the case when the stable and unstable foliation of $M$ are not orientable. We denote by $\wb M$ the bundles-orientations cover of $M$, and $\wb \phi$ the lifted flow of $\phi$ inside $\wb M$. Denote by $g\colon \wb M\to\wb M$ the deck transformation. One can apply the previous argument to $\wb M$ with the following additional property. 

Take a Birkhoff section $\Section$ of $\phi$ (with non-empty boundary) and lift it in a Birkhoff section $\wb\Section\subset \wb M$ for $\wb\phi$. For $\gamma_{\wb M}$ a boundary component of $\wb\Sigma$, $g(\gamma_{\wb M})$ is equal either to $\gamma_{\wb M}$ or to another orbit. In the first case, do one Fried-surgery along $\gamma_{\wb M}$ given by a curve invariant by the deck transformation $g$. In the second case, do two symmetric Fried-surgeries along $\gamma_{\wb M}$ and $g(\gamma_{\wb M})$. In both case, choose the slopes of the surgeries such that the multiplicity of $\Section$, after surgery, is positive along all its boundary components. 

The flow obtained by these Fried surgeries is positively twisted $\RR$-covered, and the ambient manifold admits an involution $g'$ induced by $g$. Since the orbit space of an Anosov flow is invariant by covering map, the quotient by~$g'$ is a positively twisted $\RR$-covered Anosov flow which is obtained from $\phi$ by Fried-Goodman surgeries along $\partial \Section$.

\end{proof}

%% file: SectionMarty/SectionMarty.tex
In this section, we prove Theorem~\ref{theorem:rcoveredcondition} using partial sections. We use two ideas. The first idea is the following. Given two partial section $S_1,S_2$ (assume first with disjoint boundary), the linking number between $S_1$ and $S_2$ (in the usual sens) is the algebraic intersection of $S_1$ with $\partial S_2$, or equivalently of $\partial S_1$ and $S_2$. Denote by~$\algcap$ the algebraic intersection. When $S_1$ has only positive boundary components, the algebraic intersection $\partial S_1\algcap S_2$ is positive (by our orientation convention) and so $S_1\algcap\partial  S_2$ is also positive. It implies that $S_2$ cannot have only negative boundary components. Therefore positive Birkhoff sections and negative partial sections cannot coexist. When the partial sections have some common boundary components, we add to the linking number the difference of framing along their common boundary. The generalize linking quantity, written in Lemma~\ref{Tlemma:LinkingNumberLikeFormula}, is then used to prove Theorem~\ref{Ttheorem:BirkhoffSectionRestriction}. 

For Anosov flows, it follows from works of Barbot and Fenley that there exists a large family of partial sections with exactly two boundary components which are additionally of the same sign (see Section \ref{Tsection:TwisteAnosovFlow}). Roughly speaking, given a lozenge $L$ in the orbit space, invariant by the action of an element $g\in\pi_1(M)\setminus\{1\}$, Barbot used a fundamental domain of the action $g$ on $L$ to build an immersed annulus, transverse to the flow in the interior and bounded by two periodic orbits. A technique introduced by Fried allows one to desingularize the annulus into a partial section. Additionally the sign of its boundary is the same as the one of the lozenge $L$. For $\RR$-covered positively twisted Anosov flows, every lozenge based on a periodic point (which represent a dense subset of the orbit space) is let invariant by a non-trivial element $g$, and all lozenges are positive. Then using a well-chosen family of lozenges and the corresponding positive partial sections, we can build by union and desingularization a positive Birkhoff section. For non-$\RR$-covered Anosov flows, it follows from Fenley's work that there exists invariant lozenges both positive and negative, and so positive and negative partial sections of the flow. This facts plus the discussion just above prevent the existence of a Birkhoff section whose boundary components are all of the same sign. It gives a proof of the two directions in Theorem~\ref{theorem:rcoveredcondition}.



\subsection{Partial sections and twisted flows}\label{Tsection:GeneralFlow}

We describe a restriction in the co-existence of some partial section for any flow in a oriented closed 3-dimensional manifold. 

\paragraph*{Partial section of a general flow.}

In this subsection only,~$\phi$ in a smooth flow non-necessarily Anosov. The ambient manifold~$M$ is still supposed closed, oriented and of dimension 3. Recall the definition of multiplicity given in Section~\ref{Section:PartialAndBirkhoffSection}, which extends for general flows. Given a partial section~$S$ and a parametrization~$\iota\colon\wh S\to S\subset M$, one boundary component~$\gamma^*\subset\wh S$ immersed inside a closed orbit~$\gamma$, the \emph{multiplicity} of the boundary component~$\gamma^*$ is the algebraic degree of the immersion~$\gamma^*\to\gamma$. Here~$\gamma^*$ is oriented according to the orientation on the interior of~$S$, and~$\gamma$ is oriented by the flow so that the multiplicity is positive if and only if the map~$\iota\colon\gamma^*\to\gamma$ preserves the orientation. Additionally~$\gamma^*$ is called positive (resp. negative) if its multiplicity is positive (resp. negative).

\begin{remark}
    Consider on~$\gamma^*$ the two orientations given by the flow and by the orientation on the interior of~$S$. Then~$\gamma^*$ is positive if and only if the two orientations are equal. This allows us to determine the sign of a boundary component by using only an arc of the boundary component, without computing the multiplicity.
\end{remark}

For non-Anosov flow, we cannot use the stable foliation as a reference to measure the linking number of a partial section. Hence the previous definition of the linking number does not generalize. Instead, we consider the linking number between two partial sections at a common boundary component.
Denote by~$\gamma$ a closed orbit of the flow,~$M_\gamma$ the manifold obtained by blowing up~$M$ along~$\gamma$, and by~$\pi\colon\TT_\gamma\to\gamma$ the normal bundle to~$\gamma$ inside~$M$, so that~$\partial M_\gamma=\TT_\gamma$. We fix the orientation of~$\TT_\gamma$ such that a vector normal to~$\TT_\gamma$ and going inside~$M_\gamma$ is positively transverse to~$\TT_\gamma$. 

Denote by~$\Section_1$ and~$\Section_2$ two partial sections whose boundaries contain~$\gamma$. We denote by~$\partial_\gamma \Section_i$ the union of boundary components of~$\Section$ which are immersed in~$\gamma$, and we orient~$\partial_\gamma \Section_i$ accordingly to the orientation of the interior of~$\Section_i$. The differential of the immersion~$\Section_i\to M$ induces a map~$f_i\colon\partial_\gamma \Section_i\to \TT_\gamma$. The \emph{linking number} between~$\Section_1$ and~$\Section_2$, denoted by~$\link_\gamma(\Section_1,\Section_2)$, is defined as the algebraic intersection between the image of~$f_1$ and~$f_2$, that is~$\link_\gamma(\Section_1,\Section_2)=[f_1]\algcap[f_2]$ where the algebraic intersection in taken inside~$H_1(\TT_\gamma,\ZZ)$. For this notion of linking number, it is relevant to consider all the boundary components of~$\Section_i$ immersed inside~$\gamma$. Notice that if we change the order of~$\Section_1$ and~$\Section_2$, the sign of the linking number changes. Additionally, we denote by~$\link(\Section_1,\Section_2)$ the sum of the linking numbers of~$\Section_1$ and~$\Section_2$ along all their common boundary components. 

\vspace{\baselineskip}

\paragraph{Sign of the linking number.}

The sign of the linking number of two partial sections can be determined when these partial sections have opposite signs on their boundary. We detail this sign below, and study the particular case when one of the partial sections is in fact a Birkhoff section.

Denote by~$\Section$ a partial section and by~$\gamma^*$ a boundary component of~$\Section$ immersed inside a closed orbit~$\gamma$. We use the notations~$\TT_\gamma$,~$M_\gamma$ and~$\pi_\gamma$ defined above. We view the boundary of~$\Section$ inside~$\TT_\gamma$ as the graph of a multivalued section of the circle bundle~$\pi_\gamma\colon\TT_\gamma\to\gamma$. Additionally, the flow~$\flow$ is smooth so it induces a smooth flow~$\psi$ inside~$\TT_\gamma$, whose orbits are positively transverse to each fiber of~$\pi_\gamma$. We use the orbits of~$\psi$ to speak about increasing and decreasing section of the bundle~$\TT_\gamma\to\gamma$. More precisely given a multivalued section~$\delta$ the orientation of~$\gamma$ by the flow induces an orientation on the graph of~$\delta$. Then an arc of this graph is called \emph{increasing} if its graph is positively transverse to the flow~$\psi$, \emph{decreasing} if it is negatively transverse, and \emph{constant} if its graph is an orbit arc. Notice that even when~$\Section$ is a Birkhoff section, the graph of~$\gamma^*$ may have constant arcs which are not transverse to the flow~$\psi$. To be more precise, an arc~$\delta$ is increasing if given~$x\in\delta$,~$u_\delta\in T_x\delta$ a non-zero tangent vector inducing the orientation on~$\delta$,~$u_\psi=\left(\frac{\partial \psi_t(x)}{\partial t}\right)_{|t=0}$ and~$v_{out}$ a normal vector to~$\TT_\gamma$ going outside~$M_\gamma$, the basis~$(u_\delta,u_\psi,v_{out})$ of~$T_xM_\gamma$ is positive.

\begin{lemma}\label{Tlemma:GraphMonoton}
    The graph of the multivalued section of~$\TT_\gamma\to\gamma$ induced by a boundary component of~$\Section$ is monotone. 
\end{lemma}

\begin{proof}
    The graph of~$\delta$ is smooth so if it is not monotone, it admits two points on which it is positively and negatively transverse to the flow~$\psi$. Then close to that two points inside the blowing up manifold~$M_\gamma$, there exists two other points in the interior of~$\Section$ whose orbit arcs of the flow induce two opposite co-orientations, which is not possible.
\end{proof}

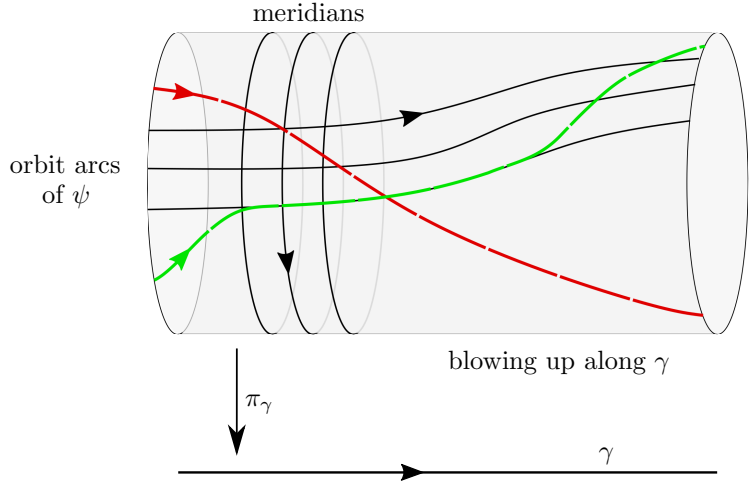
\begin{figure}
    \begin{center}
        \begin{picture}(100,65)(0,0)
        \put(80,3){$\gamma$}
        \put(60,15){blowing up along~$\gamma$}
        \put(34,61.3){meridians}
        \put(1.7,41){orbit arcs}
        \put(6,37){of~$\psi$}
        \put(33,10){$\pi_\gamma$}
        \put(20,0){\includegraphics[width=8cm]{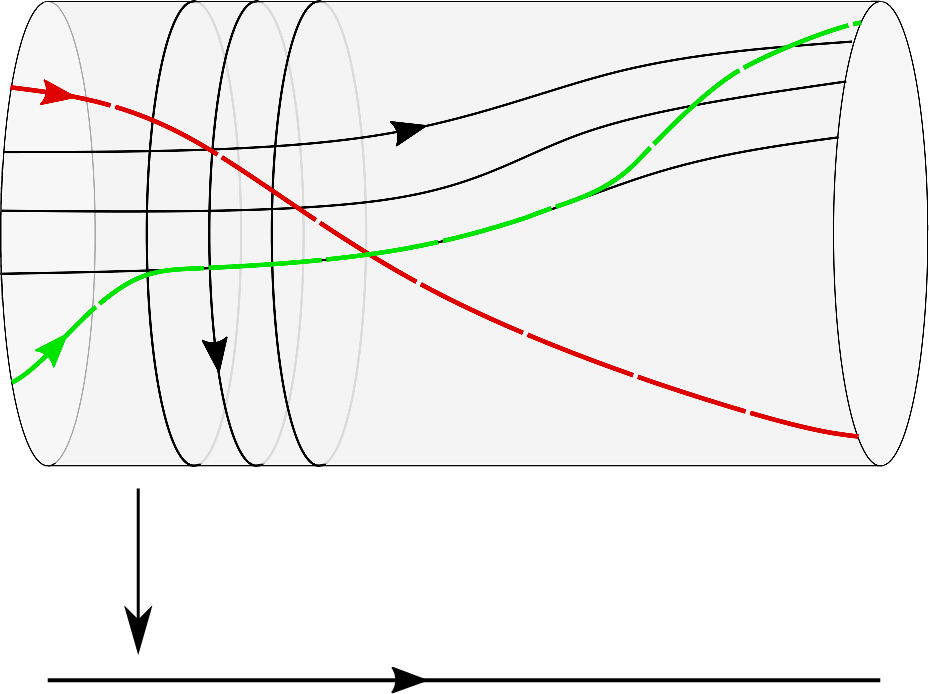}}
        \end{picture}
    \end{center}
    \caption{Boundary of two partial sections (in dashed-line) in the blowing up manifold of a closed orbit~$\gamma$. They are oriented relatively to the orientation of a meridian. The dark red curve corresponds to a negative boundary component, the light green to a positive boundary component.}
    \label{Tfigure:MonotoneCurves}
\end{figure}

\begin{lemma}\label{Tlemma:GraphMonotonSign}
    If~$\delta$ is induced by a positive (resp. negative) boundary component of~$\Section$, then its graph is non-decreasing (resp. non-increasing) relatively to the orbits of~$\psi$.
\end{lemma}

When the flow is Anosov, along a boundary component of linking number zero, the graph of~$\delta$ may be a closed orbit of the flow~$\psi$, that is one of the curves in~$\TT_\gamma$ induced by the stable and unstable directions. In this case the graph is constant.

\begin{proof}
    Suppose that~$\delta$ is induced by a positive boundary component and that there exists an arc~$\delta'\subset\delta$ which is increasing. Consider the orientation~$o$ of~$\delta'$ induced by the orientation of the interior of~$\Section$. Because the interior of~$\Section$ intersects positively the orbits of~$\flow$, the graph of~$\delta'$ oriented by~$o$ intersects negatively the orbit arcs of~$\psi$. Since the arc~$\delta'$ is increasing its algebraic intersection with an orbit arc of~$\psi$ has the same sign as its algebraic intersection with a fiber of~$\pi_\gamma$. So when oriented by~$o$ the graph of~$\delta'$ intersects negatively the fibers of~$\TT_\gamma\to\gamma$. Hence~$\pi_\gamma$ sends the orientation~$o$ to the orientation on~$\gamma$ opposite to the one induced by the flow. Hence the multiplicity of the boundary component of~$S$ corresponding to~$\delta'$ is negative, which contradict the hypothesis that~$\delta$ is induced by a positive boundary component. So the graph of~$\delta$ is non-decreasing.
\end{proof}

\begin{lemma}\label{Tlemma:GraphForBirkhoffSection}
    If~$\Section$ is a Birkhoff section, then the graph of~$\delta$ intersects every orbit of the flow~$\psi$.
\end{lemma}

\begin{proof}
    Denote by~$S^*$ the surface induced by~$\Section$ inside the blowing up manifold~$M_\gamma$, and let~$x$ be a point in~$\TT_\gamma$. Take~$(x_n)_{n\in\NN}$ be a family of points in~$M_\gamma\setminus T_\gamma$ that converge to~$x$. Since~$\Section$ is a Birkhoff section, there exists~$T>0$ and for every~$n\in\NN$ there exists~$0\leq t_n\leq T$ so that~$\psi_{t_n}(x_n)$ is inside~$S^*$. Up to taking a subfamily of~$(x_n)_n$, we can suppose that~$(t_n)_n$ converges to some~$t\in[0,T]$, so that~$\psi_t(x)\in S^*\cap\TT_\gamma$.
\end{proof}

For general flows, we need to understand the linking number with a Birkhoff section. The following lemma is a generalization of Lemma~\ref{Lemma:NonPositiveLink} for general flows. 

\begin{lemma}\label{Tlemma:LinkinNumberAlongCommonBoundaryComponent}
    Let~$\Section_1$ be a positive partial section and~$\Section_2$ be a negative partial section which have a common boundary component denoted by~$\gamma$. Then we have~$\link_\gamma(\Section_1,\Section_2)\geq 0$. 
    
    If~$\Section_1$ or~$\Section_2$ is a Birkhoff section, then we have~$\link_\gamma(\Section_1,\Section_2)> 0$.
\end{lemma}

\begin{proof}
    Denote by~$\delta_i$ the multivalued section of~$\TT_\gamma\to\gamma$ induced by the boundary of~$\Section_i$ along the orbit~$\gamma$. We oriented~$\delta_i$ via the direction of the flow, such that for any fiber~$m$ of~$\TT_\gamma\to\gamma$ we have~$[\delta_i]\algcap [m]<0$. We also denote~$\delta_i'$ a copy~$\delta_i$, but oriented using the orientation of~$\Section_i$. The section~$\Section_1$ is positive and the section~$\Section_2$ is negative, so the homology elements satisfies~$[\delta_1']=[\delta_1]$ and~$[\delta_2']=-[\delta_2]$

    According to Lemma~\ref{Tlemma:GraphMonotonSign}, the graph of~$\delta_1$ is non-decreasing and the graph of~$\delta_2$ is non-increasing, so the algebraic intersection~$[\delta_1]\algcap[\delta_2]$ is non-positive. Hence the algebraic intersection~$[\delta_1']\algcap[\delta_2']$ is non-negative and~$\link_\gamma(\Section_1,\Section_2)=[\delta_1']\algcap[\delta_2']$.

    We now additionally assume that~$\Section_1$ is a Birkhoff section, and suppose that the linking number is zero. Then according to the previous lemma, the graph of~$\delta_1$ is not a single orbit of the flow~$\psi$. Hence we can homotope~$\delta_1$ to a multivalued section~$\wt\delta_1$ whose graph is decreasing and which also intersects every orbit of the flow~$\psi$. Since the graph of~$\delta_2$ is non-decreasing, the intersection of the graphs of~$\wt\delta_1$ and~$\delta_2$ are transverse and always intersect positively. Additionally, the graph of~$\wt\delta_1$ intersects every orbit of the flow~$\psi$, so the graphs of~$\wt\delta_1$ and~$\delta_2$ intersects at less once. Indeed, if~$x_2$ is a point in the graph of~$\delta_2$ and~$x_1=\psi_t(x)$ is in the graph of~$\wt\delta_1$ with~$t>0$, because~$\wt\delta_1$ is decreasing and~$\delta_2$ is increasing, their graphs intersect in the flow box contained between the fibers of~$x_1$ and~$x_2$. So the linking number is not zero.

    The case when~$\Section_2$ is a Birkhoff section is obtained by reversing the role of~$\Section_1$ and~$\Section_2$ in the previous argument.
\end{proof}

\begin{lemma}
    \label{Lemma:MultLinkWellDefinedGeneralCase}
     Let~$\phi$ be a flow in a oriented closed 3-dimensional manifold~$M$. For a Birkhoff section~$\Section$ of~$\phi$ and a closed orbit~$\gamma$ in~$\partial\Section$, the multiplicity~$\mult(\gamma^*)$ does not depend on the choice of a boundary component~$\gamma^*$ immersed in~$\gamma$.
\end{lemma}

\begin{proof}
    Suppose there are two boundary components~$\gamma^*_1$ and~$\gamma^*_2$ of~$\Section$ immersed inside~$\gamma$. Consider the immersion~$i\colon\Section\to M_{\gamma}$ from~$\Section$ to the blowing up manifold~$M_{\gamma}$. Denote by~$\TT_\gamma$ the torus boundary of~$M_\gamma$. Consider a tubular neighborhood~$U$ of~$\gamma$ such that its neighborhood is transverse to~$\Section$. Then~$\partial U\cap\Section$ is a union of finitely many curves, which are disjoint since~$\Section$ is embedded in it interior. Hence~$i(\gamma^*_1)$ and~$i(\gamma^*_2)$ are homotopic inside~$\TT_\gamma$ as non-oriented curves. 
    
    We oriente~$\gamma^*_k$ using the orientation on~$\Int(\Section)$, which is induced by the orientation of~$M$ and the coorientation by the flow. If~$i(\gamma^*_1)$ and~$i(\gamma^*_2)$ where anti-homotopic as oriented curves, then using the arguments from the previous lemma, we would have~$i(\gamma^*_1)\algcap i(\gamma^*_2)\neq 0$, which is false since the curves are homotopic. Hence~$i(\gamma^*_1)$ and~$i(\gamma^*_2)$ where homotopic as oriented curves, and they have the same multiplicities. 
\end{proof}

\paragraph*{Linking equation for two partial sections.}

Let~$\Section_1$ and~$\Section_1$ be two partial sections of the flow~$\flow$. The intersection between these partial sections can be use to express a restriction for the co-existence of some partial section. We prove a linking type lemma and prove Theorem~\ref{Ttheorem:BirkhoffSectionRestriction}. For two partial sections~$\Section_1$ and~$\Section_2$, we denote by~$\partial \Section_2\algcap \Intt{\Section}_1$ the algebraic intersection between the boundary of~$\Section_2$ taken with multiplicity and the interior of~$\Section_1$.
  
\begin{lemma}\label{Tlemma:LinkingNumberLikeFormula}
    Let~$\flow$ be a flow of an oriented 3-dimensional manifold~$M$, and~$\Section_1$ and~$\Section_2$ be two partial sections of~$\flow$. Then we have
    $$(\partial \Section_1\algcap \Intt{\Section}_2) - (\partial \Section_2\algcap\Intt{\Section}_1)+\link(\Section_1,\Section_2)=0$$
\end{lemma}

The lemma remains true for non-transverse surfaces, as long as two boundary components are either disjoint or equal.

\begin{proof}
    
    We denote by~$M_{\partial \Section_1 \cup\partial \Section_2}$ the blow up of~$M$ along~$\partial \Section_1 \cup\partial \Section_2$. We denote by~$\Section_1^*$ and~$\Section_2^*$ the surfaces in~$M_{\partial \Section_1 \cup\partial \Section_2}$ induced by~$\Section_1$ and~$\Section_2$. Up to a small isotopy of~$\Section_1^*$ and~$\Section_2^*$ relatively to the boundary of~$M_{\partial \Section_1 \cup\partial \Section_2}$, we can suppose that~$\Section_1^*$ and~$\Section_2^*$ are embedded and intersect transversally. 
    
    Then~$\partial \Section_1\algcap \Intt{\Section}_2$ is equal to the algebraic intersection~$\partial\Section_1^*\algcap_{\partial \Section_1 \setminus\partial \Section_2}\partial\Section_2^*$ between the curves~$\partial\Section_1^*$ and~$\partial\Section_2^*$, where the algebraic intersection is taken inside the union of the torus boundary components of~$M_{\partial \Section_1 \cup\partial \Section_2}$ obtained by blowing up~$\partial \Section_1 \setminus\partial \Section_2$. 
    Similarly,~$\partial \Section_2\algcap \Intt{\Section}_1=-\partial\Section_1^*\algcap_{\partial \Section_2 \setminus\partial \Section_1}\partial\Section_2^*$ and~$\link(\Section_1,\Section_2)=\partial\Section_1^*\algcap_{\partial \Section_1 \cap\partial \Section_2}\partial\Section_2^*$. Hence~$(\partial \Section_1\algcap \Intt{\Section}_2) - (\partial \Section_2\algcap\Intt{\Section}_1)+\link(\Section_1,\Section_2)=\partial\Section_1^*\algcap\partial\Section_2^*$ so it is enough to have~$\partial\Section_1^*\algcap\partial\Section_2^*=0$, which is true since~$\partial\Section_1^*$ and~$\partial\Section_2^*$ are the boundaries of two surfaces inside~$M_{\partial \Section_1 \cup\partial \Section_2}$

\end{proof}

The previous Lemma is enough to prevent two positive and negative partial sections to intersect on their boundaries.

\begin{proof}[Proof of Theorem~\ref{Ttheorem:BirkhoffSectionRestriction}]
Denote by~$\Section_1$ a positive Birkhoff section and suppose that there exists a negative partial section~$\Section_2$. Then by Lemma~\ref{Tlemma:LinkingNumberLikeFormula}, we have
$$(\partial \Section_1\algcap \Intt{\Section}_2) - (\partial \Section_2\algcap\Intt{\Section}_1)+\link(\Section_1,\Section_2)=0$$
Since the partial section~$\Section_2$ has only negative boundary components,~$\partial\Section_2$ is a union of closed orbit with the orientation opposite to the flow. And the orientation on the interior~$\Section_1$ is chosen such that the intersection with an orbit arc oriented by the flow is positive, so the algebraic intersection~$\partial \Section_2\algcap\Intt \Section_1$ is non-positive. Similarly the algebraic intersection~$\partial \Section_1\algcap\Intt \Section_2$ is non-negative. According to Lemma~\ref{Tlemma:LinkinNumberAlongCommonBoundaryComponent}, the linking number~$\link(\Section_1, \Section_2)$ is non-negative. So in the sum the three terms are non-negative and add up to zero, so there are all zero.

Since all algebraic intersection of~$\partial \Section_2$ and~$\Intt \Section_1$ have the same sign, the geometric intersection~$\partial \Section_2\cap \Intt \Section_1$ is empty. But~$\Section$ is a Birkhoff section which intersects every orbit, so the boundary components of~$\Section_2$ are also boundary component of~$\Section_1$. For each boundary component~$\gamma$ of the partial section~$\Section_2$ (which exist by definition of a negative partial section), the linking number along~$\gamma$ satisfies~$\link_\gamma(\Section_1,\Section_2)>0$ according to Lemma~\ref{Tlemma:LinkinNumberAlongCommonBoundaryComponent}. So the sum of the linking number for all boundary components of~$\Section_2$ is non-zero, which contradicts the previous paragraph. Hence if a flow admits a positive Birkhoff section, it does not admit a negative partial section.

The second statement can be proved the same way.

\end{proof}
    
\subsection{Application to Anosov flows}\label{Tsection:TwisteAnosovFlow}

When the flow is Anosov, we can use the orbit space of the flow and some well-chosen lozenges, to construct partial sections with specific signs on their boundary. Then we can relate the existence of these partial section with the existence of a positive Birkhoff section.


\paragraph{Fried-desingularisation.}

We can construct a Birkhoff section by taking the union enough immersed partial sections and desingularizing the union. This operation is called \emph{Fried-desingularisation}. We describe this process, then construct good partial sections. Let~$S$ be an oriented smooth compact surface and~$f\colon S\to M$ be a~$\mathcal{C}^2$ immersion such that~$f$ is transverse to the flow on the interior of~$S$ and tangent to the flow on the boundary of~$S$. We denote by~$\Gamma = f(\partial S)$ which is a union of closed orbits, by~$M_{\Gamma}$ the blowing up manifold of~$M$ along~$\Gamma$ and by~$\Section^*$ is the blowing of~$S$ along~$\Int(S)\cap f^{-1}(\Gamma)$. The map~$f$ lifts to a~$\mathcal{C}^1$ map~$f^*\colon \Section^*\to M_{f(\partial S)}$. We say that~$f$ is \emph{$\delta$-strong} if the restriction~$f^*_{\partial \Section^*}$ is transverse to the flow on~$\partial M_\Gamma$. This is a technical condition which makes future argument easier. It can be replaced by finer hypothesis or weaker conclusions for the following lemma.

\begin{lemma}[Fried-desingularisation~\cite{Fri}]\label{Tlemma:Fried-desingularisation}
    We suppose~$f\colon S\to M$~$\delta$-strong. Then for all~$\epsilon>0$, there exists a~$\delta$-strong partial section~$S'$ in an~$\epsilon$-neighborhood of~$f(S)$, homologous to~$f(S)$ relatively to~$f(\partial S)$ and a bounded and piecewise continuous function~$T\colon \Int(S)\to\RR$ such that~$(x\in S)\mapsto \phi_{T(x)}(f(x))$ is a bijection from~$\Int(S)$ to~$\Int(S')$.  
\end{lemma}

Notice that the surface~$S$ is not necessarily connected. Outside small neighborhood of~$\Gamma$, the surface~$S'$ can be taken in an~$\epsilon$-flow box of~$f(S)$. Close to~$\Gamma$, it is not necessarily the case since the desingularisation may require doing bounded but non-small homotopies along the flow.

The following proof is due to Fried, but is made more precise using the~$\delta$-strong hypothesis. 

\begin{proof} 
    Denote by~$f^*\colon \Section_1^*\to M_\Gamma$ the lift of~$f$ defined above the lemma. For a function~$T\colon \Section_1^*\to\RR$, we define the map~$f^*_T\colon \Section_1^*\to M_\gamma$ by~$f^*_T(x)=\phi_{T(x)}(f^*(x))$. Because~$f\colon S\to M$ is~$\delta$-strong, we can find a function~$T_1\colon \Section_1^*\to [\epsilon/3,\epsilon/3]$ such that the image of~$f^*_{T_1}$ has only transverse self-intersections. 
    
    We can cut, paste and smooth~$\Section_1^*$ along the self-intersection curves of~$f^*_{T_1}(S)$ to obtain a smooth compact surface~$\Section_2^*$, a smooth embedding~$g^*\colon \Section_2^*\to M_\gamma$ whose image is in a~$\frac{2}{3}\epsilon$-flow box of~$f^*(\Section_1^*)$ and homologous to~$f^*(\Section_1^*)$ relatively to~$\partial M_\Gamma$, and a piecewise continuous function~$T_2\colon \Section_1^*\to [-\epsilon,\epsilon]$ such that~$f^*_{T_2}\colon \Section_1^*\to M_\Gamma$ is a bijection onto~$g^*(\Section_2^*)$. We can do the cut-and-paste operation such that~$g^*(\Section_2^*)$ is everywhere transverse to the flow. 
    
    Then we pull tight the curves in~$g^*_{T_2} (\Section_2^*)\cap \partial M_\gamma$ relatively to the fibers of~$M_\Gamma\to M$. More precisely we can find a smooth map~$T_3\colon \Section_2^*\to\RR$ such that 
    \begin{itemize}
        \item~$g^*_{T_3}$ has its image inside an~$\epsilon$-neighborhood of~$f^*_{T_2} (\Section_1^*)\cup \partial M_\gamma$,
        \item for every boundary component~$\delta$ of~$\Section_2^*$, the restriction of~$g^*_{T_3}$ to~$\delta$ either has its image is a fiber of~$\partial M_\Gamma\to\Gamma$ or is transverse to the foliation of~$\partial M_\Gamma$ by these fibers,
        \item~$g^*_{T_3}$ is an embedding.
    \end{itemize}
    We denote by~$\pi\colon M_\Gamma\to M$ the blow down projection, and by~$S'$ the blow down surface of~$S'$ along all boundary components~$\delta$ whose images by~$g^*_{T_3}$ are fibers of~$\partial M_\Gamma\to\Gamma$. Then the map~$\pi\circ g^*_{T_3}$ induces a continuous map~$g\colon S' \to M$ which is smooth outside the points induced by the previous blowing down operation. We can choose~$T_3$ such that~$g$ is everywhere smooth. Then the map~$g$ is a parametrization of a partial section. By construction, it is relatively homologous to~$f(S)$, in an~$\epsilon$-neighborhood of~$f(S)$ and it comes with a function~$T$ induced by~$T_1$,~$T_2$ and~$T_3$, which satisfies the lemma.
\end{proof}

\paragraph{From Lozenge to partial section.}

We are interested in immersed Birkhoff annulus, which correspond to `immersed partial section' with the topology of an annulus. A large family of them can be defined using the orbit space and thanks to the Fried-desingularisation, they can be used to produce good partial sections. An annulus immersed inside~$M$ is called an \emph{immersed Birkhoff annulus} if its interior is transverse to the flow and its two boundary components are closed orbits of the flow. T.Barbot studied the relation between immersed Birkhoff annuli of an Anosov flow and the orbit space of the flow. To express that relation we denote by~$A$ the Birkhoff annulus. Denote by~$\wt A$ a lift of~$\Int A$ to the universal covering~$\wt M$ of~$M$. That is~$\wt A$ is a connected component of~$P^{-1}(\Int A)$. The subset~$Q(\wt A)\subset\orbitspace$ obtained is called the \emph{trace} of~$A$, it is well-defined up to the action of~$\pi_1(M)$ on the orbit space. The trace can be defined more generally for any partial section.

\begin{theorem}[Barbot~\cite{Ba2} + $\delta$-strong]\label{Ttheorem:LozengeToAnnulus}
    Let~$L$ be a lozenge in the orbits space~$\orbitspace$ bounded by two periodic orbits. Then there exists a~$\delta$-strong immersed Birkhoff annulus~$A$ whose trace in~$\orbitspace$ is the union of~$L$ and its two corners. That is there exist an immersion~$f\colon \RR/\ZZ\times[0,1]\to M$ whose image is~$A$, and an embedding~$\wt f\colon \RR\times[0,1]\to\wt M$ lifting~$f$, such that~$Q\circ \wt f(\RR\times(0,1))=L$.
\end{theorem}

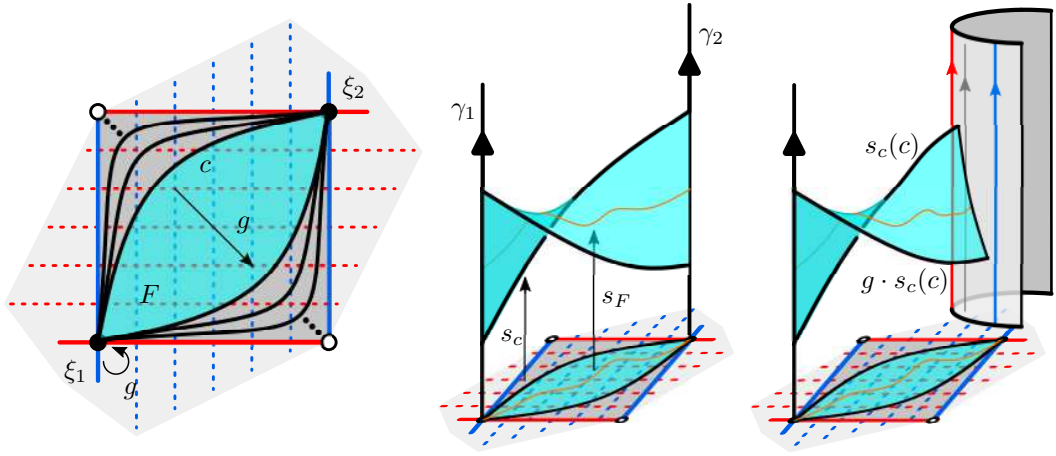
\begin{figure}
    \begin{center}
        \begin{picture}(140,60)(0,0)
            \put(-10,0){\includegraphics[height=6cm]{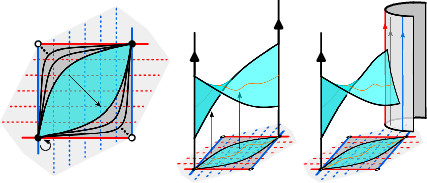}}
            \put(-2,10){$\xi_1$}
            \put(35,47.8){$\xi_2$}
            \put(6,8){$g$}
            \put(21,30){$g$}
            \put(16,37.5){$c$}
            \put(8,20){$F$}

            \put(49.5,45){$\gamma_1$}
            \put(82.4,55){$\gamma_2$}
            \put(56,15){$s_c$}
            \put(69.5,20){$s_F$}

            \put(104.5,40){$s_c(c)$}
            \put(104,22){$g\cdot s_c(c)$}
        \end{picture}
    \end{center}
    \caption{Construction of a Birkhoff annulus from a lozenge invariant by the action of an element $g\in\pi_1(M)\setminus\{1\}$.}
    \label{Tfigure:BirkhoffAnnulus}
\end{figure}

The following proof is adapted from the original proof by Barbot to obtain the additional~$\delta$-strong property. See Figure \ref{Tfigure:BirkhoffAnnulus} for an illustration.

\begin{proof}
    We construct the annulus~$A$ by using a fundamental domain of~$L$. 
    Up to a double cover, we can suppose that the stable and unstable foliations are orientable. By hypothesis, the corners~$\xi_1,\xi_2$ of~$L$ correspond to two closed orbits~$\gamma_1$ and~$\gamma_2$ of the flow~$\flow$. We fix a point~$x\in\gamma_1$ and~$g=[\gamma_1]^k\in\pi_1(M,x)$ to be the first power of the homotopy class of~$\gamma_1$ which preserves the orientation of the stable and unstable leaves of~$\gamma_1$. The action of~$g$ lets invariant the half stable and half unstable leaves of~$\xi_1$ since the action of~$g$ preserves the orientation of the stable and unstable leaves of~$\gamma_1$. Hence~$g$ preserves the lozenge~$L$ and its second corner~$\xi_2$.
    
    We consider a curve~$c$ inside~$L$ joining the two corners~$\xi_1,\xi_2$ of the lozenge~$L$, and transverse to the two stable and unstable foliations. Because~$g$ expands the leaf~$\ofoliation^s(\xi_1)$ and contracts the leaf~$\ofoliation^u(\xi_1)$ (see Proposition~\ref{Aprop:g action})~$g*c$ and~$c$ only intersect at the corners of~$L$, and they delimit inside~$L$ a set~$F\subset L$. It is a fundamental domain for the action of~$g$ on~$L$. More precisely the following holds:
    \begin{itemize}
        \item $\bigcup_{n\in\ZZ}(g^n*F)=L$,
        \item $(g^n*F)\cap(g^m*F)=\emptyset$ when~$|n-m|\geq 2$,
        \item $(g^n*F)\cap(g^{n-1}*F)=g^n*c$
    \end{itemize}

    We take the fundamental domain~$F$ closed inside~$L$ so that the interior of~$c$ and~$g* c$ are inside~$F$. 
    
    Choose a lift $\wt s\colon L\subset\orbitspace\to\wt M$ equivariant by the action of $g$. It can be built by choosing its value on $c$, pushing it by $g$ on $g*c$, interpolating it inside $F$ and pushing it again by $g^n$ for all $n\in\ZZ$. We can additionally choose it smooth.

    
    Denote by~$P\colon \wt M\to M$ the covering map and by~$F'=F/(x\simeq g* x\text{ for }x\in c)$ the quotient annulus. Then the map $\wt s$ is equivariant by $g$ so it induces a well-defined map~$f\colon F'\to M$. It is the parametrization of a smooth and immersed open annulus transverse to the flow. To obtain a~$\delta$-strong Birkhoff annulus, we carefully chose the curve~$c$ and the section~$s_F$ on neighborhoods of the corners of~$L$. 
    
    Now assume that $c$ has be chosen to be the restriction to $L$ of a smooth curve containing $\xi_1$ and$ \xi_2$, and transverse to the stable and unstable foliations at these points. Then after lifting $c$ to $\wt M$, pushing it to $M$, and blowing up $\gamma_1$, the curve corresponding to $c$ converges on the torus boundary (or the Möbius strip boundary) corresponding to $\gamma_1$ in a point which is not periodic. Thus we can choose the maps $s$ and $f$ so that the boundary of $f(F')$ is transverse to the flow on the blown up orbit of $\gamma_1$. And similarly for $\gamma_2$.
\end{proof}
    
We extend the sign of a boundary component of an immersed Birkhoff annulus with the same definition. The type of the lozenge~$L$ determines the sign of the boundary components of the induced immersed Birkhoff annulus. 

\begin{lemma}\label{Tlemma:SignOfBirkhoffAnnulus}
    Let~$L$ be a lozenge in the orbit space of~$\flow$, and~$\Section$ be an immersed Birkhoff annulus for~$\flow$ with trace~$L$, constructed in Theorem~\ref{Ttheorem:LozengeToAnnulus}. If~$L$ is a lozenge of type~$(++)$, then the two boundary components of~$\Section$ are positive. Otherwise the two boundary components of~$L$ are negative.
\end{lemma}

\begin{proof}
    We assume that~$L$ is a lozenge of type~$(++)$, the second case need the same proof with adapted signs. Denote by~$f\colon (\RR/\ZZ)\times[0,1]$ a parametrization of the immersed Birkhoff annulus~$\Section$ with trace~$L$. We denote by~$\gamma$ the image of~$f_{|(\RR/\ZZ)\times\{0\}}$, which is a closed orbit of the flow. We chose~$f$ such that~$f_{|(\RR/\ZZ)\times\{0\}}$ send the direct orientation on~$\RR/\ZZ$ to the orientation on~$\gamma$ induced by the flow. The map~$f$ is used to relate the orientation on~$\Section$ with the orientation on the orbit space. 
    
    We lift the map~$f$ to a map~$\wt f\colon \RR\times[0,1]\to\wt M$. According to Theorem~\ref{Ttheorem:LozengeToAnnulus}, we can choose the lift~$\wt f$ such that~$Q\circ\wt f$ sends~$\RR\times(0,1)$ to~$L$. Then the two boundary components~$\RR\times\{0\}$ and~$\RR\times\{1\}$ are sent to the two corners of~$L$. We suppose that~$\RR\times\{0\}$ is sent to the lower left corner of the lozenge~$L$, that we denote by~$\xi_1$.  Denote by~$x_0$ a point inside~$\gamma$ and by~$g\in\pi_1(M,x_0)$ the homotopy class of~$\gamma$, which satisfies~$g*\xi_1=\xi_1$. According to Proposition~\ref{Aprop:g action} the element~$g$ expands the stable leaf of~$\xi_1$ and contracts the unstable leaf of~$\xi_1$. Hence for any point~$x$ in the interior of~$L$, the point~$g*x$ is still inside~$L$ and is in the lower right quadrant to~$x$.
    
    Given a element~$a\in[0,1]$ the curve~$f_{|(\RR/\ZZ)\times\{a\}}$ is homotopic to~$m>0$ copies of~$\gamma$, where~$m$ is the degree of the covering~$f\colon (\RR/\ZZ)\times\{0\}\to\gamma$. So for any~$t\in\RR$ we have~$\wt f(t+1,a)=g^m\cdot \wt f(t,a)$. Hence~$Q\circ\wt f(t+1,a)=g^m * Q\circ\wt f(t,a)$ and this point is in the lower right quadrant of the point~$Q\circ\wt f(t,a)$. 
    
    \begin{claim*}
        The curve~$(t\in\RR\mapsto Q\circ\wt f(t,a)\in L)$ is an open curve which at infinity goes from the upper left corner of~$L$ to the lower right corner of~$L$.
    \end{claim*}

    To prove the claim, we use the coordinate system on~$L$ given by the homeomorphism~$q\colon L\to\ofoliation_+^s(c)\times\ofoliation_+^u(c)$ defined by~$q(\xi)=(q^s(\xi ), q^u(\xi ))=(\ofoliation_+^s(c)\cap\ofoliation^u(\xi ), \ofoliation_+^u(c)\cap\ofoliation^s(\xi ))$. Using that coordinate system, for any~$\xi~$ in the interior of~$L$, the sequences~$(n\mapsto q^s(g^n*\xi))_n$,~$(n\mapsto q^u(g^n*\xi))_n$ are monotone in the half leaves they respectively lie on. Since there is no points inside~$\ofoliation_+^s(c)$ and~$\ofoliation_+^u(c)$ invariant by~$g$, the above sequences must diverge to an infinity ends of their respective half leaves. The element~$g$ expands~$\ofoliation^s(c)$ and contracts~$\ofoliation^s(c)$, so the sequence~$(n\mapsto g^n* \xi )_{n\in\ZZ}$ diverges to the lower right corner of~$L$ when~$n\rightarrow+\infty$ and diverges to the upper left corner of~$L$ when~$n\rightarrow-\infty$. Applying this fact on a compact fundamental domain of~$\RR\to\RR/\ZZ$ finish the proof of the claim.

    We continue the proof of the lemma. For~$t\in\RR$ fixed, since the curve~$(a\mapsto f(a,t))$ goes from~$\gamma$ to the other boundary component of~$\Section$, the curve~$(a\mapsto Q\circ\wt f(a,t))$ goes from the lower left corner of~$L$ to its upper right corner. 
    
    Hence the algebraic intersection between the two curves~$(t\mapsto Q\circ\wt f(t,a_0))$ and~$(a\mapsto Q\circ\wt f(t_0,a))$ is~$+1$. The algebraic intersection between the two curves~$(t\mapsto (t,a_0))$ and~$(a\mapsto (t_0,a))$ is also~$+1$. Hence the map~$Q\circ\wt f$ send the orientation on~$\RR\times[0,1]$ to the orientation of the orbit space induced by the flow. The projection~$Q$ preserves the orientation of a transverse surface which is induced by the coorientation by the flow, so~$f$ itself is orientation-preserving. So the orientation of~$\Section$ induces on~$\gamma$ the same orientation that the flow, which implies that~$\gamma$ is a positive boundary component of~$\Section$. And by symmetry the second boundary component of~$\Section$ is also a positive boundary component of~$\Section$. 
\end{proof}


\vspace{\baselineskip}

\begin{proof}[Proof of Theorem~\ref{theorem:rcoveredcondition} implication ($\Rightarrow$)]
    Suppose that the flow~$\flow$ admits a positive Birkhoff section. According to Theorem~\ref{Ttheorem:BirkhoffSectionRestriction}, the flow~$\flow$ does not admit any negative partial section nor any Birkhoff section without boundary. Solodof proved that an Anosov flow admits a Birkhoff section without if and only if it is non-twisted~$\RR$-covered. Hence the flow~$\flow$ is not non-twisted~$\RR$-covered.
    
    Suppose that~$\flow$ is not positively twisted~$\RR$-covered. Then it is either negatively twisted~$\RR$-covered, or not~$\RR$-covered. In the first case, there exists a lozenge~$L\subset\orbitspace$ bounded by closed orbits and in the quadrant~$(+-)$ of one of its corners. According to Theorem~\ref{Ttheorem:LozengeToAnnulus}, there exists an immersed Birkhoff annulus~$A$ whose trace in~$\orbitspace$ is~$L$, and~$A$ has two negative boundary component according to Lemma~\ref{Tlemma:SignOfBirkhoffAnnulus}. Then according to Lemma~\ref{Tlemma:Fried-desingularisation}, a Fried-desingularisation of~$A$ induces a partial section~$\Section$ which is relatively homologous to~$A$. Then~$\Section$ the multiplicities of~$\Section$ along the two boundary components of~$A$ are equal to the multiplicities of~$A$, so~$\Section$ is a negative partial section which contradicts the previous paragraph.
    
    In the second case,~$\flow$ is not~$\RR$-covered and according to Theorem~\ref{Athm:adjacent lozenge} there exists a lozenge~$L$ of type~$(+-)$ bounded by two periodic orbits. For the same reason, this case is not possible. Hence~$\flow$ is positively twisted~$\RR$-covered.
\end{proof}

We now give the proof of the existence of positive Birkhoff sections for positively skewed Anosov flows. The proof has some similar arguments than Fried's original proof of the existence of Birkhoff section for transitive Anosov flows, but we do not actually use Fried results.

\begin{proof}[Proof of Theorem~\ref{theorem:rcoveredcondition} implication ($\Leftarrow$)]
    Suppose that the flow~$\flow$ is $\RR$-covered and positively twisted. We construct a finite number of immersed Birkhoff annuli with positive boundaries and whose interiors intersect every orbit. Then we use them to construct a positive Birkhoff section.

    Take~$K\subset\orbitspace$ a compact containing a lift of every orbit of~$\flow$. Since~$(M,\flow)$ is positively twisted~$\RR$-covered, for every~$\xi$ in~$K$ there exists a lozenge~$L_\xi\subset \orbitspace$ containing~$\xi$. Since~$\flow$ is~$\RR$-covered, it is transitive~\cite{Ba1}. Hence we can choose the lozenge $L_\xi$ to be bounded by closed orbits. 

    By compactness of~$K$, there exist a finite sequence~$\xi_1,\hdots,\xi_n$ inside~$\orbitspace$ such that all points inside~$K$ is in the interior of one of the lozenge~$L_{\xi_i}$. According to Theorem~\ref{Ttheorem:LozengeToAnnulus}, for every index~$i$ there exists an immersed Birkhoff annulus~$A_i\subset M$ of the flow whose trace is given by the lozenge~$L_{\xi_i}$. The flow is positively~$\RR$-covered, so the lozenge~$L_\xi$ is of type~$(++)$. Then according to Lemma~\ref{Tlemma:SignOfBirkhoffAnnulus} the boundary components of~$A_i$ are both positive. 
    
    We consider the partial section~$\Section$ obtained as the Fried-desingularisation~\cite{Fri} (see Lemma~\ref{Tlemma:Fried-desingularisation}) of the union of the immersed Birkhoff annulus~$A_i$. The surface~$\Section$ is homologous to~$\cup_i A_i$ relatively to the union of the boundary components of the~$A_i$, and all~$A_i$ have only positive boundary components, so~$\Section$ is a positive partial section. Notice that all orbits of the flow intersect the interior of one immersed Birkhoff annulus~$A_i$, so according to Lemma~\ref{Tlemma:Fried-desingularisation}, all orbits intersect~$\Section$. 
    
    We adapt a classical argument to prove that~$\Section$ is a Birkhoff section. 
    Suppose that for all~$T>0$, there exists a point~$q_T\in M$ so that~$\flow_{[0,2T]}(q_T)\cap \Section=\emptyset$. By compactness of~$M$,~$\flow_T(q_T)$ accumulates on some point~$q_\infty\in M$. If~$q_\infty$ is not in the boundary of~$\Section$, then the orbit of~$q_\infty$ intersects the interior of~$\Section$. That is there exists~$t>0$ such that~$\flow_t(q_\infty)\in\Intt \Section$. Then all orbits close to~$q_\infty$ intersect the interior of~$\Section$ in a time at most~$t+\epsilon$ for some~$\epsilon>0$, which contradicts~$\flow_{[-T,T]}(\flow_T(q_T))\cap \Section=\emptyset$ for all~$T>0$. If~$q_\infty\in \partial \Section$ the orbit of~$q_\infty$ intersects the interior of one of the annulus~$A_i$ for some index~$i$. So there exists a neighborhood~$U$ of~$q_\infty$ and a~$\epsilon>0$ and a time~$T_1>\epsilon$ such that for all point~$x\in U$ the orbit arc~$\phi_{[T_1-\epsilon,T_1+\epsilon]}(x)$ intersects the interior of~$A_i$. By Lemma~\ref{Tlemma:Fried-desingularisation}, there is a time~$T_2>0$ such that for all point~$y\in\cup_i A_i$, the orbit arc~$\phi_{[-T_2,T_2]}(y)$ intersects~$\Section$. Then for all points~$x\in U$, the orbit arc~$\phi_{[T_1-\epsilon-T_2,T_1+\epsilon+T_2]}(x)$ intersects~$\Section$, which contradict the definition of~$q_\infty$. Hence the surface~$\Section$ is a positive Birkhoff section.
\end{proof}

\begin{remark}
    The construction described above immediately yield a proof for the following statement. Given a finite family of periodic points $\xi_1\cdots\xi_{2n}$ in~$\orbitspace$, so that $\xi_{2i}$ and $\xi_{2i+1}$ are the corners of a common positive lozenge $L_i$, and so that every point $\xi\in\orbitspace$ lies in (the interior of) a copy $g\cdot L_i$ of a lozenge under the action of $g\in\pi_1(M)$, then set $\gamma_1\cdots\gamma_{2n}\in M$ of periodic orbits corresponding to $\xi_1\cdots\xi_{2i}$ bound a positive Birkhoff section.
\end{remark}


%% file: SectionAsaoka/SectionAsaoka.tex
In this section, we give the following theorem,
 which gives a second proof of the $\RR$-coveredness
 of Anosov flows which admit a positive Birkhoff section.

\begin{theorem} 
\label{Athm:A main}
Let $\flow$ be a topologically transitive Anosov flow
 on a 3-dimensional closed manifold $M$.
If $\flow$ is not $\RR$-covered, then
 any Birkhoff section of $\flow$ admits
 both positive and negative boundaries.
\end{theorem}

Through the section,
 we fix a closed three-dimensional manifold $M$,
 a topologically transitive Anosov flow $\flow$,
 and a smooth and tame Birkhoff section $S$ of $\flow$.
We also fix a coherent pair of orientations
 of $\ofoliation^s$ and $\ofoliation^u$.

The key theme is to relate the curves in the orbit space which avoid the boundary of a Birkhoff section $S$ with curves on the Birkhoff section. When such a curve $c$ starts and ends at the same point, one can continuously lift the curve on the preimage of the Birkhoff section in the universal cover $\wt M$ of $M$. The lifted curve starts and ends on two points on the same orbits, and so one is obtained from the other by pushing along the flow by a time which we denote by $\Delta(c)$ (the drift along $c$). An elementary remark is that $\Delta(c)$ depends only on the points in $\orbitspace$ corresponding to $\partial S$ which lie inside a region bounded by~$c$, and of the local invariant of the surface $S$ around these boundaries component. In particular, when $\Delta(c)$ is negative, there exists a positive boundary component of $\partial S$ associated to a point in $\orbitspace$ lying in a region bounded by $c$.

It remains to find curves with prescribed signs in the orbit space. For that take a lozenge $L$ invariant by the action of an element $g\in\pi_1(M)\setminus\{1\}$. We for $c$ the boundary of a rectangle $R$ which corresponds to a fundamental domain of the action of $g$ on $L$. Then a direct computation yields that the sign of $\Delta(\partial R)$ is opposite to the one of $L$. In particular Proposition \ref{Aprop:lift} stats that when $L$ is of positive type, there exists a positive boundary of $S$ associated to a point inside $L$. The same argument used in the previous section, that when the flow is not $\RR$-covered, there exists both positive and negative lozenges in the orbit space. And so any Birkhoff section has both a positive and a negative boundary component.

\subsection{Drifts along closed curves in the orbit space}
\label{Asec:drift}

Let $\wt S\subset M$ be the lift of $S$ in the universal cover of $M$.
We introduce the drifts
 along a closed curve in $\Int \wt{\Sec}/\wtc{O}$,
 which is a key ingredient of the proof of Theorem~\ref{Athm:A main}.
Let $\gamma\colon [0,1] \ra \Int \wt{\Section}/\wtc{O}$
 be a closed curve.
Take a lift $\wt{\gamma}\colon [0,1] \ra \Int \wt{\Section}$ 
 of $\gamma$ with respect to the covering map
 $Q_\Section\colon \Int \wt{\Section} \ra \Int \wt{\Section}/\wt{O}$ (see Proposition~\ref{Aprop:covering}).
Since the group of deck transformations is generated
 by~$\wt{\returnmap}$,
 there exists a  unique integer $\Delta(\gamma)$ satisfying
\begin{equation*}
 \wt{\gamma}(1)=\wt{\returnmap}^{\Delta(\gamma)}(\wt{\gamma}(0)).
\end{equation*}
For another choice $\wt{\eta}$ of a lift of $\gamma$,
 we have $\wt{\eta}=\wt{\returnmap}^n \circ \wt{\gamma}$
 for some $m \in \ZZ$.
Then,
\begin{align*}
 \wt{\eta}(1) & =\wt{\returnmap}^m \circ \wt{\gamma}(1)
 =\wt{\returnmap}^{m+\Delta(\gamma)} \circ \wt{\gamma}(0)
 =\wt{\returnmap}^{\Delta(\gamma)}\circ \wt{\eta}(0).
\end{align*}
Therefore, $\Delta(\gamma)$ does not depend on the choice of 
 the lift $\wt{\gamma}$.
We call the number~$\Delta(\gamma)$ the \emph{drift} along $\gamma$.

Let $H_1(\Int \wt{\Section}/\wtc{O})$ be
 the first homology group of $\Int\wt{\Section}/\wtc{O}$
 with integer coefficients.
We denote the homology class of a closed curve $\gamma$
 in  $\Int \wt{\Section}/\wtc{O}$
 by~$[\gamma]$.
\begin{lemma}
\label{Alemma:drift H1}
There exists a homomorphism
 $\Delta_*\colon H_1(\Int \wt{\Section}/\wtc{O}) \ra \ZZ$ such that
 $\Delta_*([\gamma])=\Delta(\gamma)$
 for each closed curve $\gamma$ in $\Int \wt{\Section}/\wtc{O}$.
\end{lemma}
\begin{proof}
Fix $\xi_0 \in \Int \wh{\Section}/\wtc{O}$.
Let $\gamma_1,\gamma_2\colon [0,1] \ra \Int \wt{\Section}/\wtc{O}$
 be closed curves starting at $\xi_0$.
By $\gamma_1 * \gamma_2$
 we denote the concatenation of $\gamma_1$ and $\gamma_2$.
Take lifts $\wt{\gamma}_1,\wt{\gamma}_2\colon [0,1] \ra \Int \wt{\Section}$
 of $\gamma_1,\gamma_2$ such that $\wt{\gamma}_1(1)=\wt{\gamma}_2(0)$.
Then, 
\begin{equation*}
 \wt{\gamma}_2(1)=\wt{\returnmap}^{\Delta(\gamma_2)}(\wt{\gamma}_2(0))
 =\wt{\returnmap}^{\Delta(\gamma_2)}(\wt{\gamma}_1(1))
 =\wt{\returnmap}^{\Delta(\gamma_2)+\Delta(\gamma_1)}(\wt{\gamma}_1(0)).
\end{equation*}
This implies that
 $\Delta(\gamma_1 * \gamma_2)=\Delta(\gamma_1)+\Delta(\gamma_2)$.

Let $\eta_1,\eta_2\colon [0,1] \ra \Int \wt{\Section}/\wtc{O}$
 be closed curves starting at $\xi_0$
 which are homotopic as closed curves starting at $\xi_0$.
Take lifts $\wt{\eta}_1,\wt{\eta}_2$ of $\eta_1,\eta_2$
 such that $\wt{\eta}_1(0)=\wt{\eta}_2(0)$.
Since $\eta_1$ and $\eta_2$ are homotopic,
 we have $\wt{\eta}_1(1)=\wt{\eta}_2(1)$.
This implies that $\Delta(\eta_1)=\Delta(\eta_2)$.

We have proved that
 $\Delta$ is additive with respect to concatenation of curves
 and depends only on the homotopy class as a closed curve starting at $\xi_0$.
This implies that the function $\Delta$ induce a homeomorphism
 from the fundamental group of $\Int \wt{\Section}/\wtc{O}$ to $\ZZ$.
Since $\ZZ$ is a commutative group,
 the kernel contains the commutator subgroup of
 the fundamental group.
Therefore, $\Delta$ induces a homeomorphism
 from $H_1(\Int \wt{\Section}/\wtc{O})$ to $\ZZ$.
\end{proof}

Recall that $\orbitspaceS$ is the blow up of $\orbitspace$
 at the discrete subset $\del\wh{\Sec}/\wtc{O}$.
We denote the projection from $\orbitspaceS$ to $\orbitspace$
 by $\Pi_{\del S}$.
For $\xi \in \del\wt{\Sec}/\wtc{O}$,
 the orientation of $\orbitspace$ induces
 an orientation of the boundary curve $\Pi_{\del S}^{-1}(\xi)$.
Let $c_\xi$ be the homology class
 in $\Int{\wt{\Sec}}/\wtc{O}$ such that
 $-c_\xi$ represents $\Pi_{\del \Sec}^{-1}(\xi)$ as an oriented closed curve.
In other words, $c_\xi$ represents a small circle surrounding $\xi$
 with the counter-clockwise orientation.
We call $\Delta_*(c_\xi)$ the \emph{local drift} at $\xi$.
The following formula relates the drift of a simple closed curve
 to the local drifts.

\begin{proposition}
\label{Aprop:local drift}
Let $D$ be an  open simply connected domain
 and $\eta_D$ the boundary curve of $D$.
Suppose that $\eta_D$ does not intersect 
 with $\del \wt{\Sec}/\wtc{O}$. 
Then, $D \cap \del\wt{\Sec}/\wtc{O}$ consists of finitely many points and 
\begin{equation*} 
 \Delta(\eta_D)=\sum_{\xi \in D \cap \del\wt{\Sec}/\wtc{O}} \Delta_*(c_\xi),
\end{equation*}
 where the curve $\eta_D$ is oriented as the boundary curve
 of the domain $D$ in the oriented plane $\orbitspace$.
In particular, if $\Delta(\eta_D)$ is negative then
 $D \cap \del\wt{\Sec}/\wtc{O}$ contains
 a point whose local drift is negative.
\end{proposition}
\begin{proof}
Since $\del \wt{\Sec}/\wtc{O}$ is a discrete subset of $\orbitspace$
 and $D$ is a bounded domain,
 the set $D \cap \del \wt{\Sec}/\wtc{O}$ consists of finitely many points.
Put $D \cap \del\wt{\Sec}/\wtc{O}=\{p_1,\dots,p_k\}$. 
The domain $D$ is bounded by $\eta_D$
 and $\Pi_{\del S}^{-1}(\xi_1), \cdots,\Pi_{\del S}^{-1}(\xi_k)$
 as a domain in the blow up $\orbitspace$.
Notice that the curve $\eta_D$ is oriented as a boundary component of~$D$
 and the orientations
 of $\Pi_{\del S}^{-1}(\xi_1), \cdots,\Pi_{\del S}^{-1}(\xi_k)$
 are opposite to the orientations as boundary components of $D$.
Hence, we have
 $[\eta_D]-(c_{\xi_1}+\dots+c_{\xi_k})=0$ in~$H_1(\wt{\Sec}/\wtc{O})$.
Since $\Delta_*:H_1(\wt{\Sec}/\wtc{O}) \ra \RR$
 is a homomorphism,
\begin{equation*}
 \Delta(\eta_D)-(\Delta_*(c_{\xi_1})+\dots+ \Delta_*(c_{\xi_k}))=0
\qedhere
\end{equation*}
\end{proof}

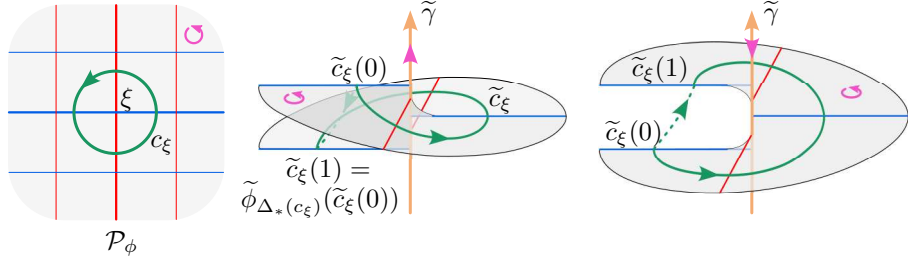
\begin{figure}
    \begin{center}
        \begin{picture}(120,28)(0,0)
            \put(0,0){\includegraphics[width=12cm]{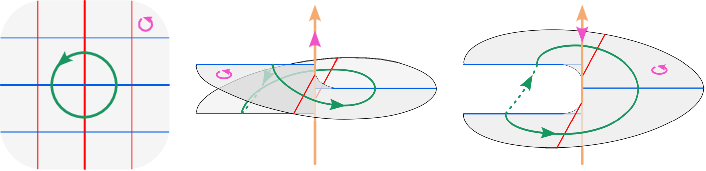}}
            \put(13,-3.5){$\orbitspace$}
            \put(15,15.5){$\xi$}
            \put(19,10){$c_\xi$}
            \put(55,27){$\wt\gamma$}
            \put(100.5,27){$\wt\gamma$}
            \put(64,15.5){$\wt c_\xi$}
            \put(43,19.5){$\wt c_\xi(0)$}
            
            \put(37,6){$\wt c_\xi(1)=$}
            \put(31,2){$\wt\phi_{\Delta_*(c_\xi)}(\wt c_\xi(0))$}
            
            \put(79.5,11){$\wt c_\xi(0)$}
            \put(83,19.5){$\wt c_\xi(1)$}
        \end{picture}
    \end{center}
    \caption{Signs of $\Delta_*(c_\xi)$ for $\xi\in\partial S/\wt{\mathcal{O}}$ corresponding to a positive boundary component in the middle and a negative boundary component on the right. The orientations and signs of boundary are represented in pink.}
    \label{Tfigure:shiftOrbit}
\end{figure}

The signature of the local drift $\Delta_*(c_\xi)$ is determined
 by the multiplicity of the corresponding boundary component
 of $\del\wt{\Sec}$ as follows (see Figure \ref{Tfigure:shiftOrbit} for an illustration).

\begin{proposition}
\label{Aprop:sign}
Let $\xi$ be an element of $\del\wt{\Sec}/\wtc{O}$
 and $\gamma=P(Q^{-1}(\xi))$ the correspondent closed orbit in $\del S$.
Then, the signs of $\Delta_*(c_\xi)$ and $\mult(\gamma)$ are opposite.
\end{proposition}

\begin{proof}

Take a boundary component $\gamma^*$ immersed in $\gamma$.
Recall that the Birkhoff section is smooth and tame.
The intersection of the lift $\pi_{\del \Sec}^{-1}(\foliation^s(\gamma))$
 of the stable manifold of $\gamma$ to the blow up $M_{\del \Sec}$
 and the resolution $S^*$ of $S$ is transverse.
Hence, there exists a curve $C$
 in $\Sec^* \cap \pi_{\del \Sec}^{-1}(\foliation^s(\gamma))$
 and $x_* \in \del C$ such that $C \cap \TT_{\gamma} =\{x_*\}$.
Let $\tau\colon C \ra \RR_+$ and $f\colon C \ra C$
 be the return time function and its
 associated return map of the flow
 $\flow^{\del \Sec}$ on $C$, that is
$$\tau(x)=\inf\{t>0 \mid \flow^{\del \Sec}_t(x) \in C\}$$
 and $f(x)=\flow^{\del \Sec}_{\tau(x)}(x)$. 
Put $T=\tau(x^*)$.
Then, $T$ is equal to
 the period of the closed orbit $\gamma^*$ of $\flow^{\del \Sec}$
 since $\{x_*\}=C \cap \TT_{\gamma}$.

Set $m=|\mult(\gamma^*)|$. We can take a smooth immersion $\psi\colon [0,1] \times [0,mT] \ra S^*$
 such that $\pi_{\Sec \circ }\psi(0,t)=\flow_t \circ \pi_\Sec(x_*)$
 for any $t \in [0,1]$, $\psi(s,t) \in \Int \Sec$ if $s>0$,
 and $\psi(s,mT)=f^m(\psi(s,0))$ for any $s\in [0,1]$.
Let $\wt{\psi}\colon [0,1] \times [0,mT] \ra \wt{S}^*$ be the lift of $\psi$,
 where $\wt{S}^*$ is the lift of $S^*$ to the universal covering of $M$. 
Then, $\wt{\psi}(0,t)=\wt{\flow}_t\vphantom{(}^{\hspace{-1.2mm}\del\Sec}(\wt{\psi}(0,0))$,
 and hence,
 $\wt{\psi}(0,0)$ and $\wt{\psi}(0,mT)$ lie
 in the same $\wt{\flow}\vphantom{(}^{\del\Sec}$-orbit.
Since $\psi(s,0)$ and $\psi(s,mT)$ lie in the same $\flow^{\del \Sec}$-orbit
 for any $s \in [0,1]$, 
 the lifts $\wt{\psi}(s,0)$ and $\wt{\psi}(s,mT)$
 lie in the same $\wt{\flow}\vphantom{(}^{\del\Sec}$-orbit 
 by the continuity of $\wt{\psi}$.
Hence, the projection of the curve $\wt{\psi}(s,\cdot)$ to $\orbitspaceS$ is
 a closed curve for each~$s$.
In particular,
 the projections $Q_\Sec \circ \wt{\psi}(0,\cdot)$
 and $Q_\Sec \circ \wt{\psi}(s,\cdot)$ are homotopic
 as curves in $\orbitspaceS$ for any $s$.

Put $\kappa=1$ if $\mult(\gamma^*)>0$
 and $\kappa=-1$ if $\mult(\gamma^*)<0$.
The closed curve $\psi(0,\cdot)$ is a parametrization
 of the curve $\gamma^*$ whose orientation
 is the same as $\gamma^*$ if $\kappa>0$ and opposite if $\kappa<0$.
Hence, the curve $\psi(0,\cdot)$ represents the homology class
\begin{equation*}
\kappa [\gamma^*] = \kappa (\mult(\gamma^*)[\lambda_{\gamma}]
+ \link(\gamma^*) [\mu_{\gamma}])
\end{equation*}
 in $H_1(\TT_{\gamma_\xi})$.
Let $\wt{\lambda}_\gamma$ and $\wt{\mu}_\gamma$
 be the lifts of $\lambda_\gamma$ and $\mu_\gamma$ to $\wt{M}_{\del \Sec}$,
 respectively.
Since $\lambda_\gamma$ is tangent to the flow $\flow^{\del S}$,
 the projection of the curve $\wt{\lambda}_\gamma$ to $\orbitspaceS$
 is constant.
The lift $\wt{\mu}_\gamma$ is a closed curve
 whose projection to $\orbitspaceS$ represents the homology class $c_\xi$.
This means that the projection of the lift of $[\gamma^*]$
 is a closed curve homotopic to
 $\link(\gamma^*)$-multiple of $c_\xi$.
In particular,
 the projection of the curve $\wt{\psi}(0,\cdot)$ to $\orbitspaceS$
 represents the homology class $\kappa\; \link(\gamma^*) \, c_\xi$.
Since the projections $Q_\Sec \circ \wt{\psi}(0,\cdot)$
 and $Q_\Sec \circ \wt{\psi}(s,\cdot)$ are homotopic
 as curves in $\orbitspaceS$,
 we have $[Q_S \circ \wt{\psi}(s,\cdot)]=\kappa\; \link(\gamma^*)\,c_\xi$.
Hence,
\begin{equation*}
 \kappa\, \Delta_*(\link(\gamma^*)\, c_\xi)
 = \Delta(Q_\Sec \circ \wt{\psi}(s,\cdot))>0.
\end{equation*}
Since the linking number of $\gamma^*$ is negative,
 the local drift $\Delta_*(c_\xi)$ is negative if $\kappa=\mult(\gamma^*)>0$
 and positive if $\kappa=\mult(\gamma^*)<0$.
\end{proof}

\subsection{Lifts of leaves of $\ofoliation^s$}

\begin{lemma}
\label{Alemma:landing} 
Let $\xi$ be an element of $\orbitspace$
 and $\tilde{\alpha}\colon \RR_+ \ra \Int \wt{\Sec}$ a continuous map
 such that $Q \circ \wt{\alpha}(\RR_+) \subset \ofoliation^s(\xi)$
 and $\lim_{u \ra 0}Q \circ \wt{\alpha}(u)=\xi$.
Then, the limit $\lim_{u \ra 0}\wt{\alpha}(u)$ exists
 and is contained in $Q^{-1}(\xi)$.
\end{lemma}
\begin{proof}
Let $\beta\colon \RR_+ \cup \{0\} \ra \orbitspace$ be
 the continuous extension of $Q \circ \wt{\alpha}$ with $\alpha(0)=\xi$.
When $\xi$ is a point of $\Int \wt{\Sec}/\wtc{O}$,
 $\beta$ lifts to a curve
 $\wt{\beta}\colon  \RR_+\cup \{0\} \ra \Int \wt{\Sec}$
 such that $\wt{\beta}(1)=\wt{\alpha}(1)$
Then, $\wt{\beta}(u)=\wt{\alpha}(u)$ for any $u \in \RR_+$,
 and hence,
\begin{equation*}
 \lim_{u \ra 0}\wt{\alpha}(u)=\wt{\beta}(0) \in Q^{-1}(\beta(0))=Q^{-1}(\xi).
\end{equation*}
Suppose that $\xi$ is a point of $\del\wt{\Sec}/\wtc{O}$
 and put $\wt{\gamma}=Q^{-1}(\xi)$.
Since the Birkhoff section $S$ is tame,
 the lift $\lfoliation^s(\wt{\gamma})^*$ of
 the leaf of $\lfoliation^s$ containing $\wt{\gamma}$ to
 the blow up $\wt{M}_{\del S}$ is transverse to the resolution
 $\wt{S}^*$ of $\wt{S}$.
This means that $\lfoliation^s(\wt{\gamma})^* \cap \wt{S}^*$
 is a union of mutually disjoint curves whose
 endpoint is contained in the boundary cylinder $\wt{\TT}_{\wt{\gamma}}$
 corresponding to $\wt{\gamma}$.
In particular, each connected component
 of $\lfoliation^s(\wt{\gamma}) \cap \Int\wt{S}$ lands
 to a point of $\wt{\gamma}$.
Since $\wt{\alpha}(\RR_+)$ is contained in
 $\lfoliation^s(\wt{\gamma}) \cap \Int \wt{\Sec}$,
 this implies that $\lim_{u \ra 0}\wt{\alpha}(u)$ exists
 and it is a point of $\wt{\gamma}$.
\end{proof}

\begin{proposition}
\label{Aprop:lift} 
Let $g$ be a non-trivial element of $\pi_1(M)$,
 $\xi$ a point of $\Fix(g)_+$, and
 $\alpha\colon \RR \ra \ofoliation^s(\xi)$ a homeomorphism
 given by Proposition \ref{Aprop:g action},
 {\it i.e.}, a homeomorphism
 satisfying $\alpha(0)=\xi$ and $g * \alpha(u)=\alpha(2u)$.
Then, $\alpha(\RR_+)$ is contained in $\Int \wt{\Sec}/\wtc{O}$
 and there exists a unique integer $n\geq 1$ such that
 $g \cdot \wt{\alpha}(u)=\wt{\returnmap}^n \circ \wt{\alpha}(2u)$
 for any lift $\wt{\alpha}\colon \RR_+ \ra \Int \wt{\Sec}$ of
 the restriction of $\alpha$ to $\RR_+$ for
 the covering map $Q_\Sec\colon \Int\wt{\Section} \ra \wt{\Sec}/\wtc{O}$.
\end{proposition}
\begin{proof}
Put $\wt{\gamma}=Q^{-1}(\xi)$.
Since $g * \alpha(u)=\alpha(2u)$, $\xi=\alpha(0)$
 is a fixed point of the action of $g$,
 and hence, $P(\wt{\gamma})$ is the unique periodic point of $\flow$
 in $\foliation^s(P(\wt{\gamma}))$.
This implies that
 $Q^{-1}(\alpha(\RR_+))$ does not intersect with $\del\wt{\Section}$,
 and hence, $\alpha(\RR_+)$ is a subset of $\Int \wt{\Section}/\wtc{O}$.
Denote the restriction of $\alpha$ to $\RR_+$ by $\alpha_+$.
Take a lift $\wt{\alpha}_0\colon \RR_+ \ra \Int \wt{\Section}/\wtc{O}$ of $\alpha_+$.
Let $\lambda_2\colon \RR_+ \ra \RR_+$ be the map given by
 $\lambda_2(u)=2u$.
Since $Q \circ (g \cdot \wt{\alpha}_0)=g * \alpha_+ =\alpha_+\circ \lambda_2$,
 the curve $g \cdot \wt{\alpha}_0$ is a lift of $\alpha_+ \circ \lambda_2$.
Since $\wt{\alpha}_0 \circ \lambda_2$ also is a lift of
 $\alpha_+ \circ \lambda_2$,
 there exists an integer $n$ such that
 $g \cdot \wt{\alpha}_0=\wt{\returnmap}^n(\wt{\alpha}_0) \circ \lambda_2$.
By Lemma \ref{Alemma:landing},
 there exists $\wt{p} \in Q^{-1}(\xi)$
 such that $\wt{p}=\lim_{u \ra 0}\wt{\alpha}_0(u)$.
Since $\xi$ is an element of $\Fix(g)_+$,
 there exists $T>0$ such that $g \cdot \wt{p}=\wt{\flow}_T(\wt{p})$.
The equation
 $\wt{\returnmap}^n(\wt{\alpha}_0(u))
 =\wt{\flow}_{\wt{\tau}_n(\wt{\alpha}_0(u))}(\wt{\alpha}_0(u))$
 implies that
 $\wt{\tau}_n(\wt{\alpha}_0(u))$ converges to $T$  as $u$ goes to zero,
 and hence, $n$ is positive.
For another lift $\wt{\alpha}\colon \RR_+ \ra \Int \wt{\Section}/\wtc{O}$
 of~$\alpha_+$,
 there exists an integer $m$ such that
 $\wt{\alpha}=\wt{\returnmap}^m \circ \wt{\alpha}_0$
 since the group of deck transformations of the covering $Q_\Section$
 is generated by $\wt{\returnmap}$.
Recall that the action of $g$ on $\orbitspace$ commutes with
 $\wt{\returnmap}$
 since $\wt{\returnmap}$ is the lift of the return map
 $\returnmap$ to the universal cover $\orbitspace$.
Hence, we have
\begin{equation*}
 g \cdot\wt{\alpha}(u)= g \cdot \wt{\returnmap}^m \circ \wt{\alpha}_0(u)
 =\wt{\returnmap}^m(g \cdot \wt{\alpha}_0(u))
 =\wt{\returnmap}^{m+n}(\wt\alpha_0(2u))
 =\wt{\returnmap}^n \circ \wt{\alpha}(2u).
\qedhere
\end{equation*}
\end{proof}

\subsection{A rectangle associated with a lozenge}
Theorem \ref{Athm:A main} follows from the following proposition,
 which relates lozenges of type $(++)$ and closed curves
 in $\Int \wt{\Sec}/\wtc{O}$ with negative drift.
\begin{proposition}
\label{Aprop:drift lozenge} 
Suppose that $\flow$ admits a lozenge of type $(++)$
 bounded by a closed orbit.
Then, there exists a closed curve $\gamma$ in $\Int\wt{\Sec}/\wtc{O}$
 which represents the boundary of a simply connected domain 
 as an oriented curve and satisfies $\Delta(\gamma)<0$.
\end{proposition}
\begin{corollary}
Suppose that $\flow$ admits a lozenge of type $(++)$ ({\it resp.}, $(+-)$)
 bounded by a closed orbit.
Then, any Birkhoff section $S$ admits a positive ({\it resp.} negative)
 boundary component.
\end{corollary}
\begin{proof}
Suppose that $\flow$ admits a lozenge of type $(++)$.
Remark that the case of type $(+-)$ can be deduced to the case of type $(++)$
 by reversing the orientation of the manifold $M$.

Let $\gamma$ be the closed curve obtained by the above proposition
 and $D$ the simply connected domain bounded by $\gamma$.
Since $\Delta(\gamma)<0$,
 there exists $\xi \in D \cap \del\wt{\Sec}/\wtc{O}$
 such that $\Delta_*(c_\xi)<0$ by Proposition \ref{Aprop:local drift}.
Let $\gamma_\xi$ be the closed curve $P(Q^{-1}(\{\xi\}))$ in $\del \Sec$.
Then, Proposition \ref{Aprop:sign} implies
 that $\mult(\gamma_\xi)>0$, and hence, 
 any boundary component immersed in $\gamma_\xi$ is positive.

\end{proof}
Theorem \ref{Athm:A main} follows from
 this corollary and Corollary \ref{Acor:lozenge R-covered}.

We start the proof of Proposition \ref{Aprop:drift lozenge},
 with some lemmas on the intersection of leaves
 of $\ofoliation^s$ and $\ofoliation^u$.
\begin{lemma}
For any $\xi_1,\xi_2 \in \orbitspace$,
 the intersection $\ofoliation^s(\xi_1) \cap \ofoliation^u(\xi_2)$
 contains at most one point.
\end{lemma}
\begin{proof}
If $\ofoliation^s(\xi_1) \cap \ofoliation^u(\xi_2)$ contains
 two distinct points $\eta_1$ and $\eta_2$,
 the union of segments in $\ofoliation^s(\xi_1)$
 and $\ofoliation^u(\xi_1)$ which connect $\xi_1$ and $\xi_2$
 is a simple closed curve.
By smoothing the curve, we can obtain a smooth simple closed curve
 transverse to $\ofoliation^s$.
It contradicts Poincar\'e-Bendixon's theorem
 since $\orbitspace$ is homeomorphism to $\RR^2$.
\end{proof}

\begin{lemma}
\label{Alemma:intersection} 
Suppose that $\xi_1, \xi_2 \in \orbitspace$ satisfy
 that $\ofoliation^s_+(\xi_1) \cap \ofoliation^u(\xi_2) \neq \emptyset$
 and $\ofoliation^s(\xi_2) \cap \ofoliation^u(\xi_1) \neq \emptyset$.
Then,  $\ofoliation^s_-(\xi_2) \cap \ofoliation^u(\xi_1) \neq \emptyset$.
Similarly,
 if $\ofoliation^u_+(\xi_1) \cap \ofoliation^s(\xi_2) \neq \emptyset$
 and $\ofoliation^u(\xi_2) \cap \ofoliation^s(\xi_1) \neq \emptyset$,
 then $\ofoliation^u_-(\xi_2) \cap \ofoliation^s(\xi_1) \neq \emptyset$.
\end{lemma}
\begin{proof}
We only show the former since the proof of the latter is same as the former.
If $\ofoliation^u(\xi_1)$ contains $\xi_2$,
 then, $\ofoliation^u(\xi_1)=\ofoliation^u(\xi_2)$
 intersects with $\ofoliation^s(\xi_1)$ at $\xi_1$
 and another point in $\ofoliation^s_+(\xi_1)$.
It contradicts the previous lemma.
Hence, $\ofoliation^u(\xi_1)$ does not contain $\xi_2$.
If $\ofoliation^s_+(\xi_2) \cap \ofoliation^u(\xi_1) \neq \emptyset$,
 then there exists a simple closed curve $\gamma\colon [0,4] \ra \orbitspace$
 such that $\gamma((0,1)) \subset \ofoliation^s_+(\xi_1)$,
 $\gamma([1,2]) \subset \ofoliation^u(\xi_2)$,
 $\gamma((2,3)) \subset \ofoliation^s_+(\xi_2)$,
 $\gamma([3,4]) \subset \ofoliation^u(\xi_1)$.
By smoothing this curve, we can obtain a smooth simple closed curve
 which is transverse to $\ofoliation^u$.
It contradicts Poincar\'e-Bendixon's theorem,
 and hence, $\ofoliation^s_+(\xi_2) \cap \ofoliation^u(\xi_1) = \emptyset$.
Since $\ofoliation^s(\xi_2) \cap \ofoliation^u(\xi_1) \neq \emptyset$,
 we have $\ofoliation^s_-(\xi_2) \cap \ofoliation^u(\xi_1) \neq \emptyset$.
\end{proof}

Next, we introduce a rectangle in $\orbitspace$
 bounded by segments in leaves of $\ofoliation^s$ and $\ofoliation^u$.
A continuous map $\surectangle\colon [0,1]^2 \ra \orbitspace$
 is called an \emph{$su$-rectangle}
 if it is a continuous embedding such that
\begin{alignat*}{5}
\surectangle((0,1]\times\{0\})
 & \subset \ofoliation^s_+(\surectangle(0,0)), & \hsp
\surectangle(\{0\} \times (0,1])
 & \subset \ofoliation^u_+(\surectangle(0,0)), \\
\surectangle([0,1)\times\{1\})
 & \subset \ofoliation^s_-(\surectangle(1,1)), & \hsp
\surectangle(\{1\} \times [0,1))
 & \subset \ofoliation^u_-(\surectangle(1,1)).
\end{alignat*}
We define the \emph{boundary curve}
 $\gamma_{\surectangle}\colon  [0,4] \ra \surectangle(\del [0,1]^2)$
 by
\begin{alignat*}{5}
 \gamma_{\surectangle}(t)=
\begin{cases}
\surectangle(t,0) & (0 \leq t<1),\\
\surectangle(1,t-1) & (1 \leq t<2),\\
\surectangle(3-t,1) & (2 \leq t<3),\\
\surectangle(0,4-t) & (3 \leq t \leq 4).
\end{cases}
\end{alignat*}

\begin{lemma}
For an $su$-rectangle $\surectangle$,
 the boundary curve $\gamma_{\surectangle}$ is a positively oriented
 parametrization of the boundary of the domain $\surectangle((0,1)^2)$.
\end{lemma}
\begin{proof}
Take $0<t<1$ and put $\eta=\surectangle(t,0)$.
First, we show that
 neither $\ofoliation^u_+(\eta)$ nor $\ofoliation^u_-(\eta)$ intersect with
 $\gamma_{\surectangle}([0,2] \cup [3,4])$.
Since $\ofoliation^s(\surectangle(0,0))$ intersects
 with $\ofoliation^u(\eta)$ at $\eta$
 and $\gamma_{\surectangle}([0,1])$
 is a subset of $\ofoliation^s(\surectangle(0,0))$,
 the uniqueness of intersection of leaves of $\ofoliation^s$
 and $\ofoliation^s$ implies that
 neither $\ofoliation^u_+(\eta)$ nor $\ofoliation^u_-(\eta)$
 intersect with $\gamma_{\surectangle}([0,1])$.
The leaves $\ofoliation^u(\surectangle(0,0))$
 and $\ofoliation^u(\surectangle(1,1))$ intersect with
 $\ofoliation^s(\eta)$ at $\surectangle(0,0)$
 and $\surectangle(1,0)$ respectively.
Since the leaf $\ofoliation^s(\eta)$ intersects with
 $\ofoliation^u(\eta)$ at $\eta \neq \surectangle(0,0),\surectangle(1,0)$,
 the uniqueness of intersection of leaves
 of $\ofoliation^s$ and $\foliation^u$
 implies that $\ofoliation^u(\eta)$ coincides with
 neither $\ofoliation^u(\surectangle(0,0))$ nor
 $\ofoliation^u(\surectangle(1,1))$.
Hence, $\ofoliation^u(\eta)$ does not
 intersect with $\gamma_{\surectangle}([1,2] \cup [3,4])$.
Therefore,
 neither $\ofoliation^u_+(\eta)$ nor $\ofoliation^u_-(\eta)$
 intersect with $\gamma_{\surectangle}([0,2] \cup [3,4])$.

Since $\ofoliation^s(\eta)$ intersects with
 $\ofoliation^u_-(\surectangle(0,1))$ at $\surectangle(0,0)$, 
Lemma \ref{Alemma:intersection} implies that
 $\ofoliation^u_-(\eta)$ does not intersect
 with $\ofoliation^s(\surectangle(0,1))$.
In particular, $\ofoliation^u_-(\eta)$ does not intersect with
 $\gamma_{\surectangle}([2,3])$, which is a subset
 of $\ofoliation^s(\surectangle(1,1))=\ofoliation^s(\surectangle(0,1))$,
Therefore, $\ofoliation^u_-(\eta) \cap \surectangle(\del [0,1]^2)=\emptyset$.
This means that $\ofoliation^u_-(\eta)$ does not intersect with
 $\surectangle([0,1]^2)$
 since $\surectangle([0,1]^2)$ is bounded simply connected domain.
Therefore, `the right-hand side' of the curve $\gamma_{\surectangle}$
 is outside $\surectangle([0,1]^2)$,
 and hence,
 the orientation of the curve $\gamma_\surectangle$
 coincides with the orientation as a boundary curve of
 the simply connected domain $\surectangle([0,1]^2)$.
\end{proof}

Let $L$ be a lozenge of type $(++)$ with corners $(\xi_1,\xi_2)$.
Suppose that the action of
 a non-trivial element $g$ of $\pi_1(M)$
 fixes both $\xi_1$ and $\xi_2$ and $\xi_1 \in \Fix(g)_+$.
\begin{lemma}
\label{Alemma:rectangle}
The corner $\xi_2$ is contained in $\Fix(g^{-1})_+$
 and there exist $\eta_1 \in \ofoliation^s_+(\xi_1)$,
 $\eta_2 \in \ofoliation^s_-(\xi_2)$, and
 an $su$-rectangle $\surectangle\colon [0,1]^2 \ra \orbitspace$
 such that $\surectangle(0,0)=\eta_1$,
 $\surectangle(1,0)=g*\eta_1$, $\surectangle(1,1)=g * \eta_2$, 
 and $\surectangle(0,1)=\eta_2$.
\end{lemma}

\begin{proof}
Take $\eta \in L$.
Then, $\ofoliation^s(\eta)$ intersects with
 both $\ofoliation^u_+(\xi_1)$ and $\ofoliation^u_-(\xi_2),$
 and $\ofoliation^u(\eta)$ intersects with
 both $\ofoliation^s_+(\xi_1)$ and $\ofoliation^s_-(\xi_2)$.
By Lemma \ref{Alemma:intersection},
 $\ofoliation^u_-(\eta)$ intersects with $\ofoliation^s(\xi_1)$.
The uniqueness of intersection point implies
 that $\ofoliation^u_-(\eta)$ intersects with $\ofoliation^s_+(\xi_1)$.
We denote this point by $\eta_1$. 
Since $\xi_1 \in \Fix(g)_+$ and
 $\eta_1 \in \ofoliation_+^s(\xi_1)$,
 Proposition \ref{Aprop:g action} implies that
 $g* \eta_1 \in \ofoliation^s_+(\eta_1)$.
In the same way, we can obtain the unique intersection point
 $\eta_2$ of  $\ofoliation^s_-(\xi_2)$ and $\ofoliation^u_+(\eta)$.
This means that $\ofoliation^s_-(\xi_2)$ and $\ofoliation^u_+(\eta_1)$
 intersects at $\eta_2$.
Since the action of $g$ preserves the $\ofoliation^s$
 and its orientation and $\xi_2$ is a fixed point of the action of $g$,
 $g*\eta_2$ is contained in both $\ofoliation^u_+(g*\eta_1)$
 and $\ofoliation^s(\eta_2)=\ofoliation^s(\xi_2)$.
The half leaf $\ofoliation^s_+(\eta_1)$
 intersects with $\ofoliation^u(g*\eta_2)$ at $g*\eta_1$
 and the leaf $\ofoliation^s(g*\eta_2)$ intersects
 with $\ofoliation^u_+(\eta_1)$ at $\eta_2$.
By Lemma \ref{Alemma:intersection},
 $\ofoliation^s_-(g*\eta_2)$ intersects with
 $\ofoliation^u_+(\eta_1)$ at $\eta_2$.
Therefore, there exists an $su$-rectangle $\surectangle$
 such that $\surectangle(0,0)=\eta_1$, $\surectangle(0,1)=g * \eta_1$,
 $\surectangle(1,1)=g*\eta_2$, $\surectangle(1,0)=\eta_2$,
Since $\eta_2$ is contained in both $\ofoliation^s_-(g*\eta_2)$
 and $\ofoliation^s_-(\xi_2)$, Proposition \ref{Aprop:g action}
 implies that $\xi_2$ is a point of $\Fix(g^{-1})_+$.
\end{proof}

We show that the boundary curve of the $su$-rectangle $\surectangle$
 in the above lemma satisfies the required property
 in Proposition \ref{Aprop:drift lozenge}.

\begin{proof}[Proof of Proposition \ref{Aprop:drift lozenge}]
Let $L$ be a lozenge of type $(++)$ with corners $(\xi_1,\xi_2)$
 such that $\xi_1$ and $\xi_2$ are fixed by the action
 of a common non-trivial element $g$ of $\pi_1(M)$.
We may assume that $\xi_1$ is contained in $\Fix(g)_+$.
Let $\surectangle$ be the $su$-rectangle
 obtained by Lemma \ref{Alemma:rectangle}
 and $\gamma$ its boundary curve.
Put $\eta_1=\surectangle(0,0)$ and $\eta_2=\surectangle(0,1)$.
Then, $\surectangle(1,0)=g* \eta_1$
 and $\surectangle(0,1)=g* \eta_2$.

We construct a lift $\wt{\gamma}\colon [0,4] \ra \Int\wt{\Sec}/\wtc{O}$
 of the boundary curve $\gamma$
 and show that $\wt{\gamma}(4)=\wt{\Phi}^{-N}(\wt{\gamma}(0))$
 for some $N \geq 1$.
By Proposition \ref{Aprop:g action},
 we can take homomorphisms $\alpha_1\colon \RR \ra \ofoliation^s(\xi_1)$
 and $\alpha_2\colon \RR \ra \ofoliation^s(\xi_2)$ such that
 $\alpha_i(0)=\xi_i$, 
 $\alpha_1(\RR_+)=\ofoliation^s_+(\xi_1)$,
 $\alpha_1(\RR_+)=\ofoliation^s_-(\xi_2)$,
 and $g * \alpha_1(u)=\alpha_1(2u)$,
 $g^{-1}* \alpha_2(u)=\alpha_2(2u)$ for any $u \in\RR$.
Put $u_i=\alpha_i^{-1}(\eta_i)$ for $i=1,2$.
Then, $u_i>0$,
 $g * \eta_1=g *\alpha_1(u_1)=\alpha_1(2u_1)$,
 and $g * \eta_2=g * \alpha_2(u_2)=\alpha_2(u_2/2)$.
Let $\wt{\beta}\colon [0,1] \ra \Int\wt{\Sec}$ be the lift
 of the curve $t \mapsto \surectangle(0,t)$.
Since $\surectangle(0,0)=\eta_1=\alpha_1(u_1)$,
 there exist a lift $\wt{\alpha}_1\colon \RR_+ \ra \Int \wt{\Sec}$ of 
 the restriction of $\alpha_1$ to $\RR_+$
 such that $\wt{\alpha}_1(u_1)=\wt{\beta}(0)$.
Similarly, we can take a lift $\wt{\alpha}_2\colon \RR_+ \ra \Int \wt{\Sec}$ of
 the restriction of $\alpha_2$ to $\RR_+$
 such that $\wt{\alpha}_2(u_2)=\wt{\beta}(1)$.
Notice that the restriction of $\wt{\alpha}_1$ to $[u_1,2u_1]$
 is a reparametrization of a lift of the curve $t \mapsto \surectangle(t,0)$
 and the restriction of $\wt{\alpha}_2$ to $[u_2/2,u_2]$
 is a reparametrization of a lift of the curve $t \mapsto \surectangle(1-t,1)$.
By Proposition \ref{Aprop:lift},
 there exist $N_1,N_2 \geq 1$ such that
 $g \cdot \wt{\alpha}_1(u)
 =\wt{\returnmap}^{N_1} \circ \wt{\alpha}_1(2u)$
 and 
 $g^{-1} \cdot \wt{\alpha}_2(u)
 =\wt{\returnmap}^{N_2} \circ \wt{\alpha}_2(2u)$
 for any $u \in \RR_+$.
Therefore,
\begin{align*}
 g \cdot \wt{\beta}(0) & = g \cdot \wt{\alpha}_1(u_1)=
 \wt{\returnmap}^{N_1}(\wt{\alpha}_1(2u_1)),\\
 g \cdot \wt{\beta}(1) & = g \cdot \wt{\alpha}_2(u_2)=
 \wt{\returnmap}^{-N_2}(\wt{\alpha}_2(u_2/2)).
\end{align*}
This implies that
 $g \cdot \wt{\beta}$ is the lift of 
 the curve $t \mapsto \surectangle(1,t)$ which connects
 $\wt{\returnmap}^{N_1}(\wt{\alpha}_1(2u_1))$
 and $\wt{\returnmap}^{-N_2}(\wt{\alpha}_2(u_2/2))$.


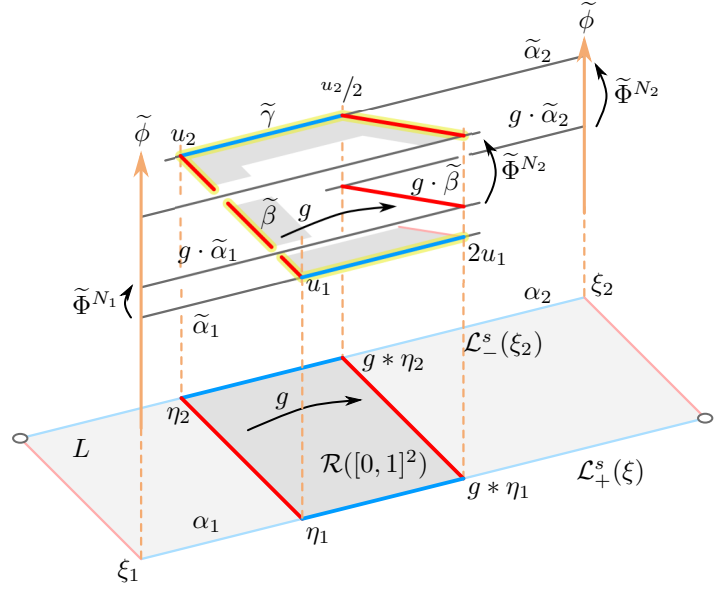
\begin{figure}
    \begin{center}
        \begin{picture}(93,60)(0,0)
            \put(0,0){\includegraphics[height=7cm]{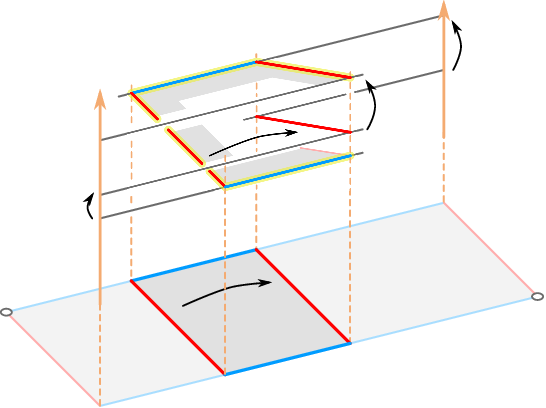}}
            \put(33, 58.5){$\wt\gamma$}
            \put(33,45){$\wt\beta$}
            \put(53,49.5){$g\cdot\wt\beta$}

            \put(14,-2){$\xi_1$}
            \put(77,36){$\xi_2$}

            \put(24, 4.5){$\alpha_1$}
            \put(75,11){$\ofoliation^s_+(\xi)$}
            \put(68,35){$\alpha_2$}
            \put(60,28){$\ofoliation^s_-(\xi_2)$}

            \put(20.5,19){$\eta_2$}
            \put(39, 3){$\eta_1$}
            \put(46.5,26){$g*\eta_2$}
            \put(60.5, 9){$g*\eta_1$}

            \put(39,36){$u_1$}
            \put(60.5,40){$2u_1$}
            \put(21,56){$u_2$}
            \put(41,62){$\sfrac{u_2}{2}$}

            \put(24, 30.5){$\wt\alpha_1$}
            \put(22, 40.5){$g\cdot\wt\alpha_1$}
            \put(68, 67){$\wt\alpha_2$}
            \put(66, 58.5){$g\cdot\wt\alpha_2$}

            \put(16.2,56){$\wt\phi$}
            \put(75.3,71){$\wt\phi$}
            \put(8,33){$\wt\Phi^{N_1}$}
            \put(80,61){$\wt\Phi^{N_2}$}
            \put(65,51){$\wt\Phi^{N_2}$}

            \put(35,21){$g$}
            \put(38,46){$g$}
            \put(41,12){$\surectangle([0,1]^2)$}
            \put(8,14){$L$}
        \end{picture}
    \end{center}
    \caption{Proof of Proposition \ref{Aprop:drift lozenge}. The curve $\wt\gamma$ is highlighted in yellow.}
    \label{Tfigure:lozenge}
\end{figure}

Define a curve $\wt{\gamma}\colon [0,4] \ra \orbitspace$ by
\begin{equation*}
 \wt{\gamma}(t)=
\begin{cases}
\wt{\returnmap}^{N_2}(g \cdot \wt{\beta}(t)) & (0 \leq t \leq 1),\\
\wt{\alpha}_2(tu_2/2) & (1 <t \leq 2),\\
\wt{\beta}(3-t) & (2<t \leq 3),\\
\wt{\alpha}_1((t-2)u_1) & (3 < t \leq 4).
\end{cases}
\end{equation*}
Then, $\wt{\gamma}$ is a lift of
 the boundary curve $\gamma$ of the domain $\surectangle((0,1)^2)$.
Since
\begin{equation*}
 \wt{\gamma}(4)=\wt{\alpha}_1(2u_1)
 =\wt{\returnmap}^{-N_1}(g \cdot \wt{\beta}(0))
 =\wt{\returnmap}^{-(N_1+N_2)}\wt{\gamma}(0)
\end{equation*}
 and $N_1,N_2$ are positive, the drift $\Delta({\gamma})$ is negative.
\end{proof}

\begin{remark}
    The proof of Proposition \ref{Aprop:drift lozenge} can be used to reprove one of Barbot's theorem: an Anosov flow which preserves a transverse contact structure is $\RR$-covered and twisted \cite{Ba01}. To see that, use the contact distribution (a plane distribution transverse to the flow) for the drift operation instead of the Birkhoff section. That is, given a curve $c$ in the orbit space, lift it to $\wt M$ keeping the lift tangent to the contact structure. It yields a drifting quantity $\Delta(c)$ for the curve~$c$. Similarly to the proof of Proposition \ref{Aprop:drift lozenge}, the sign of the drift of the boundary of a $su$-rectangle associated to a lozenge $L$ depends only on the sign of the lozenge. The natural sign that comes with the contact structure ensure that the drifts of all these $su$-rectangles have the same sign, and so all lozenges have the same sign. Then using arguments similar to ours, the flow is $\RR$-covered and skewed (it cannot be the reparametrization of a suspension). 

    The idea of drifting a curve was actually inspired by Barbot's article. 
\end{remark}

%% file: SectionBonatti/SectionBonatti.tex
In this section, we give a third proof of Theorem \ref{theorem:rcoveredcondition} the approach is that the foliations on the Birkhoff section $\Section$ allow us to describe the holonomies of the foliations on the bi-foliated plane associated to the flow (see Theorem~\ref{t.holonomies}). 

Then for proving Theorem~\ref{theorem:rcoveredcondition}, we will use the following characterization of positively twisted $\RR$-covered Anosov flows: 

Let $\flow$  be an Anosov flow on an oriented $3$-manifold  and $(\orbitspace,\ofoliation^s,\ofoliation^u)$ its bi-foliated plane, endowed with two coherent orientations of $\ofoliation^s$ and $\ofoliation^u$. 
We say that the quadrant $C_{+,+}(x)$ is \emph{complete} (or satisfy the \emph{complete holonomy property}) if $\ofoliation^u(y) \cap \ofoliation^s(z) \neq \emptyset$ for any $y\in \ofoliation^s_+(x)$ and $z\in \ofoliation^u_+(x)$.
By Corollary 2.1 in~\cite{BoIa}, if $\flow$ is either non-$\RR$-covered or negatively twisted, then there exists $x \in \orbitspace$ such that the quadrant $C_{+,+}(x)$ is not complete.

Proposition~\ref{p.complete} asserts that the positivity of a Birkhoff section implies that all $C_{+,+}$-quadrants are complete. To see that, we take a point $y\in\ofoliation^s_+(x)$, consider the stable segment $[x,y]^s$ and push it along the unstable foliation in the positive direction. The goal is to prove that it may be pushed to a stable segment on all stable leaves $\ofoliation^s(z)$ for $z\in\ofoliation^u_+(x)$. Morally speaking, when pushing the segment along the unstable foliation, it will cross boundary points of the Birkhoff section. When it crosses a point, the length (as explained below) of the segment may increase. If it increases too much, it may cross more and more boundary points, and its length may diverges in finite time.

The flow induces a pseudo-Anosov map on the Birkhoff section (after contracting every boundary component to a point), which naturally comes with a measure on stable leaves. The length of a stable segment mentioned above is given by integrating this measure. 

The proof strategy goes as follows. Given a stable segment $[z,w]^s$ obtained as above by pushing $[x,y]^s$ along the unstable foliation, there is a natural lift of the segment to the Birkhoff section (with eventually one singularity). We measure the length of the lift for the measure given by the pseudo-Anosov map. The key point is that when the segment $[z,w]^s$ goes through a positive boundary of the Birkhoff section, the length of the lifted segment decreases. Thus, when the flow admits a positive Birkhoff section, the length of $[z,w]^s$ is decreasing in $z$ and so bounded above. Therefore we may push the segment $[x,y]^s$ along the unstable foliation indefinitely, so that the quadrant $C^{+,+}(x)$ is complete. 

\subsection{Comparison of the holonomies of $\ofoliation^s,\ofoliation^u$ and $\ofoliation^s_\returnmap,\ofoliation^u_\returnmap$}

\paragraph*{Foliations on $\wh\Section_\partial$ and holonomies.}

Let $\flow$ be a transitive Anosov flows with oriented stable and unstable bundles. Fix a Birkhoff section $\Section$. Recall that we denote by $\wh\Section_\partial$ the surface obtained by blowing down each boundary component of $\Section$, by $\Si\subset \wh\Section_\partial$ the image of $\partial \Section$, and by $\returnmap:\wh\Section_\partial\to\wh\Section_\partial$ the first return map on $\wh\Section_\partial$. The map $\returnmap$ is a pseudo-Anosov homeomorphism, so it leaves invariant two transverse transversely measured foliations (singular on $\Sigma$) $\ofoliation^s_\returnmap,\ofoliation^u_\returnmap$ on $\wh\Section_\partial$. Note that a point in $\Sigma$ may have a neighborhood bi-foliated as a 2-prong, so that it is in fact a regular point of the (topological) bi-foliation on $\wh\Section_\partial$. By abuse of language, we still call these points singular in the rest of the section. Let $\mu^u,\mu^s$ be the transverse measures to $\ofoliation^s_\returnmap$ and $\ofoliation^u_\returnmap$, respectively. There is a dilation rate $\lambda>1$ so that 
$$\mu^u(\returnmap(I^u))=\lambda \mu^u(I^u)$$ 
for any segment $I^u$ of a leaf of $\ofoliation^u_\phi$ and in the same way, 
$$\mu^s(\returnmap(I^s))=\lambda^{-1} \mu^s(I^s).$$
Note that $\ofoliation^s_\returnmap,\ofoliation^u_\returnmap$ are induced by the intersection of the weak foliations $\foliation^{s},\foliation^{u}$ of $\flow$ with $\Section\setminus\partial \Section$. We fix the orientations on the foliations $\ofoliation^s_\returnmap$ and $\ofoliation^u_\returnmap$ induced by the orientations of $\foliation^{ss},\foliation^{uu}$. We suppose that the orientation on $\ofoliation^s_\returnmap$ plus the one on $\ofoliation^u_\returnmap$ induces the same orientation on $\wh\Section_\partial$ that the one induced by the flow.

Recall that there is a flat metric on $\wh\Section_\partial$, singular only at $\Si$, for which for any stable and unstable segments $I^u$ and $I^s$, $I^s$ and $I^u$ are orthogonal geodesics of lengths respectively $\mu^s(I^s)$ and $\mu^u(I^u)$. Up to now the use of this metric will be implicit. 

\begin{figure}
    \begin{center}
        \begin{picture}(80,42)(0,0)
            \definecolor{myblue}{rgb}{0, 0.6, 1}
            \put(0,0){\includegraphics[width=8cm]{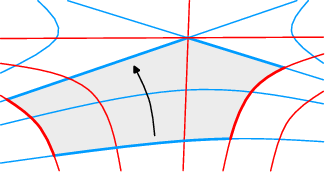}}
            \put(11.5,1.5){$p$}
            \put(0.5,15){$q$}
            \put(35,3.4){$I^s_+(p)$}
            \put(21,28.4){$I^s_+(q)$}
        \end{picture}
    \end{center}
    \caption{Holonomy between two stable curves (in blue) along the unstable foliation $\ofoliation^u_\returnmap$ (in dark red).}
    \label{Cfigure:Holonomy}
\end{figure}

Recalling that, by definition of transversely measured foliations, a transverse measure is holonomy invariant. Meaning if $I_1^s,I_2^s$ are two stable segments isotopic along the foliation $\ofoliation^u_\returnmap$, then $\mu^s(I_1^s)=\mu^s(I_2^s)$.
Next  observation is tautological but crucial for our argument (illustrated in Figure \ref{Cfigure:Holonomy}): 

\begin{proposition}\label{p.tauto}
Let $I^u= [p,q]^u$ be an unstable segment and $\ell>0$ so that every stable oriented segment of the form
$$I^s_+(x)= [x, y]^s, \; x\in[p,q]^s,\mbox{ and }\mu^s(I^s_+(x))= \ell$$ is disjoint from $\Si$.  Then the holonomy $H^u_{p,q}\colon \ofoliation^s_\returnmap(p)\to \ofoliation^s_\returnmap(q)$ of $\ofoliation^u_\returnmap$ is well-defined on $I^s_+(p)$, its image is $I^s_+(q)$ and it is the unique isometry $I^s_+(p)\to I^s_+(q)$ mapping $p$ on $q$. 

The same statement holds for the family $I^s_-(x)= [z,x]^s$, $x\in [p,q]^u$ and $\mu^s(I^s_+(x))= \ell$. 
\end{proposition}

One may extend by continuity this statement to the case where the segments $I^s_\pm(p)$ and/or  $I^s_{\pm}(q)$ contain a point of $\Si$, as follows
\begin{corollary}\label{c.tauto}
Assume that every stable segment of the form
$$I^s_+(x)= [x, y]^s, x\in(p,q)^s,\mbox{ and }\mu^s(I^s_+(x))= \ell$$ is disjoint from $\Si$. Let $I^s_+(p)$ and $I^s_+(q)$ denote the limits of $I^s_+(x)$ when $x$ tend to $p$ and $q$ respectively. 
 
Then
\begin{itemize} \item $I^s_+(p)$ and $I^s_+(q)$ are stable segment of $\mu^s$-length $\ell$,
\item the holonomy $H^u_{p,q}\colon I^s_+(p)\to I^s_+(q)$ of $\ofoliation^u_\returnmap$ is well-defined  and is the unique isometry $I^s_+(p)\to I^s_+(q)$ mapping $p$ on $q$. 
\end{itemize}

 The same statement holds for the family $I^s_-(x)= [z,x]^s$, $x\in (p,q)^u$.  
\end{corollary}


Corollary~\ref{c.tauto} looks straightforward but indeed the isometry would be wrong if we admitted $1$-prong singular points. Our hypothesis on the orientability of the stable/unstable foliations of $\flow$ is mostly motivated by this fact. We suppose the foliations orientable to make the proof more conformable, but the ideas hold for the non-orientable case.

\begin{remark}\label{r.bounded} A stable or unstable segment meets $\Si$ on at most $1$ point. 
\end{remark}

\begin{proof} If $I^s_+(p)$ and $I^s_+(q)$ are disjoint from $\Si$, the item $1$ just follows from the continuous dependence of the leaves. Assume that $I^s_+(p)\cap \Si=\sigma$ (the intersection is at most $1$ point). Let $J^u$ be the union of the local unstable prongs of $\sigma$ in a small enough neighborhood: thus the intersection of the segment $I^s_+(x)$  with $J^u$ is at most $1$ point $t(x)$. We split $I^s_+(x)= [x,t(x)]^s\cup[t(x),y]^s$.  Each of this segment has a bonded length, the sum of both lengths is $\ell$ and each tend to a compact segment in a prong of $\sigma$ (with the limit lengths, whose sum is $\ell$): it remains to see that the two limits are on different prongs, and therefore their intersection is just $\sigma$. 

Now the second item is straightforward. 
\end{proof}

In other words, the holonomies of $\ofoliation^u_\returnmap$ on the stable leaves, following a given unstable leaf $\ofoliation^u_0$ consist in flowing isometrically the stable segment starting (or ending) at $x\in \ofoliation^u_0$ and this procedure holds as long as the stable segment does not cross a singular point. 

One deduces easily, arguing by contradiction, the following property: 
\begin{corollary}\label{c.bounded} Given $\ell>0$ there is $\delta>0$ so  that, for any $[p,q]^u$ if 
there are $x,y\in [p,q]^u$  and $\sigma_x,\sigma_y\in \Si$ with $\sigma_x\in \ofoliation^s_\returnmap(x)$, $\sigma_y\in \ofoliation^s_\returnmap(y)$ and 
$$\mu^s([x,\sigma_x]^s)\leq \ell \mbox{ and } \mu^s([y,\sigma_y]^s)\leq \ell$$ then 
$$\mu^u([p,q]^u)>\delta.$$ 
\end{corollary}

\paragraph*{Parametrizations of the leaves of $\ofoliation^s$, $\ofoliation^u$ by the ones of $\ofoliation^s_\returnmap$, $\ofoliation^u_\returnmap$: the lift/projection procedure.}

Notice that $(M\setminus \partial \Section, \flow)$ is orbitally equivalent to the suspension flow of the restriction of  $\returnmap$ to $Int(\Section)$ or equivalently of 
$\returnmap$ to $\wh\Section_\partial\setminus \Si$, by an equivalence inducing the identity map on $\wh\Section_\partial\setminus \Si$. 

Let $P\colon \wt M\to M$ be the universal cover. 
Consider the lift $\wt \Section$ of $\Section$ on $(\wt M,\wt \flow)\simeq (\RR^3, \partial/\partial x)$. Note that the first return $\returnmap$ admits a canonical lift $\wt \returnmap$ on $\wt \Section\setminus \partial \wt \Section$ which is the first return map of $\wt \flow$. Denote by  $\wt \Ga$ the points of $\orbitspace$ corresponding to the periodic orbits of $\flow$ in $\partial \Section$, that is $\wt\Ga = \partial\wt{S}/\wt{\mcal{O}}$.

Recall Lemma \ref{Aprop:covering}: the projection $\Pi\colon\wt \Section\setminus \partial \wt \Section\to \orbitspace\setminus \wt \Ga$ along the $\wt \flow$-orbits is a infinite cyclic cover whose deck transformation group is generated by $\wt \returnmap$. This allows us to compares the holonomies of $\ofoliation^s$ and $\ofoliation^u$ with the holonomies of $\ofoliation^s_\returnmap$ and $\ofoliation^u_\returnmap$.

More precisely, let $l^u$ be an oriented leaf of $\ofoliation^u$ disjoint from $\wt \Ga$. Consider
a lift $\wt l^u_\returnmap$  on $\wt \Section$. Choose a point $O\in\wt l^u_\returnmap$ that we will call \emph{the origin}, 
and its projection $l^u_\returnmap=P(\wt l^u_\returnmap)$ on $\Section$: it is a leaf of $\ofoliation^u_\returnmap$ disjoint from $\Si$. 
One gets a parametrization $$\theta\colon l^u_\returnmap\to l^u $$
of $l^u$ by $l^u_P$ by considering the composition of  $P^{-1}\colon l^u_\returnmap\to \wt l^u_\returnmap$ composed with the projection $\Pi\colon \wt l^u_\returnmap\to l^u$. In other words, 

$$\theta=\Pi\circ (P|_{\wt l^u_\returnmap})^{-1}.$$ 

Note that $\theta$ only depends on the origin $O$: the origin $O$ determine $\wt l^u_\returnmap$, $ l^u_\returnmap=P(\wt l^u_\returnmap)$ , $l^u=\Pi(\wt l^u_\returnmap)$, $(P_{\wt l^u_\returnmap})^{-1})$ \dots.

Now, at any point $x\in l^u_\returnmap$ for which $\ofoliation^s_\returnmap(x)$ is disjoint from $\Si$,  the lift of $\ofoliation^s_\returnmap(x)$ on $\wt \Section$ through  $(P|_{\wt \ofoliation^u_\returnmap})^{-1}(x)\in \wt l^u_\returnmap$ is well-defined and provides by projection by $\Pi$ a parametrization 

$$\theta^s_x\colon \ofoliation^s_\returnmap(x)\to  \ofoliation^s(\theta(x)).$$ 
Notice that $\theta^s_x$ only depends on the origin $O$. It is called the map obtained by \emph{the lift/projection procedure with origin at $O$}. 

The lift/projection procedure with origin at $O$ can be applied to any path $\gamma\colon [0,1]\to \Section\setminus\partial \Section$ with $\gamma(0)=P(O)$, leading to a map $\theta_\gamma\colon [0,1]\to \orbitspace$ with $\theta_\gamma(0)=\Pi(O)$, and $\theta_\gamma(1)$ only depends on $O$ and of the homotopy class, with fixed ends, of $\gamma$. 
More generally the lift/projection procedure with origin at $O$ can be applied
to any continuous map $h\colon \De\to \Section\setminus \partial \Section$, $h(o)=P(O)$ where  $\De$ is a path connected simply connected set $\De$, with an origin $o$. 

\begin{figure}
    \begin{center}
        \begin{picture}(120,33)(0,0)
        \put(-10,0){\includegraphics[width=14cm]{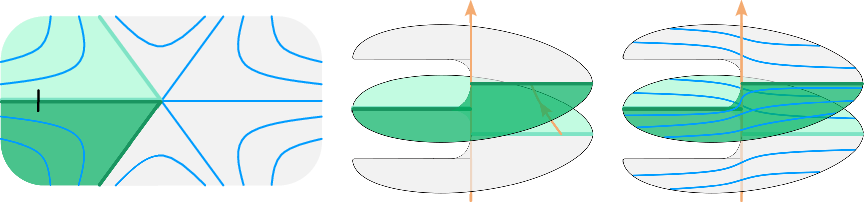}}
        \put(0,-1){$\ofoliation^s_\returnmap(x_-)$}
        \put(0,31.5){$\ofoliation^s_\returnmap(x_+)$}
        \put(-2.5,17){$x$}
        \put(80,-1){$\wh\Section_\partial$}
        \end{picture}
    \end{center}
    \caption{On the left: the leaf $\ofoliation^s_\returnmap(x_-)$ in dark green and $\ofoliation^s_\returnmap(x_+)$ in light green (green areas are visual indications). In the middle is represented their lifts to the Birkhoff section embedded in $M$. The boundary component is a negative boundary of $\Section$. One the right is represented in blue the trace of the weak stable foliation on the Birkhoff section.}
    \label{figure:LiftProjection}
\end{figure}
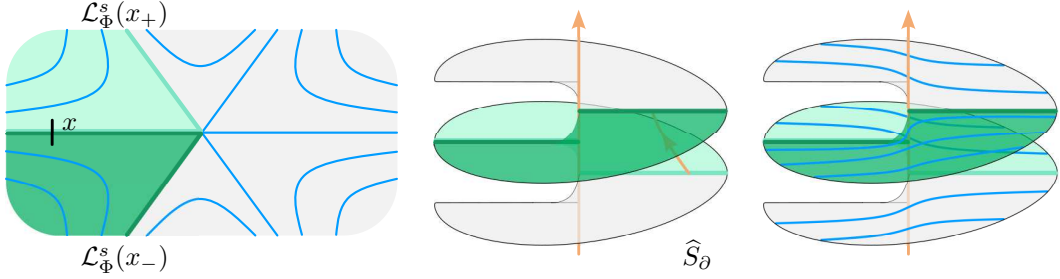

Consider now $x\in l^u_\returnmap$ so that $\ofoliation^s_\returnmap(x)\cap\Si\neq \emptyset$.  This intersection is a unique singular point $\sigma$, because there are no stable segments joining singular points.  We define the two following curves in $\wh\Section_\partial$, illustrated in Figure \ref{figure:LiftProjection}:
$$\ofoliation^s_\returnmap(x_-)=\lim_{y\to x, y<x, \ofoliation^s_\returnmap(y)\cap\Si=\emptyset} \ofoliation^s_\returnmap(y)$$
\centerline{and}
$$\ofoliation^s_\returnmap(x_+)=\lim_{y\to x, y>x, \ofoliation^s_\returnmap(y)\cap\Si=\emptyset} \ofoliation^s_\returnmap(y)$$

In the above definition, $y$ varies in $l^u_\returnmap$.
Both $\ofoliation^s_\returnmap(x_-)$ and $\ofoliation^s_\returnmap(x_+)$ consist in the union of $2$ separatrixes of $\sigma$: the one containing $x$ and another which is adjacent for the natural cyclic order of  the prongs at $\sigma$.  

Consider a small disc $D$ around $\sigma$  so that $D\cap \Si=\sigma$ and  $\partial D$ cuts each local prong in a unique point. Let $\alpha_+\colon [0,1]\to \partial D$ (resp. $\alpha_-\colon [0,1]\to\partial D$) be the arc of $\partial D$ joining the two prongs constituting $\ofoliation^s_\returnmap(x_+)$ (resp.$\ofoliation^s_\returnmap(x_-)$). 

Consider  the disjoint union of $l^u_\returnmap$, $\ofoliation^s_\returnmap(x_+)\setminus\{\sigma\}$ and $\alpha_+$.  Consider the quotient of the disjoint union by the following identification: 

 $$x\in l^u_\returnmap\simeq x\in \ofoliation^s_\returnmap(x_+),\quad
\alpha_+(0)\in \ofoliation^s_\returnmap(x_+) \simeq\alpha_+(0)\in \alpha_+$$ $$\mbox{and } 
\alpha_+(1)\in \ofoliation^s_\returnmap(x_+)\simeq \alpha_+(1)\in \alpha_+$$

One gets a simply connected set containing $P(O)$. Thus we can apply the lift/projection procedure with origin at $O$ leading to a parametrization

$$\theta^s_{x_+}\colon \ofoliation^s_\returnmap(x_+)\to \ofoliation^s(\theta(x)).$$ 
One easily check (by using the classical lift of homotopies on cover) that $\theta^s_{x+}$ does not depend on the choice of the disc $D$ around $\sigma$. 

We proceed in the same  way for $\ofoliation^s_\returnmap(x_-)$ by using $\alpha_-$ instead of $\alpha_+$ and the lift/projection procedure with origin at $O$ leads to a parametrization

$$\theta^s_{x_-}\colon \ofoliation^s_\returnmap(x_-)\to \ofoliation^s(\theta(x)).$$ 

\paragraph*{The crossing holonomies.}

We have seen that the lift/projection procedure with origin at $O$ applied to $\ofoliation^s_\returnmap(x_-)$ and $\ofoliation^s_\returnmap(x_+)$, $x\in \ofoliation^u_\returnmap$ provides two parametrization of the same leaf $\ofoliation^s(\theta(x))$ of $\orbitspace$, more precisely
$$\theta^s_{x_-}\colon \ofoliation^s_P(x_-)\to \ofoliation^s(\theta(x)), \quad \mbox{and }\quad \theta^s_{x_+}\colon \ofoliation^s_P(x_+)\to \ofoliation^s(\theta(x), $$

Notice that $\ofoliation^s_\returnmap(x_-)$ and $\ofoliation^s_\returnmap(x_+)$ each consists in the concatenations of $2$ successive prongs (for the cyclic  order) of $\sigma$: the one containing $x$, that we will denote by $L_x$ and another that we will denote simply $L_-\subset \ofoliation^s_\returnmap(x_-)$ and $L_+\subset \ofoliation^s_\returnmap(x_+)$.

\begin{definition} We denote by $h^s_x\colon \ofoliation^s_\returnmap(x_-)\to \ofoliation^s_\returnmap(x_+)$ the homeomorphism 
$$h^s_x:(\theta^s_{x_+})^{-1}\circ \theta^s_{x_-}$$
and we call it \emph{the crossing holonomy} of $x$. 
\end{definition}

The restriction of $\theta^s_{x_-}$ and of $\theta^s_{x,+}$ on $L_x$ is just the map $\theta_x$ obtained by the lift/projection procedure with origin at $O$ applied to $L_x$. Thus
$$h^s_x\equiv id \mbox{ on } L_x.$$

Now  $L_-$ and $L_+$ are mapped homeomorphically by $\theta^s_{x,-}$ and of $\theta^s_{x,+}$ on the same half leaf of $\ofoliation^s$.  Consider their lifts $\wt L_-$ and $\wt L_+$ on $\wt \Section\setminus \partial\wt \Section$, which is the first step of the lift/projection procedure. 

On $\wt M$ a leaf of $\wt \ofoliation^s_\returnmap$ cuts a $\wt \flow$-orbit in at most $1$ point. Thus one gets that for every point $\wt y\in \wt L_-$ there is exactly $1$ point $\wt z\in \wt L_+$ in the same $\wt \flow$ orbit: the point whose projection by $\Pi$ is the same as $\wt y$. Thus $\wt z$ is the image of $\wt y$ by a power of the first return map $\wt \returnmap$ of $\wt \flow$ on $\wt \Section$. In other words 
$$\exists k(x)\in\ZZ, \wt L_+ =\wt \returnmap^k(\wt L_-)$$
\centerline{and} 
$$ h^s_x\colon L_-\to L_+ \equiv \returnmap^{k(x)}\colon L_-\to L_+ .$$
In other words, $h^s_x\colon L_-\to L_+$ is the maps which sends $\sigma\in L_-$ on $\sigma\in L_+$ and which maps $y\in L_-$ on unique the point $z\in L_+$ so that 
$$ \mu^s([\sigma,z]^s)=\lambda^{k(x)}\mu^s([\sigma,y]^s),$$
(here the segments $[\sigma,z]^s$ and $[\sigma,y]^s$ are non oriented).

For calculating the crossing holonomies, it remains to calculate the integer $k$ above. It will depend on
\begin{itemize}
 \item the linking number and multiplicity $(\link(\sigma),\mult(\sigma))$ associated to $\sigma$,
 \item the period $p(\sigma)$ of $\sigma$ for $\returnmap$, or in other words, the number of boundary components of $\Section$ which are mapped on the periodic orbit $\gamma$ corresponding to $\sigma$,
 \item if $\sigma$ is on the positive or on the negative side of $x$ in its stable leaf. 
\end{itemize}

This calculation is the aim of next lemma which also summarize the arguments above. 

\begin{lemma} \label{l.crossing}

\begin{enumerate} \item Assume that the singular point $\sigma$ is at the right of $x$; in other words the stable segment $[x,\sigma]^s$ is positively oriented.

Then the crossing holonomy $h^s_x$ is
\begin{itemize}\item the identity map on the stable prong $L_x^s$ of $\sigma$ containing $x$, 
\item and is an affine map (for $\mu^s$) of ratio $\lambda^{k(x)}$  from the other prong $L_-\subset \ofoliation^s_{\returnmap}(x_-)$ to $L_+\subset \ofoliation^s_{\returnmap}(x_+)$ 
where 
$$k(x)=-\mult(\sigma)\cdot p(\sigma).$$ 
\end{itemize}
In particular, if $\mult(\sigma)>0$ then the crossing holonomy is a contraction on $L_-$ and if $\mult(\sigma)<0$ the crossing holonomy is an expansion. 

\item If  $\sigma$ is at the left of $x$ (i.e  $[\sigma,s]^s$ is positively oriented).

Then the crossing holonomy $h^s_x$ is
\begin{itemize}\item the identity map on the stable prong $L_x^s$ of $\sigma$ containing $x$, 
\item an affine map (for $\mu^s$) of ratio $\lambda^{k(x)}$  from  $L_-\subset \ofoliation^s_{\returnmap}(x_-)$ to $L_+\subset \ofoliation^s_{\returnmap}(x_+)$ 
where 
$$k(x)=\mult(\sigma)\cdot p(\sigma)$$ 
\end{itemize}
In particular, if $\mult(\sigma)>0$ then the crossing holonomy is a expansion on $L_-$ and if $\mult(\sigma)<0$ the crossing holonomy is a contraction. 

\end{enumerate}
\end{lemma}

\begin{proof} 
Consider   $L_-$ is a stable separatrix (prong) of $\sigma$.  Its $\returnmap$-orbit contains $|\link(\sigma)|$ stable prongs at $\sigma$ (and $|\link(\sigma)|p(\sigma)$ prongs on $\wh\Section_\partial$) . These prongs at~$\sigma$ are naturally cyclically ordered (the counterclockwise order).

Assume that $\sigma$ is at the right of $x$. Then 
the prong in $L_-$ follows $L$ which follows $L_+$ for the cyclic order on the orbit of $L_-$. We consider a simple segment $J$ joining a point $z$ in  $L_+$ to the point $y=(h^u_x)^{-1}(z)\in L_-$, so that~$J$ is disjoint from the orbit of the prong.

Let us consider a lift $\wt J$ on $\wt \Section\subset \wt M$,
The projection of $\wt J$ on $\orbitspace$ is a positively oriented circle around a point $\wt \gamma$ corresponding to $\sigma$.   This means that there is an orbit segment $\wt I$ of $\wt \flow$ joining $\wt J(1)$ to $\wt J(0)$. The concatenation of $\wt J$ with $\wt I$ is an embedded circle homotopic, in $\wt M\setminus \partial\wt \Section$ to a meridian of the $\wt \flow$ orbit $\wt \gamma$.  Its goes down on $M$ as a meridian of the periodic orbit $\gamma$ corresponding to $\sigma$.  The intersection number of a meridian with $\Section$  is $\mult(\sigma)p(\sigma)$ (here we fix the sign so that the intersection of  $\flow$ orbits with $\Section$ is positive). So the algebraic intersection $\wt\phi_{\epsilon}(\wt I)\algcap \wt \Section$ is equal to $\mult(\sigma)p(\sigma)$. Since $\wt I$ is an orbit ar from a lift of $y$ to a lift of $z$, this means that $$z= \wt \returnmap^{\mult(\sigma)p(\sigma)}(y).$$ 

On the stable prongs, $\wt\returnmap$ acts by multiplying $\mu^s$ by $\lambda^{-1}$, concluding.

When $\sigma$ is at the left of $x$ the proof follows the same way, but now the projection of the  path $\wt J$ on $\orbitspace$ is a negatively oriented meridian, explaining the change of sign. 
\end{proof}

By a useful  convention, we define $h^u_x$ by the identity map of $\ofoliation^s_\returnmap(x)$ when $\ofoliation^s_\returnmap(x)$ is a regular leaf of $\ofoliation^s_\returnmap$. 

\paragraph*{Generalized holonomies on $\Section$ and holonomies on $\orbitspace$.}

Let $I^u=[p,q]^u$ be an oriented unstable segment of $\ofoliation^u_\returnmap$. First assume that the stable leaves through $p$ and $q$ are disjoint from $\Si$.
We call \emph{generalized holonomy} $\cH^u_{p,q}\colon \ofoliation^s_\returnmap(p)\to \ofoliation^s_\returnmap(q)$ the map partially defined as follows:

Consider $x\in \ofoliation^s_\returnmap(p)$. Then $\cH^u_{p,q}(x)$ is defined and $\cH^u_{p,q}(x)=y\in \ofoliation^s_\returnmap(q)$ if and only if there is an ordered finite sequence 
$p=p_0,p_1,\dots ,p_k=q$ of points in $[p,q]^u$ and points $x_i\in \ofoliation^s_\returnmap(p_i)$, $x_0=x$,  and $x_k=y$ so that, for every $i\in \{0,\dots,k-1\}$, the (usual) holonomy $H^u_{p_i,p_{i+1}}(x_i)$ is defined and 
 $$h^u_{p_{i+1}}\circ H^u_{p_i,p_{i+1}}(x_i)=x_{i+1}.$$

In general the generalized holonomy is defined only on a sub-interval of $\ofoliation^s_\returnmap(p)$. In the case where the stable leaves through $p$ and $q$ are not assumed to be disjoint from $\Si$ we define in the same way the generalized holonomy
 $$\cH^u_{p,q}\colon \ofoliation^s_\returnmap(p_+)\to \ofoliation^s_\returnmap(q_+)$$

 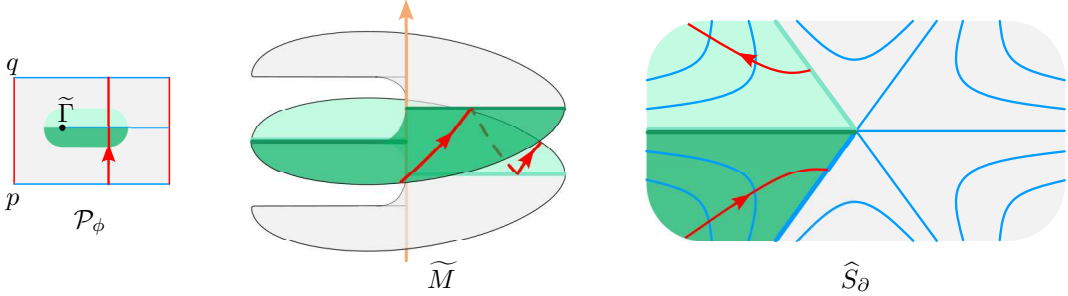
\begin{figure}
     \begin{center}
         \begin{picture}(120,35)(0,0) 
         \put(-10,0){\includegraphics[width=14cm]{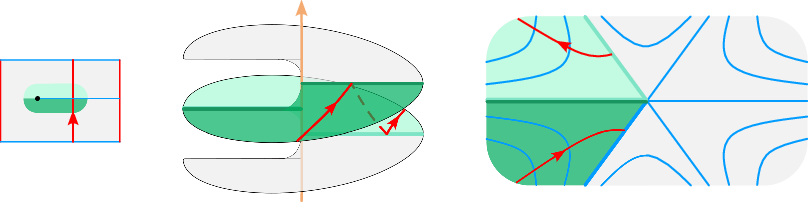}}
         \put(-4.1,19){$\wt\Gamma$}
         \put(-11,8){$p$}
         \put(-11,26){$q$}
         \put(-2,5){$\orbitspace$}
         \put(45,-3){$\wt M$}
         \put(100,-3){$\wh \Section_\partial$}
         \end{picture}
     \end{center}
     \caption{The holonomy in the orbit space on the left, inducing a generalized holonomy on the Birkhoff section in the middle and on the right. $\ofoliation^s_\returnmap(x_-)$ is in dark green and $\ofoliation^s_\returnmap(x_+)$ is in light green. The boundary component is a negative boundary of $\Section$.}
     \label{Cfigure:OrbitSpaceHolonomy}
 \end{figure}
 
\begin{remark} 
 \begin{enumerate}
\item  Note that the points $p_{i+1}\in [p,q]^u$ are all the points for which the segment $[p_{i+1}, H^u_{p_i,p_{i+1}}(x_{i})]^s$ may intersect (but not necessarily)  $\Si$. So we can keep the same family $p=p_0,p_1,\dots ,p_k=q$
 for all the points $y\in[p,x]^s$.  This leads to the next remark:  
 
 \item Let $I^u=[p,q]^u$ be an oriented unstable segment of $\ofoliation^u_\returnmap$.
  Consider $x\in  \ofoliation^s_\returnmap(p_+)$ and assume that the generalized holonomy $\cH^u_{p,q}$ is defined at $x$ and let $p=p_0,p_1,\dots ,p_k=q$ of points in $[p,q]^u$ the points used in the definition of $\cH^u_{p,q}$.  Then  the restriction of $\cH^u_{p,q}$ to  the (non-oriented) stable segment $[p,x]^s$ is:
  $$\cH^u_{p,q}|_{[p,x]^s}= \Pi_{i=0}^{k-1}h^u_{p_{i+1}}\circ H^u_{p_i,p_{i+1}} $$
  \end{enumerate}
 \end{remark}

 Next statement, whose proof is now almost straightforward, is one of the main objectives of this section: 
 
 \begin{theorem}\label{t.holonomies} Let $[p^0,q^0]^u\subset \orbitspace$ be an unstable segment disjoint from $\Si$  and let $[p,q]^u=\theta^{-1}([p_0,q_0]^u)\subset \Section $ be a corresponding unstable segment by the lift/projection procedure. We denote 
 $$\theta_p\colon \ofoliation^s_\returnmap(p_+)\to \ofoliation^s(p_0) \mbox{ and }\theta_q\colon \ofoliation^s_\returnmap(q_+)\to \ofoliation^s(q_0)$$ the parametrization given by the lift/projection procedure. 
 
 Then 
 $$H^u_{p^0,q^0}=\theta_{q_+}\circ \cH^u_{p,q}\circ \theta_{p_+}^{-1}$$
  
 \end{theorem}

 \begin{proof}
    Consider a point $x^0\in \ofoliation^s(p^0)$ for which $H^u_{p^0,q^0}(x^0)$ is defined and let $x\in \ofoliation^s_{\returnmap}(p_+)=\theta_{p_+}^{-1}(x)$. 
 
    Then, on $[p_0,x_0]^s$ we have: 
    $$
    \begin{array}{cl}
    H^u_{p_0,q_0}&=\theta_{q_+}\circ \cH^u_{p,q}\circ \theta_{p_+}^{-1}\\
    &=\theta_{q_+}\circ \left(\Pi_{i=0}^{k-1}h^u_{p_{i+1}}\circ H^u_{p_i,p_{i+1}}\right)\circ \theta_{p_+}^{-1}\\
    &= \Pi_{i=0}^{k-1}  \theta_{(p_{i+1})_+} \circ h^u_{p_{i+1}}\circ H^u_{p_i,p_{i+1}}\circ \theta_{(p_i)_+}^{-1}\\
    &=\Pi_{i=0}^{k-1}  \theta_{(p_{i+1})_-} \circ H^u_{p_i,p_{i+1}}\circ \theta_{(p_i)_+}^{-1} \mbox{, by definition of } h^u_{p_i+1}.\\
    \end{array}
    $$
    
    In this composition, each term $\theta_{(p_{i+1})_-} \circ H^u_{p_i,p_{i+1}}\circ \theta_{(p_i)_+}^{-1}$ is
    $$\theta_{(p_{i+1})_-} \circ H^u_{p_i,p_{i+1}}\circ \theta_{(p_i)_+}^{-1}=H^u_{p^0_i,p^0_{i+1}} \mbox{, where } p^0_i= \theta(p_i)$$
    
    In particular, each term is a holonomy map for $\ofoliation^u$, so the composition is a holonomy map for $\ofoliation^u$, concluding. 
 \end{proof}

In other words, one can read the holonomies of $\ofoliation^u$ between leaves of $\ofoliation^s$ as follows: we look at the leaves on the surface  $\Section$. We flow the segments $[p,x]^s$ isometrically along the leaves of $\ofoliation^u_P$ until it crosses $\Si$ (if it does not cross $\Si$ one just gets the usual holonomy). Then one apply the crossing holonomy and the $\mu^s$-length of the segments is changed according to Lemma~\ref{l.crossing}; then one proceeds in the same way. 

Next lemma answers to the natural question: how this process stops?: 

\begin{lemma}\label{l.infinite} The generalized holonomy $\cH^u_{p,q}$ is not defined at $x\in \ofoliation^s_P(p_+)$ if and only if there is an infinite increasing sequence $p_0,\cdots p_i,\dots,\in [p,q]^u$ so that $\cH^u_{p,p_i}(x)$ is defined and 
$$\lim \mu^s\left([p_i,\cH^u_{p,p_i}(x)]^s\right)=+\infty$$
\end{lemma}
\begin{proof}Corollary~\ref{c.bounded} implies that, given $\ell$ there is $\delta(\ell)>0$ so that if we have $\mu^s\left([p_i,\cH^u_{p,p_i}(x)]^s\right) \leq \ell$ then $\mu^u([p_i,p_{i+1}]^u)>\delta$. Thus if $\mu^s\left([p_i,\cH^u_{p,p_i}(x)]^s\right)$ remains bounded there are finitely many non-trivial crossing holonomies and the generalized holonomy $\cH^u_{p,q}(x)$ is well-defined. 
\end{proof}

\subsection{Positive Birkhoff section and $\RR$-covered flows}

We give the third proof of the implication $\Rightarrow$ in Theorem \ref{theorem:rcoveredcondition}. The main step for the proof is Proposition~\ref{p.complete} below:

\begin{proposition}\label{p.complete} If $\flow$ admits a  Birkhoff section so that every linking number  $\mult(\sigma)$ is positive, then every $C_{+,+}$ quadrant is complete. 
 
\end{proposition}

\begin{proof}According to Lemma~\ref{l.crossing} every crossing holonomy in the positive side is contracting  so that  

$$\mu^s(\cH^u_{p,t}([p,x]^s)) \leq \mu^s([p,x]^s) $$ 
for $p,t,x\in \Section$ so that $[p,t]^u$ and $[p,x]^s$ positively oriented.  Thus Lemma~\ref{l.infinite} implies that, for every $x\in \ofoliation^s_{\returnmap,+}(p)$  the generalized holonomy $\cH^u_{p,q}(x)$ is defined for every $q\in \ofoliation^u_{\returnmap,+}(p)$. 

Now Theorem~\ref{t.holonomies} implies that for every $p_0$ in $\orbitspace$ so that $\ofoliation^u(x)\cap \Si=\emptyset$ one has that the holonomy $H^u_{p,q}(x)$ is well-defined for every $q\in \ofoliation^u_+(p)$ and $x\in \ofoliation^s_+(x)$.  In other words, $C_{+,+}(p)$ is complete.  As the hypothesis holds for $p$ in a dense open subset of $\orbitspace$ and as the completeness is a closed property, one gets that every quadrant $C_{+,+}(x)$ is complete, ending the proof of the proposition.
\end{proof}

\begin{proof}[Proof of Theorem~\ref{theorem:rcoveredcondition}]
 Proposition \ref{p.complete} proves that the flow $\flow$ is $\RR$-covered, either suspension or positively twisted. 

The periodic orbits of a suspension flow all correspond to a positive homology class in $H_1(M,\ZZ)=\ZZ$.  If $\Section$ is positive or negative, $\partial \Section$ consists in homology classes all of the same sign, but the sum vanishes: this implies that $\Section$ is indeed a transverse torus. 

Otherwise, $\flow$ is positively twisted, ending the proof of the theorem in the case where the stable and unstable bundle of  $\flow$ are orientable. In the non-orientable case, one considers the lift $\flow^{Or}$ of $\flow$ to the $2$-folds orientation cover~$M^{Or}$ of its stable (and unstable) bundle. The positive (resp. negative) Birkhoff section~$\Section$ lifts on  $M^{Or}$ as a positive (resp. negative) Birkhoff section of $\flow^{Or}$. Therefore~$\flow^{Or}$ is a positively (resp. negatively) twisted $\RR$-covered Anosov flows, and this property goes down to $\flow$, ending the proof. 
\end{proof}